\documentclass{amsart}

\usepackage{geometry}
\usepackage{srcltx,graphicx}
\usepackage{amsmath,amssymb,amsthm,amscd,mathrsfs,bm}
\usepackage[unicode,linkcolor=blue,hyperindex]{hyperref}
\usepackage{indentfirst} 
\usepackage{enumerate} 
\usepackage{cases} 
\usepackage[mathscr]{euscript}
\usepackage{url,breakurl}
\usepackage{algorithm}
\usepackage{algpseudocode}
\usepackage{tikz-cd}
\usepackage{stmaryrd}
\usepackage{tcolorbox}
\usepackage{arydshln}
\usepackage{listings}
\usepackage{subfig}
\usepackage{colortbl}
\usepackage{leftidx}
\usepackage{array,multirow}
\usepackage{multicol}
\usepackage{epstopdf}
\usepackage{undertilde}

\numberwithin{equation}{section} \numberwithin{figure}{section}
{\theoremstyle{remark} \newtheorem{remark}{Remark}}

\newtheorem{proposition}{Proposition}[section]

\newtheorem{theorem}{Theorem}[section]
\newtheorem{corollary}{Corollary}[section]
\newtheorem{lemma}{Lemma}[section]

\newcommand\Laps{(-\Delta)^s}

\newcommand{\cT}{\mathcal{T}}
\newcommand{\cN}{\mathcal{N}}
\newcommand{\cB}{\mathcal{B}}
\newcommand{\cR}{\mathcal{R}}
\newcommand{\cJ}{\mathcal{J}}
\newcommand{\cI}{\mathcal{I}}
\newcommand{\cP}{\mathcal{P}}

\newcommand{\cV}{\mathcal{V}}

\newcommand{\polR}{\mathbb{R}}
\newcommand{\polRd}{{\mathbb{R}^d}}
\newcommand{\vA}{{\textbf A}}
\newcommand{\vB}{{\textbf B}}
\newcommand{\vI}{{\textbf I}}

\newcommand{\wt}{\widetilde}
\newcommand{\wh}{\widehat}
\newcommand{\supp}{\mathrm{supp }}
\newcommand{\cond}{\mathrm{cond }}

\newcommand{\eps}{\varepsilon}

\begin{document}
\title[{Robust BPX Preconditioner for Fractional
Laplacians}]{Robust BPX Preconditioner for Fractional
Laplacians on Bounded Lipschitz Domains}

\author[J.P.~Borthagaray]{Juan Pablo~Borthagaray}
\address[J.P.~Borthagaray]{Departamento de Matem\'atica y
Estad\'istica del Litoral, Universidad de la Rep\'ublica, Salto,
Uruguay}
\email{jpborthagaray@unorte.edu.uy}
\thanks{JPB has been supported in part by NSF grant DMS-1411808 and
Fondo Vaz Ferreira grant 2019-068}

\author[R.H.~Nochetto]{Ricardo H.~Nochetto}
\address[R.H.~Nochetto]{Department of Mathematics and Institute for
Physical Science and Technology, University of Maryland, College
Park, MD 20742, USA}
\email{rhn@math.umd.edu}
\thanks{RHN has been supported in part by NSF grant DMS-1411808}

\author[S.~Wu]{Shuonan~Wu}
\address[S.~Wu]{School of Mathematical Sciences, Peking University,
  Beijing, 100871, P.R. China}
\email{snwu@math.pku.edu.cn}
\thanks{SW has been supported in part by the startup grant from Peking
University and the National Natural Science Foundation of China grant
No. 11901016.}

\author[J.~Xu]{Jinchao~Xu}
\address[J.~Xu]{Department of Mathematics, Pennsylvania State
University, University Park, PA 16802, USA}
\email{xu@math.psu.edu}
\thanks{JX has been supported in part by NSF grant DMS-1819157}
%\date{}

\begin{abstract}
We propose and analyze a robust BPX preconditioner for 
the integral fractional Laplacian on bounded Lipschitz domains.
For either quasi-uniform grids or graded bisection grids, we show that
the condition numbers of the resulting systems remain uniformly
bounded with respect to both the number of levels and the fractional
power. The results apply also to the spectral and censored
fractional Laplacians.
\end{abstract}

\maketitle

%------------------------------------------------------------------------
%------------------------------------------------------------------------

\section{Introduction}
%------------------------------------------------------------------------
%------------------------------------------------------------------------
Given $s \in (0,1)$, the fractional Laplacian of order $s$ in $\polRd$
is the pseudodifferential operator with symbol $|\xi|^{2s}$. That is,
denoting the Fourier transform by $\mathcal{F}$, for every function
$v:\polRd \to \polR$ in the Schwartz class $\mathcal{S}$ it holds that
\begin{equation} \label{eq:fourier}
\mathcal{F}\left(\Laps v\right) (\xi) = |\xi|^{2s} \mathcal{F}(v) (\xi).
\end{equation}
Upon inverting the Fourier transform, one obtains the following
equivalent expression:
\begin{equation}
\label{eq:defofLaps}
\Laps v(x) = C(d,s) \mbox{ p.v.} \int_\polRd
\frac{v(x)-v(y)}{|x-y|^{d+2s}} d y, \qquad C(d,s) = \frac{2^{2s} s
\Gamma(s+\frac{d}{2})}{\pi^{d/2}\Gamma(1-s)}.
\end{equation}
The constant $C(d,s) \simeq s(1-s)$ compensates the singular behavior
of the integrals for $s\to0$ (as $|y| \to \infty$) and for $s\to1$ (as
$y\to x$), and yields \cite[Proposition 4.4]{Nezza_BSM2012}
\begin{equation}\label{eq:limits_laps}
\lim_{s \to 0} \Laps v(x) = v(x), \quad \lim_{s \to 1} \Laps v(x) =
-\Delta v(x), \quad \forall v \in C^\infty_0(\polRd).
\end{equation}

From a probabilistic point of view, the fractional Laplacian is
related to a simple random walk with arbitrarily long jumps
\cite{Valdinoci09}, and is the infinitesimal generator of a 
$2s$-stable process \cite{Bertoin96}.  Thus, the fractional Laplacian
has been widely utilized to model jump processes arising in social and
physical environments, such as finance \cite{CoTa04}, predator search
patterns \cite{Sims08}, or ground-water solute transport
\cite{BeWhMe00}.

There exist several nonequivalent definitions of a
fractional Laplace operator $\Laps$ on a bounded domain $\Omega
\subset \polRd$ (see \cite{Bonforte18, BBNOS18}). Our emphasis is on the
homogeneous Dirichlet problem for the
{\em integral} Laplacian:  given $f \colon \Omega
\to \polR$, one seeks $u \colon \polRd \to \polR$ such that
\begin{equation} \label{eq:Dirichlet}
\left\lbrace \begin{array}{rl}
(-\Delta)^s u  = f & \mbox{in } \Omega, \\
u = 0 & \mbox{in } \Omega^c,
\end{array} \right.
\end{equation}
where the pointwise definition of $\Laps u(x)$ is given by \eqref{eq:defofLaps}
for $x\in\Omega$. Consequently, the integral fractional Laplacian on $\Omega$
maintains the probabilistic interpretation and corresponds to a {\em killed}
L\'evy process \cite{Bertoin96, ChenSong05}.  It is noteworthy that,
as the underlying stochastic process admits jumps of arbitrary length,
volume constraints for this operator need to be defined in the
complement of the domain $\Omega$. 

Weak solutions to \eqref{eq:Dirichlet} are the minima of the
functional $v \mapsto \frac12 |v|^2_{H^s(\polRd)} - \int_\Omega f v$
on the zero-extension space $\widetilde{H}^s(\Omega)$ (see Section
\ref{sec:variational-formulation}).  In accordance with \eqref{eq:limits_laps}
restricted to any $v\in C_0^\infty(\Omega)$, it holds that if $v \in
\widetilde{H}^\sigma(\Omega)$ for some $\sigma > 0$, then
\cite{MaSh02}
\begin{equation} \label{eq:limit_norms_0}
\lim_{s \to 0^+} |v|_{H^s(\polRd)} = \| v \|_{L^2(\Omega)},
\end{equation}
while if $v\in L^2(\polRd)$ is such that $\mbox{supp } v \subset
\overline\Omega$ and $\lim_{s \to 1^-} |v|_{H^s(\polRd)}$ exists and
is finite, then $v \in H^1_0(\Omega)$ and \cite{BoBrMi01}
\begin{equation} \label{eq:limit_norms_1}
\lim_{s \to 1^-} |v|_{H^s(\polRd)} = | v |_{H^1(\Omega)}.
\end{equation}

Consider a discretization of \eqref{eq:Dirichlet} using standard
linear Lagrangian finite elements (see details in Section
\ref{sec:discrete_setting}) on a mesh $\mathcal{T}$ whose elements
have maximum and minimum size $h_{\max}$ and $h_{\min}$ respectively,
and denote by $\vA$ the corresponding stiffness matrix. Then, as shown
in \cite{AiMcTr99}, the condition number of $\vA$ obeys the relation
\begin{equation} \label{eq:cond_stiffness}
  \cond \, (\vA) \simeq ({\rm dim} \mathbb{V}(\mathcal{T}))^{2s/d} \left(
  \frac{h_{\max}}{h_{\min}} \right)^{d-2s}
\end{equation}
for $0 < s < 1$ with $2s \le d$, and one can remove the factor
involving $\frac{h_{\max}}{h_{\min}}$ by preconditioning ${\bf A}$ by
a diagonal scaling. On non quasi-uniform grids, the hidden constant in
the critical case $2s=d$ is worse by a logarithmic factor.

In recent years, efficient finite element discretizations of
\eqref{eq:Dirichlet} have been examined in several papers. Adaptive
algorithms have been considered in
\cite{du2013convergent,ainsworth2017aspects,gimperlein2019space},
and a posteriori error analysis has been addressed in
\cite{nochetto2010posteriori, faustmann2019quasi}. Standard finite
element discretizations of the fractional Laplacian give rise to full
stiffness matrices; matrix compression techniques have been proposed
and studied in
\cite{zhao2017adaptive,ainsworth2018towards,karkulik2018mathcal}.  For
the efficient resolution of the discrete problems, operator
preconditioners have been considered in \cite{gimperlein2019optimal}.

%In this work, we propose a multilevel BPX preconditioner (cf.
%\cite{xu1989thesis, bramble1990parallel}) $\vB$
%for the solution of \eqref{eq:Dirichlet}, that yields $\cond(\vB \vA)
%\lesssim 1$. In general, our result follows from the general theory
%for multigrid preconditioners (cf. \cite{xu1992iterative,
%griebel1995abstract, xu1996auxiliary, xu2002method}), although the
%nonlocal nature of the problem forces us to make modifications on some
%steps. For example, \rhn{the strengthened Cauchy-Schwarz inequality
%encodes the separation of scales and follows from the Fourier representation
%\eqref{eq:fourier} of the fractional Laplacian \eqref{eq:defofLaps} rather than
%elementwise integration by parts 
%(cf. for example, \cite[Lemma 6.1]{xu1992iterative} and
%\cite[Lemma (7.2.16)]{BrSc07}).}

In this work, we propose a multilevel BPX preconditioner (cf.
\cite{xu1989thesis, bramble1990parallel}) $\vB$ for the solution of
\eqref{eq:Dirichlet}, that yields $\cond(\vB \vA) \lesssim 1$. In
general, our result follows from the general theory for multigrid
preconditioners (cf. \cite{xu1992iterative, griebel1995abstract,
xu1996auxiliary, xu2002method}).  An important consequence of
\eqref{eq:limit_norms_0} and \eqref{eq:limit_norms_1} is that, on any
given grid, the stiffness matrices associated with integral fractional
Laplacians of order $s$ approach either the standard mass matrix (as
$s \to 0$) or the stiffness matrix corresponding to the Laplacian (as
$s \to 1$), the latter because the canonical basis functions of
$\mathbb{V}(\mathcal{T})$ are Lipschitz and $W^{1,\infty}_0(\Omega)
\subset \wt{H}^s(\Omega)$. This is consistent with
\eqref{eq:cond_stiffness}: for example, on quasi-uniform grids of size
$h$, such a formula yields $\cond (\vA) \simeq h^{-2s}$. 

Based on the above observations, one of our main goals is to
obtain a preconditioner that is {\it uniform with respect to $s$ as
well as with respect to the number of levels $\bar{J}$.} For such a
purpose, we need to weight the contributions of the coarser levels
differently to the finest level.  On a family of quasi-uniform grids
$\{\overline{\cT}_k\}_{k=0}^{\bar{J}}$ with size $\bar{h}_k$, we shall
consider a preconditioner in the operator form (cf.
\eqref{eq:BPX-FPDE} below)
\begin{equation}\label{eq:bpx-intro}
  \overline{B} = \overline{I}_{\bar{J}}
  \bar{h}_{\bar{J}}^{2s} \overline{Q}_{\bar{J}} +
  (1-\wt{\gamma}^{s})\sum_{k=0}^{\bar{J} - 1} \overline{I}_k \bar{h}_k^{2s}
  \overline{Q}_k,
\end{equation}
where the arbitrary parameter $\wt{\gamma}\in (0,1)$. Above,
$\overline{Q}_k$ and $\overline{I}_k$ are suitable $L^2$-projection
and inclusion operators, respectively.  Clearly, if $s \in (0,1)$ is
fixed, then the factor $1-\wt{\gamma}^{s}$ is equivalent to a
constant.  However, such a factor tends to $0$ as $s \to 0$, and this
correction is fundamental for the resulting condition number to be
uniformly bounded with respect to $s$.

We now present a simple numerical example to illustrate this point.
Let $\Omega = (-1,1)^2$, $f = 1$, $s = 10^{-1}, 10^{-2}$, and choose
either $\wt{\gamma} = 0$ (i.e., no correction) and $\wt{\gamma} =
\frac12$ in the preconditioner above to compute finite element
solutions to \eqref{eq:Dirichlet} on a sequence of nested grids. The
left panel in Table \ref{tab:comparison} shows the number of
iterations needed to solve the resulting linear system by using a
Preconditioned Conjugate Gradient (PCG) method with a fixed tolerance.
It is apparent that setting $\wt{\gamma} = \frac12$ gives rise to a
more robust behavior with respect to either $s$ and the number of
levels $\bar{J}$.
\begin{table}[!ht]
\centering
\begin{tabular}{|c||c|c||c|c|}
  \hline
  \multicolumn{5}{|c|}{Uniform grids} \\ \hline
 \multirow{2}{*}{DOFs} 
  & \multicolumn{2}{c||}{$s=10^{-1}$}
  & \multicolumn{2}{c|}{$s=10^{-2}$} \\ \cline{2-5}
  &  $\tilde{\gamma} = 0$ & $\tilde{\gamma} = \frac12$  &
  $\tilde{\gamma} = 0$  & $\tilde{\gamma} = \frac12$ \\ \hline
 225 & 14 & 10 & 16 & 10 \\ 
961 & 17 & 10 & 18 & 10 \\ 
3969 & 19 & 10 & 21 & 10 \\ 
16129 & 20 & 10 & 23 & 9 \\ 
\hline  
\end{tabular}
\qquad
\begin{tabular}{|c||c|c||c|c|}
  \hline
  \multicolumn{5}{|c|}{Graded bisection grids} \\ \hline
 \multirow{2}{*}{DOFs} 
  & \multicolumn{2}{c||}{$s=10^{-1}$}
  & \multicolumn{2}{c|}{$s=10^{-2}$} \\ \cline{2-5}
  &  $\tilde{\gamma} = 0$ & $\tilde{\gamma} = \frac12$  &
  $\tilde{\gamma} = 0$  & $\tilde{\gamma} = \frac12$ \\ \hline
161 & 13 & 10 & 15 & 11 \\ 
853 & 17 & 12 & 19 & 13 \\ 
2265 & 20 & 12 & 22 & 14 \\ 
9397 & 22 & 13 & 25 & 14 \\ 
\hline  
\end{tabular}
\medskip
\caption{Number of iterations needed when using a PCG method with BPX
  preconditioner without ($\wt{\gamma} = 0$) and with
  ($\wt{\gamma} = \frac12$) a correction factor. We display results
  on a family of uniformly refined grids (left panel), and on a
  sequence of suitably graded bisection grids (right panel).}
\label{tab:comparison}
\end{table}

Another aspect to take into account in our problem is the low
regularity of solutions. As we discuss in Section
\ref{sec:discrete_setting}, by using uniform grids one can only expect
convergence in the energy norm with order $\mathcal O \left({\rm dim }
\mathbb{V}(\mathcal{T})^{-1/(2d)} \right)$ (up to logarithmic factors)
independently of the smoothness of the data. The reason for such a low
regularity of solutions is boundary behavior; exploiting the a priori
knowledge of this behavior by means of suitably refined grids leads to
convergence with order $\mathcal O \left({\rm dim }
\mathbb{V}(\mathcal{T})^{-1/(2d-2)} \right)$ if $d \ge 2$ and
$\mathcal O \left({\rm dim} \mathbb{V}(\mathcal{T})^{s-2} \right)$ if
$d=1$.  In spite of this advantage, graded grids give rise to
worse-conditioned matrices, as described by \eqref{eq:cond_stiffness}.
This work also addresses preconditioning on graded bisection
grids, that can be employed to obtain the refinement as needed.
Our algorithm on graded bisection grids builds on the
subspace decomposition introduced in \cite{chen2012optimal}, that
leads to optimal multilevel methods for classical ($s=1$) problems. 
Our theory on graded bisection grids, however, differs from
\cite{chen2012optimal} to account for the uniformity with respect to
$s$. As illustrated by the right panel in Table \ref{tab:comparison},
including a correction factor on the coarser scales leads to a more
robust preconditioner.

The integral (or restricted) Laplacian operator
\eqref{eq:defofLaps} turns out to be {\it spectrally equivalent} to
the spectral Laplacian uniformly with respect to $s$ on bounded
Lipschitz domains \cite{ChenSong05}.  This property is also a
consequence of our multilevel space decomposition.  Such a
spectral equivalence can be extended to the censored (or
regional) Laplacian for $s \in (\frac12, 1)$,
%except for $s = \frac{1}{2}$ \cite{grisvard2011elliptic,ChenKimSong10}, 
and the uniformity of equivalence constant holds when $s \to 1$.
However, the three operators have a strikingly different boundary
behavior \cite{Bonforte18,Grubb15,ROSe14}.  We present their
definitions along with their properties in Section
\ref{sec:spectral-equivalence}.  Consequently, the BPX preconditioner
\eqref{eq:bpx-intro} for the integral Laplacian on quasi-uniform grids
and its counterpart on graded bisection grids apply as well to the
spectral and censored Laplacians except for the censored one when
$s \to \frac{1}{2}$.  Further, the uniformity with respect to $s$
holds for both the integral or spectral Laplacians.  

This paper is organized as follows. Section \ref{sec:preliminaries}
collects preliminary material about problem \eqref{eq:Dirichlet}, in
particular regarding its variational formulation, regularity of
solutions and its approximation by the finite element method. Next, in
Section \ref{sec:Robust-BPX} we review some additional tools that we
need to develop the theory of a robust BPX preconditioner; we discuss
general aspects of the method of  subspace corrections and introduce
an $s$-uniform decomposition that plays a central role in our
analysis.  We introduce a BPX preconditioner for quasi-uniform grids
in Section \ref{sec:BPX-uniform}, and prove that it leads to condition
numbers uniformly bounded with respect to the number of refinements
$\bar{J}$ and the fractional power $s$.  Afterwards, we delve into the
preconditioning of systems arising from graded bisection grids. For
that purpose, Section \ref{S:bisection} offers a review of the
bisection method with novel twists, while Section \ref{sec:BPX-bisect}
proposes and studies a BPX preconditioner on graded bisection grids.
Section \ref{sec:experiments} presents some numerical experiments that
illustrate the performance of the BPX preconditioners.  The paper
concludes with three Appendices that collect and prove a few technical
results.

%------------------------------------------------------------------------
%------------------------------------------------------------------------
\section{Preliminaries} \label{sec:preliminaries}
%----------------------------------------------------------------------
%----------------------------------------------------------------------
In this section we set the notation used in the rest of the paper
regarding Sobolev spaces and recall some preliminary results about
their interpolation. We are particularly concerned with the
zero-extension Sobolev space $\widetilde{H}^\sigma(\Omega):=
\overline{C^\infty_0(\Omega)}^{\|\cdot\|_{H^\sigma(\polRd)}}$, which
is the set of functions in $H^\sigma(\mathbb{R}^d)$ whose support is
contained in $\Omega$.  Given $u,v \in \widetilde H^\sigma (\Omega)$,
we define below the (scaled) inner product $(u, v)_\sigma = (u,
v)_{H^\sigma(\polR^d)}$ in $\widetilde{H}^\sigma(\Omega)$, the
corresponding norm $|u|_\sigma = ( u, u )_\sigma^{\frac12} =
|u|_{H^\sigma(\polR^d)}$, and let $\|u\|_0 := \|u\|_{L^2(\Omega)}$.
Moreover, we discuss regularity of solutions to \eqref{eq:Dirichlet}
and a priori error estimates for finite element approximations.

For convenience, we write $X \lesssim Y$ (resp. $X \gtrsim Y$) to
indicate $X \leq C Y$ (resp. $CX \geq Y$), where $C$ denotes, if not
specified, a generic positive constant that may stand for different
values at its different occurrences but is independent of the number
of levels or fractional power. The notation $X \simeq Y$ means both $X
\lesssim Y$ and $X \gtrsim Y$ hold.

%----------------------------------------------------------------------
\subsection{Variational formulation}
\label{sec:variational-formulation}
%----------------------------------------------------------------------

The natural setting to study the variational formulation of fractional
diffusion problems such as \eqref{eq:Dirichlet} is in Sobolev spaces
$\widetilde{H}^\sigma(\Omega)$ of non-integer order $\sigma$.  We
refer to \cite{BoLeNo20} for basic definitions and the notation we use
here. We consider the symmetric bilinear form $(\cdot, \cdot)_\sigma
\colon \widetilde H^\sigma(\Omega) \times \widetilde H^\sigma(\Omega)
\to \polR$,
\begin{equation} \label{eq:bilinear}
( u , v )_\sigma := \frac{C(d,\sigma)}{2}
\iint_{\polR^d\times\polR^d}
\frac{(u(x)-u(y))(v(x)-v(y))}{|x-y|^{d+2\sigma}} \, dx \, dy, 
\end{equation}
where $C(d,\sigma)$ is the constant from \eqref{eq:defofLaps}. We
point out that, because functions in $\widetilde{H}^\sigma(\Omega)$
vanish in $\Omega^c$, the integration takes place in $\left(\Omega
\times \polR^d\right) \cup \left(\polR^d \times \Omega \right)$.

Since a Poincar\'e inequality is valid in $\widetilde{H}^s(\Omega)$
(cf. \cite[Prop. 2.4]{AcBo17}, for example), the map $u \mapsto
(u,u)_s$ is an inner product $\wt{H}^s(\Omega)$.  Given $f \in
H^{-s}(\Omega)$, the dual of $\wt{H}^s(\Omega)$, the weak formulation
of the homogeneous Dirichlet problem \eqref{eq:Dirichlet} reads: find
$u \in \widetilde H^s(\Omega)$ such that
\begin{equation} \label{eq:weak}
a(u,v) := (u,v)_s = \langle f, v \rangle_{s,\Omega}  \quad \forall v
\in \widetilde H^s(\Omega),
\end{equation}
where $\langle \cdot, \cdot \rangle_{s,\Omega}$ stands for the duality
pairing between $H^{-s}(\Omega)$ and $\widetilde H^s(\Omega)$.
Existence and uniqueness of solutions of \eqref{eq:weak} is a consequence
of the Riesz representation theorem.

%----------------------------------------------------------------------
\subsection{Interpolation and fractional Sobolev spaces}
\label{sec:interpolation}
%----------------------------------------------------------------------
An important feature of the fractional Sobolev scale is that it can be
equivalently defined by {\it interpolation} of integer-order spaces.
This along with the observation that the norm equivalence constants are
uniform with respect to $s$ is fundamental for our work.  In view of
applications below, we now recall the abstract setting for two Hilbert
spaces $X^1 \subset X^0$ with $X^1$ continuously embedded and dense in
$X^0$.  Following \cite[Section 2.1]{lions1972hilbert}, the inner
product in $X^1$ can be represented by a self-adjoint and coercive
operator $S:D(S)\to X^0$ with domain $D(S)\subset X^1$ dense in $X_0$,
i.e. $(v,w)_{X^1} = (Sv,w)_{X^0}$ for all $v\in D(S), w\in X^1$.
Invoking the spectral decomposition of self-adjoint operators
\cite{yosida1965functional}, we let $\Lambda : X^1\to X^0$ be the
square root of $S$, which in turn is self-adjoint, coercive, and
satisfies
\begin{equation} \label{eq:operator-Lambda}
  (v,w)_{X^1} = (\Lambda v,\Lambda w)_{X^0}
  \quad \forall v,w \in X^1.
\end{equation}
Suppose further that the spectrum $\{\lambda_k\}_{k=1}^\infty$ of
$\Lambda$ is discrete and the corresponding eigenvectors
$\{\varphi_k\}_{k=1}^\infty$ form a complete orthonormal basis for
$X^0$; hence $\Lambda v = \sum_{k=1}^\infty \lambda_k v_k\varphi_k$
for all $v = \sum_{k=1}^\infty v_k\varphi_k\in X^1$. Then,  we can
define a fractional power $s \in (0,1)$ of $\Lambda$ as follows:
\begin{equation} \label{eq:Lambda-s}
  \Lambda^s v := \sum_{k=1}^\infty \lambda_k^s v_k\varphi_k \quad \text{if}
  \quad \|\Lambda^s v\|_{X_0}^2 :=  \sum_{k=1}^{\infty} \lambda_k^{2s} v_k^2
  < \infty.
\end{equation}

On the other hand, following \cite[Appendix B]{mclean2000strongly}, we
consider a variant of the classical $K$-method which, for simplicity,
we write for $L^2$-based interpolation. Decompose $v = v^0 + v^1$ with
$v^0\in X^0, v^1\in X^1$, take $t>0$, and set
\begin{equation} \label{eq:modified-K}
K_2(t,v) := \inf_{\substack{v^0 \in X^0, v^1 \in X^1 \\ v = v^0 +
  v^1}}  \big( \|v^0 \|_{X^0}^2 + t^2 \| v^1 \|_{X^1}^2
\big)^{\frac12}.
\end{equation}
It is immediate to verify that $K_2$ is equivalent to the usual
$K$-functional:
\[
K_2(t,v) \le K(t,v) \le \sqrt{2} K_2(t,v) \quad \forall v, t.
\]
Given $s \in (0,1)$, we consider the interpolation space $(X^0,X^1)_{s,2}$
with norm
\begin{equation} \label{eq:interpolation_norms}
\|v\|_{(X^0,X^1)_{s,2}} := \left( \frac{2 \sin(\pi s)}{\pi}
  \int_0^\infty t^{-1-2s} K_2(v,t)^2 dt \right)^{\frac12}. 
\end{equation}

The following theorem gives an intrinsic spectral equivalence between
the interpolation by $K$-method and spectral theory; see \cite[Theorem
  15.1]{lions1972hilbert} for a more general statement.

\begin{theorem}(intrinsic spectral equivalence) \label{tm:Lambda-K}
Let $X^1 \subset X^0$ be two Hilbert spaces with $X^1$
continuously embedded and dense in $X^0$. Let the self-adjoint and coercive
operator $\Lambda: X^1 \to X^0$ satisfy
\eqref{eq:operator-Lambda} and have a discrete, complete and orthonormal
set of eigenpairs $(\lambda_k,\varphi_k)_{k=1}^\infty$ in $X^0$. Given $s \in (0,1)$,
for any $v \in X^0$ with $\|\Lambda^s v\|_{X^0} < \infty$ we have
\begin{equation} \label{eq:Lambda-K}
 \|\Lambda^s v\|_{X_0} = \|v\|_{(X^0, X^1)_{s,2}}.
\end{equation}
\end{theorem}
\begin{proof}
Given $v = \sum_{k=1}^{\infty} v_k\varphi_k \in X^0$ we split it as $v
  = v^0 + v^1$, with
\[
    v^0 = \sum_{k=1}^{\infty} (1 - a_k)v_k\varphi_k, \qquad
    v^1 = \sum_{k=0}^{\infty} a_k v_k\varphi_k,
\]
and $\{a_k\}_{k=1}^\infty$ to be determined. Combining \eqref{eq:operator-Lambda} with the
definition \eqref{eq:modified-K} of the $K_2$-functional and the
orthonormality of $\{\varphi_k\}_{k=1}^\infty$ in $X^0$ yields
\[
\begin{aligned}
  K_2(t,v)^2 & = \inf_{\{a_k\}_{k=1}^\infty} \Big\| \sum_{k=1}^{\infty}
  (1 - a_k)v_k \varphi_k \Big\|_{X^0}^2
  + t^2 \Big\|\sum_{k=1}^{\infty} a_kv_k\varphi_k \Big\|_{X^1}^2 \\
  & = \inf_{\{a_k\}_{k=1}^\infty} \Big\| \sum_{k=1}^{\infty}
  (1 - a_k)v_k \varphi_k \Big\|_{X^0}^2 + t^2 \Big\|\sum_{k=1}^\infty
  a_k \lambda_k v_k \varphi_k \Big\|_{X^0}^2 \\
  &= \inf_{\{a_k\}_{k=1}^\infty} 
  \sum_{k=1}^{\infty} \Big( (1-a_k)^2 + t^2 a_k^2 \lambda_k^2
  \Big)v_k^2.
\end{aligned}
\]
We choose $a_k = (1 + \lambda_k^2 t^2)^{-1}$, which
minimizes the terms in parenthesis above, to obtain
\[
  K_2(t,v)^2 = \sum_{k=1}^{\infty}
  \frac{\lambda_k^2 t^2}{1+\lambda_k^2 t^2}v_k^2.
\]
Recalling \eqref{eq:interpolation_norms} and
applying the change of variables $\theta = 
\lambda_k t$, we end up with
\begin{align*}
  \|v\|_{(X^0,X^1)_{s,2}}^2 & = \sum_{k=1}^{\infty} \frac{2
  \sin(\pi s)}{\pi} 
  \int_0^\infty t^{-1-2s}  \frac{\lambda_k^2 t^2}{1+ \lambda_k^2t^2} v_k^2
\ dt \\
& = \sum_{k=0}^{\bar{J}} \frac{2 \sin(\pi s)}{\pi} \left(
\int_0^\infty \frac{\theta^{1-2s}}{1+\theta^2} d\theta \right)
  \lambda_k^{2s} v_k^2 = \sum_{k=1}^{\infty}
  \lambda_k^{2s} v_k^2 = \|\Lambda^s v\|_{X_0}^2,
\end{align*}
because $\int_0^\infty \frac{\theta^{1-2s}}{1+\theta^2} d\theta =
\frac{\pi}{2 \sin(\pi s)}$ (see \cite[Exercise B.5]{mclean2000strongly}).
This concludes the proof.
\end{proof}

We now apply Theorem \ref{tm:Lambda-K} (intrinsic spectral
equivalence) to $L^2$-based Sobolev spaces. Let $X^0 = \widetilde
L^2(\Omega)$ and $X^1 = \widetilde H^1_0(\Omega)$ denote the spaces of
functions in $L^2(\Omega)$ and $H^1_0(\Omega)$ extended by zero to
$\Omega^c$, respectively, and let the inner product in $X^1$ be given
by $(v,w)_{X^1} = \int_{\mathbb{R}^d} \nabla v \cdot \nabla w =
\int_\Omega \nabla v \cdot \nabla w$.  The corresponding operator $S$
equals the Laplacian $-\Delta$ with zero Dirichlet condition and
$\Lambda = (-\Delta)^{\frac12}$. Therefore, the $k$-th eigenvalues
$\wh{\lambda}_k$ of $-\Delta$ and $\lambda_k$ of $(-\Delta)^{\frac12}$
satisfy $\wh{\lambda}_k = \lambda_k^2$ whereas the $k$-th
eigenfunctions are the same, whence 
\[
\|\Lambda^s v\|_0^2 = \|(-\Delta)^{\frac{s}{2}}v\|_0^2
= \sum_{k=1}^\infty \lambda_k^{2s} v_k^2 = \sum_{k=1}^\infty
\wh{\lambda}_k^s v_k^2
\]
is the norm square of the interpolation space
\[
\widetilde H^s(\Omega) = \left(\widetilde L^2(\Omega), \widetilde
H^1_0(\Omega) \right)_{s,2}.
\]
Since this norm is equivalent to the Gagliardo norm $|\cdot|_s$
induced by \eqref{eq:bilinear} for $\sigma=s$ with a constant
independent of $s$ (cf. \cite[Theorem B.8, Theorem
B.9]{mclean2000strongly} and \cite{chandler2015interpolation}), we
deduce
\begin{equation}\label{eq:Hs-equiv-Ls}
|v|_s^2 \simeq \|\Lambda^s v\|_0^2 = \sum_{k=1}^\infty \wh{\lambda}_k v_k^2
\quad\forall v\in \widetilde H^s(\Omega).
\end{equation}

%------------------------------------------------------------------------
\subsection{Regularity of solutions} \label{sec:regularity}
%------------------------------------------------------------------------
We next discuss the regularity of solutions to \eqref{eq:Dirichlet} in
either standard or suitably weighted Sobolev spaces.  Grubb's
\cite{Grubb15} accurate elliptic regularity estimates, expressed in
terms of H\"ormander $\mu$-spaces, can be interpreted in the Sobolev
scale but they require the domain to be smooth.  Regularity estimates
valid for arbitrary bounded Lipschitz domains and a right-hand side
function $f \in L^2(\Omega)$ are derived in \cite{BoNo19}.

Reference \cite{ROSe14} studies H\"older regularity of solutions by
using a boundary Harnack method and establishes that, if $f \in
L^\infty(\Omega)$, then the solution to \eqref{eq:Dirichlet} satisfies
$u \in C^s(\overline\Omega)$. This is consistent with the boundary
behavior \cite{Bonforte18,Grubb15,ROSe14}
\begin{equation}\label{eq:boundary-behavior}
  u \sim d(\cdot,\partial\Omega)^s,
\end{equation}
where $d(\cdot,\partial\Omega)$ denotes the distance to $\partial\Omega$.
The sharp characterization
of boundary behavior in \cite{ROSe14} serves as a guide to derive Sobolev
regularity estimates in \cite{AcBo17}.

\begin{proposition}[regularity on Lipschitz domains] \label{prop:reg-Lipschitz}
Let $s \in (0,1)$ and $\Omega$ be a bounded Lipschitz domain
satisfying the exterior ball condition. If $s \in (0,\frac12)$, let $f \in
C^{\frac12-s}(\overline\Omega)$; if $s = \frac12$, let $f \in
L^\infty(\Omega)$; and if $s \in (\frac12,1)$, let $f \in
C^{\beta}(\overline\Omega)$ for some $\beta>0$. Then, for every
$\eps>0$, the solution $u$ to \eqref{eq:Dirichlet} satisfies $u \in
\widetilde H^{s+\frac12-\eps}(\Omega)$, with
\[
|u|_{s+\frac12-\eps} \lesssim \frac{1}{\eps} \|f\|_{\star}.
\]
Above, $\| \cdot \|_\star$ denotes the $C^{\frac12-s}(\overline\Omega)$,
$L^\infty(\Omega)$ or $C^{\beta}(\overline\Omega)$, for $s < \frac12$,
$s=\frac12$ or $s > \frac12$, respectively, and the hidden constant
depends on $\Omega, d$ and $s$. 
\end{proposition}

\begin{remark}[sharpness] \label{rmk:example}
Proposition \ref{prop:reg-Lipschitz} is sharp according to the following
example \cite{Getoor61}. Let $\Omega = B(0,1) \subset \polRd$ and $f =
1$. Then, the solution to \eqref{eq:Dirichlet} is
\[
u(x) =
\frac{\Gamma\left(\frac{d}{2}\right)}{2^{2s}\Gamma\left(\frac{d+2s}{2}\right)\Gamma\left(1+s\right)}
\, (1 - |x|^2)^s_+.
\]
\end{remark}

In view of Proposition \ref{prop:reg-Lipschitz}, we expect that
conforming finite element approximations over quasi-uniform grids
would converge with order $\frac12$ in the energy norm. To mitigate such a
low convergence rate we could increase the mesh grading towards
$\partial\Omega$ and compensate for \eqref{eq:boundary-behavior}.
This idea was exploited in \cite{AcBo17} (see also
\cite{BBNOS18,BoNoSa18}), where the regularity of the solution is
characterized in weighted Sobolev spaces, with the weight being a
power of $d(\cdot,\partial\Omega)$. We refer to either of these
references for a definition of the spaces $\widetilde
H^{t}_{\alpha}(\Omega)$.

\begin{proposition}[regularity in weighted spaces] \label{prop:bdry-regularity} 
Let $\Omega$ be a bounded, Lipschitz domain satisfying the exterior
ball condition. Let $f \in C^{\beta}(\overline\Omega)$ for some $\beta
\in (0,2-2s)$, $\alpha \ge 0$, $t < \min \{ \beta + 2s, \alpha + s +
\frac12 \}$ and $u$ be the solution of \eqref{eq:Dirichlet}. Then, we
  have $u \in \widetilde H^{t}_{\alpha}(\Omega)$, with
\[
\|u\|_{\widetilde H^{t}_{\alpha}(\Omega)} \le \frac{C(\Omega,d,s)}{\sqrt{(\beta + 2s - t) \, (1+2(\alpha + s - t))}} \| f \|_{C^{\beta}(\overline\Omega)}.
\]
\end{proposition}

\begin{remark}[optimal parameters] \label{rm:optimal-parameters} The
  optimal choice of parameters $t$ and $\alpha$ for finite element
  applications depends on either the smoothness of $f$ and the
  dimension of the space. For instance, in dimension $d = 2$ and a
  sufficiently smooth right-hand side $f$, one can set $t = 1 + s -
  2\eps$ and $\alpha = \frac12 - \eps$ for an arbitrary $\eps \in
  (0,\frac12)$. The resulting constant scales as $\eps^{-\frac12}$ and one
  obtains a linear convergence rate (up to logarithmic terms) with
  respect to the dimension of the finite element spaces, which is
  optimal for approximations on shape-regular meshes.  We refer to
  \cite{BoLeNo20} for a thorough discussion on this aspect.
\end{remark}

%------------------------------------------------------------------------
\subsection{Finite element discretization} \label{sec:discrete_setting}
%------------------------------------------------------------------------
Given a conforming and shape-regular triangulation $\cT$ of $\Omega$,
we consider discrete spaces consisting of continuous piecewise linear
functions that vanish on $\partial \Omega$,
\[
\mathbb{V}(\cT) = \{ v_h \in C(\overline\Omega) \colon v_h|_T \in P_1(T) \
\forall T \in \cT, \ v_h |_{\partial\Omega} = 0 \} .
\]
It is clear that $\mathbb{V}(\cT) \subset \widetilde H^s(\Omega)$,
independently of the value of $s$. Therefore, we can pose a conforming
discretization of \eqref{eq:weak}: we seek $u_h \in \mathbb{V}(\cT)$
such that
\begin{equation} \label{eq:disc-weak}
  a(u_h, v_h) = \langle f,v_h \rangle_{s,\Omega} \quad \forall v_h \in
  \mathbb{V}(\cT).
\end{equation}
Thus, the finite element solution is the elliptic projection of the
solution $u$ to \eqref{eq:Dirichlet} onto the discrete space
$\mathbb{V}(\cT)$, 
\[
| u - u_h |_s = \inf_{v_h \in \mathbb{V}(\cT)} | u - v_h |_s \quad
\forall v_h \in \mathbb{V}(\cT).
\]

Convergence rates in the energy norm are derived by combining the
estimate above with suitable interpolation estimates
\cite{CiarletJr,AcBo17} and the regularity described in Proposition
\ref{prop:reg-Lipschitz} (regularity on Lipschitz domains); cf.
\cite[Theorem 3.7]{BBNOS18}.

\begin{proposition}[convergence rates in uniform meshes]
  \label{prop:convergence_uniform}
Assume $s \in (0,1)$ and $\Omega$ is a bounded Lipschitz domain. Let
$u$ denote the solution to \eqref{eq:weak} and denote by $u_h \in
\mathbb{V}(\cT)$ the solution of the discrete problem
\eqref{eq:disc-weak}, computed over a mesh $\cT$ consisting of
elements with maximum diameter $h$. 
{Under the hypotheses of Proposition \ref{prop:reg-Lipschitz}, we have
\[
|u - u_h |_s  \lesssim h^{\frac12} |\log h|^{1+\kappa} \|f\|_{\star}.
\]
Above, $\| \cdot \|_\star$ denotes the $C^{\frac12-s}(\overline\Omega)$,
$L^\infty(\Omega)$ or $C^{\beta}(\overline\Omega)$, depending on
whether $s < \frac12$, $s=\frac12$ or $s > \frac12$, and $\kappa = 1$ if $s=\frac12$
and zero otherwise.}
\end{proposition}

We point out that, on quasi-uniform grids, this approximation rate is
optimal due to \eqref{eq:boundary-behavior}.  As shown in
\cite{AcBo17}, when solving \eqref{eq:Dirichlet} it is possible to
increase the a priori convergence rates by utilizing suitably graded
grids and making use of the regularity estimate from Proposition
\ref{prop:bdry-regularity} (regularity in weighted spaces).  More
precisely, given a grading parameter $\mu \ge 1$ and a mesh size
parameter $h$, assume that the element size $h_T$ satisfies
\begin{equation} \label{eq:grading}
 h_T \simeq \left\lbrace
\begin{array}{ll}
 h^{\mu} &
\mbox{if } S_T\cap \partial \Omega\neq \varnothing, \\
h  \, d(T,\partial \Omega)^{(\mu-1)/\mu} & \mbox{otherwise.}
 \end{array}
 \right.
 \end{equation}
 
We wish the mesh size parameter $h$ to be such that $\mbox{dim}
\mathbb{V}(\cT) \simeq h^{-d}$. Shape regularity limits the
range for $\mu$ when $d \ge 2$: we have
\begin{equation*} \label{eq:dofs}
\mbox{dim} \mathbb{V}(\cT) \simeq \left\lbrace \begin{array}{ll}
h^{-d}, & \mbox{ if } \mu \in \big[1,\frac{d}{d-1}\big), \\
h^{-d} |\log h| & \mbox{ if } \mu = \frac{d}{d-1}, \\
h^{(1-d)\mu}  & \mbox{ if } \mu > \frac{d}{d-1}. \\
\end{array} \right.
\end{equation*}
Thus, for a sufficiently smooth right hand side $f$, the optimal
choice for $\mu$ is $\frac{d}{d-1}$, and one obtains the following
convergence rates with this strategy.

\begin{proposition}[convergence rates in graded meshes]
\label{prop:convergence-graded}
Let $s \in (0,1)$ and $\Omega \subset \polR^d$ be a bounded
Lipschitz domain satisfying the exterior ball condition,  $u$
be the solution to \eqref{eq:weak} and denote by $u_h \in
\mathbb{V}(\cT)$ the solution of the discrete problem
\eqref{eq:disc-weak}. Let $\beta >  0$ be such that
\begin{equation*} \label{eq:beta_mu}
\beta \ge \left\lbrace
\begin{array}{rl}
2-2s & \mbox{if } d  = 1, \\
\frac{d}{2(d-1)} - s & \mbox{if } d \ge 2,
\end{array}
\right.
\quad \mbox{and } \quad
\mu = \left\lbrace
\begin{array}{rl}
2-s & \mbox{if } d  = 1, \\
\frac{d}{d-1} & \mbox{if } d \ge 2.
\end{array} \right.
\end{equation*}

Then, if $f \in C^{\beta}(\overline{\Omega})$, and the mesh $\cT$ is
graded according to \eqref{eq:grading} , we have
\[
\|u - u_h \|_s  \lesssim 
\left\lbrace
\begin{array}{rl}
h^{2-s} | \log h|^\kappa \|f\|_{C^{\beta}(\overline{\Omega})} & \mbox{if } d = 1, \\
h^{\frac{d}{2(d-1)}} |\log h|^{1+\kappa} \|f\|_{C^{\beta}(\overline{\Omega})} & \mbox{if } d \ge 2.
\end{array} \right.
\]
where the hidden constant depends on $\Omega$, $s$ and the shape
regularity of $\cT$ and $\kappa = 1$ if $s=\frac12$ and zero otherwise.
\end{proposition}

%------------------------------------------------------------------------
\section{Robust Additive Multilevel Preconditioning}\label{sec:Robust-BPX}
%------------------------------------------------------------------------
Let $(\cdot,\cdot)$ be the $L^2$-inner product in $\Omega$ and $V:=\mathbb{V}(\cT)$
denote the discrete space. Let $A:V\to V$ be the symmetric positive definite (SPD)
operator defined by $(Au, v) := a(u, v)$ for any $u, v \in V$, and let
$\wt{f} \in V$ be given by $(\wt{f},v) = \langle f, v
\rangle_{s,\Omega}$ for any $v\in V$.  With this notation at hand, the
discretization \eqref{eq:disc-weak} leads to the following linear
equation in $V$
\begin{equation} \label{eq:operator-eq}
  Au = \wt{f}.
\end{equation}
In this section, we give some general and basic results that will be used to
construct the additive multilevel preconditioners for \eqref{eq:operator-eq}.

%------------------------------------------------------------------------
\subsection{Space decomposition}\label{S:space-decomp}
%------------------------------------------------------------------------
We now invoke the {\it method of subspace corrections}
\cite{xu1992iterative,xu1996auxiliary, xu2002method}.
We first decompose the space $V$ as the sum of subspaces $V_j \subset V$
\[ 
V = \sum_{j=0}^J V_j.
\] 
For $j = 0, 1, \ldots, J$, we consider the following operators: 
\begin{itemize} 
\item $Q_j \colon V \to V_j$ is the $L^2$-projection operator defined by
  $(Q_jv,v_j)=(v,v_j)$ for all $v\in V, v_j\in V_j$;
\item $I_j \colon V_j \to V$ is the natural inclusion operator given by $I_j v_j = v_j$
  for all $v_j\in V_j$; 
  \item $R_j \colon V_j \to V_j$ is an approximate inverse of the restriction
    of $A$ to $V_j$ (often known as smoother); we set 
    $\|v_j\|_{R_j^{-1}}^2 := (R_j^{-1} v_j,v_j)$ for all $v_j \in V_j$
    provided that $R_j$ is SPD on $V_j$.
\end{itemize}
A straightforward calculation shows that $Q_j = I_j^t$ because
$(Q_j v, v_j) = (v, I_j v_j) = (I_j^tv, v_j)$ for all $v\in V, v_j\in V_j$.
Let the fictitious space be $\utilde{V} = V_0
\times V_1 \times \ldots \times V_J$.  Then, the 
Parallel Subspace Correction (PSC) preconditioner $B \colon V \to V$ is
defined by
\begin{equation} \label{eq:fictitious}
B := \sum_{j=0}^J I_jR_jQ_j = \sum_{j=0}^J I_jR_jI_j^t .
\end{equation}

The next two lemmas follow from the general theory of preconditioning
techniques based on fictitious or auxiliary spaces
\cite{nepomnyaschikh1992decomposition, griebel1995abstract,
xu1992iterative, xu1996auxiliary, Xu97, xu2002method}. For
completeness, we give their proofs in Appendix \ref{sec:PSC-lemmas}.

\begin{lemma}[identity for PSC] \label{lm:PSC} 
If $R_j$ is SPD on $V_j$ for $j = 0, 1, \ldots, J$, then $B$ defined
in \eqref{eq:fictitious} is also SPD under the inner product
$(\cdot,\cdot)$.  Furthermore, 
\begin{equation} \label{eq:PSC}
  (B^{-1}v, v) = \inf_{\sum_{j=0}^J v_j = v} \sum_{j=0}^J
  (R_j^{-1}v_j, v_j) 
  \quad \forall v \in V.
\end{equation}
\end{lemma}

\begin{lemma}[estimate on $\cond (BA)$] \label{lm:auxiliary}
If the operator $B$ in \eqref{eq:fictitious} satisfies
\begin{enumerate}[]
\item (A1) Stable decomposition: for every $v\in V$, there exists
$(v_j)_{j=0}^J \in \utilde{V}$ such that
$\sum_{j=0}^J v_j = v$ and
\begin{equation} \label{eq:stable-decomp}
\sum_{j=0}^J \|v_j\|_{R_j^{-1}}^2 \leq c_0 \|v\|_{A}^2,
\end{equation}
where $\|v\|_{A}^2 = (Av,v)$, then $\lambda_{\min}(BA) \geq c_0^{-1}$;
\item (A2) Boundedness: For every $(v_j)_{j=0}^J  \in \utilde{V}$ there holds
\begin{equation} \label{eq:boundedness}
\left\| \sum_{j=0}^J v_j \right\|_{A}^2 \leq c_1 \sum_{j=0}^J
\|v_j\|_{R_j^{-1}}^2,
\end{equation}
then $\lambda_{\max}(BA) \leq c_1$.
\end{enumerate}
Consequently, if $B$ satisfies (A1) and (A2), then
$\cond (BA) \leq c_0c_1$.
\end{lemma}
  
%-------------------------------------------------------------------------
\subsection{Instrumental tools for $s$-uniform preconditoner}
%-------------------------------------------------------------------------
We assume that the spaces $\{V_j\}_{j=0}^J$ are nested, i.e.
\begin{equation*}\label{nested}  
V_{j-1}\subset V_j \quad \forall 1\le j \le J.
\end{equation*}
With the convention that $Q_{-1} = 0$, we consider the $L^2$-slicing operators
\begin{equation*}\label{slicing}
\wt{Q}_j \colon V \to
V_j: \quad \wt{Q}_j := Q_j - Q_{j-1} \quad(j = 0, 1,\ldots,J).
\end{equation*}
Clearly, the $L^2$-orthogonality implies that $Q_k Q_{j} = Q_{k\wedge
j}$, where $k \wedge j := \min\{k, j\}$.  Hence,
\begin{equation} \label{eq:Q-identity}
\wt{Q}_j Q_k = 
Q_k \wt{Q}_j = 
  \begin{cases}
    \wt{Q}_j & j \leq k, \\
    0 &  j > k,
  \end{cases}
\qquad 
  \wt{Q}_k \wt{Q}_j = \delta_{k j} \wt{Q}_k.
\end{equation}

The following lemma plays a key role in the analysis of an $s$-uniform
preconditioner, which is obtained by using the identity of PSC
\eqref{eq:PSC} and reordering the BPX preconditioner
\cite{xu1989thesis, bramble1990parallel}.

\begin{lemma}[$s$-uniform decomposition] \label{lm:s-decomp}
Given $\gamma \in (0,1)$, $0<s\le1$, it holds that, for every $v\in V$,  
\begin{equation*} \label{eq:s-decomp} 
\sum_{j=0}^J \gamma^{-2sj}\|(Q_j - Q_{j-1})v\|_0^2 =
\inf_{\substack{v_j \in V_j \\ 
\sum_{j=0}^J v_j = v}} \left[ \gamma^{-2sJ}\|v_J\|_0^2 +  
\sum_{j = 0}^{J-1} \frac{\gamma^{-2sj}}{1 - \gamma^{2s}}\|v_j\|_0^2
\right].
\end{equation*} 
\end{lemma}
\begin{proof} 
This proof is an application of Lemma \ref{lm:PSC} (identity for PSC).
Taking 
$$
B = \sum_{j=0}^J \gamma^{2sj}(Q_j - Q_{j-1}),
$$ 
by the $L^2$-orthogonality \eqref{eq:Q-identity}, we easily see that
$B^{-1} = \sum_{j=0}^J \gamma^{-2sj}(Q_j - Q_{j-1})$ and
$$ 
  (B^{-1}v, v) = \sum_{j=0}^J \gamma^{-2sj} \|(Q_j - Q_{j-1}) v\|_0^2.
$$ 
On the other hand, to identify $R_j$ we reorder the sum in the
definition of $B$ 
$$ 
B = \sum_{j=0}^J \gamma^{2sj}(Q_j - Q_{j-1}) = \gamma^{2sJ}Q_J +
  \sum_{j=0}^{J-1}(1-\gamma^{2s})\gamma^{2sj} Q_j = \sum_{j=0}^J I_j
  R_j Q_j, 
$$ 
where
$$ 
  R_j v_j := 
  \begin{cases}
    (1-\gamma^{2s})\gamma^{2sj} v_j & \quad  j = 0,\ldots, J-1, \\
     \gamma^{2sj}v_j & \quad j = J,
  \end{cases}
$$
for all $v_j \in V_j$.
Finally, the identity \eqref{eq:PSC} of PSC gives the desired result.
\end{proof}

The next lemma is an application of space interpolation theory and is
crucial to obtain stable decompositions in fractional-order
norms. We postpone its proof to Appendix \ref{sec:auxiliary-lemmas}.

\begin{lemma}[$s$-uniform interpolation] \label{lm:norm-equiv1}
Assume that the spaces $\{V_j\}_{j=0}^{J}$ are nested, and
\begin{equation} \label{eq:interpolation-s1}
 \sum_{j=0}^J \gamma^{-2j} \|(Q_j-Q_{j-1}) v\|_0^2 \lesssim |v|_1^2 \qquad \forall v \in V.
\end{equation}
Then, the following inequality holds, with the hidden constant
independent of $s$ and $J$,
\begin{equation} \label{eq:norm-equiv1}
\sum_{j=0}^J \gamma^{-2sj} \|(Q_j-Q_{j-1}) v\|_0^2 \lesssim |v|_s^2 \qquad \forall v \in V.
\end{equation}
\end{lemma}

We conclude this section with a standard {\it local inverse estimate}
valid on graded grids $\cT$. The proof is elementary, and we
give it in Appendix \ref{sec:auxiliary-lemmas}. Given
$\tau\in\mathcal{T}$, we define
\[
S_\tau := \bigcup \{\tau' \in \mathcal{T} \colon \overline{\tau'}\cap
  \overline{\tau} \neq \varnothing\}.
\]

\begin{lemma}[local inverse inequality] \label{lm:inv_ineq}
Let $\sigma \in [0,3/2)$ and $\mu \in [0, \sigma]$. Then,
\begin{equation} \label{eq:inv_ineq}
  |v|_\sigma \lesssim \left(\sum_{\tau\in\mathcal{T}} h_\tau^{2(\mu-\sigma)}
|v|_{H^\mu(S_\tau)}^2 \right)^{\frac12}\quad \forall
v \in \mathbb{V}(\cT),
\end{equation}
where the hidden constant only blows up as $\sigma
\to 3/2$.
\end{lemma}

%------------------------------------------------------------------------
%------------------------------------------------------------------------
\section{Robust BPX preconditioner for quasi-uniform grids} 
\label{sec:BPX-uniform}
%------------------------------------------------------------------------
%------------------------------------------------------------------------
In this section, we propose and study a BPX preconditioner
\cite{xu1989thesis,bramble1990parallel, Melenk:2019} for the solution
of the systems arising from the finite element discretizations
\eqref{eq:disc-weak} on quasi-uniform grids. We emphasize that, in
contrast to \cite{Melenk:2019}, the proposed preconditioner is {\it
uniform} with respect to both the number of levels and the order $s$.
To this end, we introduce a new scaling for coarse spaces which
differs from the original BPX preconditioners given in
\cite{bramble1990parallel, Melenk:2019}.  Our theory applies as well
to the spectral and censored fractional Laplacians except for the
censored one when $s\to \frac12$; see Section
\ref{sec:spectral-equivalence}.

Consider a family of uniformly refined grids
$\{\overline{\cT}_k\}_{k=0}^{\bar{J}}$ on $\Omega$, where
$\overline{\cT}_0 = \cT_0$ is a quasi-uniform initial triangulation.
On each of these grids we define the space $\overline{V}_k :=
\mathbb{V}(\overline{\cT}_k)$. Let $\overline{V} =
\overline{V}_{\bar{J}}$ and $\overline{A}$ be the SPD operator on
$\overline{V}$ associated with $a(\cdot, \cdot)$: $(\overline{A}
v,w)=a(v,w)$ for all $v,w\in\overline{V}$. Let the grid size be
$\bar{h}_k \simeq \gamma^k$, where $\gamma \in (0,1)$ is a fixed
constant. For instance, we have $\gamma = \frac12$ for uniform refinement,
in which each simplex is refined into $2^d$ children, and $\gamma =
(\frac12)^{1/d}$ for uniform bisection, in which each simplex is refined
into $2$ children.

Let $\overline{Q}_k: \overline{V}\to \overline{V}_k$ and
$\overline{I}_k: \overline{V}_k \to \overline{V}$ be the
$L^2$-projection and inclusion operators defined in Section
\ref{S:space-decomp}, and let $\overline{Q}_{-1} := 0$.
The standard BPX preconditioner reads
\cite{bramble1990parallel, Melenk:2019}
\begin{equation}\label{eq:standard-BPX}
\overline{B} = \sum_{k=0}^{\bar{J}} \overline{I}_k \bar{h}_k^{2s}
  \overline{Q}_k: \overline{V} \to \overline{V}.
\end{equation}    
A rough analysis of \eqref{eq:standard-BPX} proceeds as follows. Let
$\overline{S}_k: \overline{V} \to \overline{V}_k$ be the Scott-Zhang
interpolation operator \cite{scott1990finite} and let $\overline{S}_k
- \overline{S}_{k-1}\colon \overline{V} \to \overline{V}_k$ be the
slicing operator for all $k = 0, 1, \ldots, \bar{J}$ with
$\overline{S}_{-1} := 0$.  We thus have the decomposition of any $v\in
\overline{V}$:
\[
  v=\sum_{k=0}^{\bar{J}}v_k, \qquad v_k: = (\overline{S}_k -
  \overline{S}_{k-1}) v.
\]
If $\overline{R}_k v_k := \bar{h}_k^{2s} v_k$, then the stable
decomposition \eqref{eq:stable-decomp} is a consequence of
\begin{align*}
\sum_{k=0}^{\bar{J}} \|v_k\|_{\overline{R}_k^{-1}}^2 
  = \sum_{k=0}^{\bar{J}} \bar{h}_k^{-2s}\|(\overline{S}_k -
  \overline{S}_{k-1}) v\|_{0}^2
   \lesssim \sum_{k=0}^{\bar{J}} \bar{h}_k^{-2s}\|v-
  \overline{S}_k v\|_{0}^2 \lesssim \sum_{k=0}^{\bar{J}}
  |v|_s^2 \le \bar{J} |v|_s^2,
\end{align*}
whence $c_0\lesssim \bar{J}$. On the other hand, the boundedness
\eqref{eq:boundedness} follows from an inverse estimate
\begin{equation} \label{eq:boundedness-J-dependent}
\Big| \sum_{k=0}^{\bar{J}} v_k \Big|^2_s
\leq \bar{J} \sum_{k=0}^{\bar{J}}
|v_k |_{s}^2 \lesssim \bar{J}
\sum_{k=0}^{\bar{J}} \bar{h}_{k}^{-2s}\|v_k\|_0^2 
= \bar{J} \sum_{k=0}^{\bar{J}} \|v_k\|_{\overline{R}_k^{-1}}^2,
\end{equation}
whence $c_1 \lesssim \bar{J}$. Therefore, Lemma \ref{lm:auxiliary}
(estimate of $\cond (BA)$) yields the condition number estimate
$\cond (\overline{B} \, \overline{A}) \lesssim \bar{J}^2$ but
independent of $s$. To remove the dependence on $\bar{J}$ we deal
below with the slicing $L^2$-projectors
$\overline{Q}_k-\overline{Q}_{k-1}$.  However, a
naive replacement of $\overline{S}_k - \overline{S}_{k-1}$ by
$\overline{Q}_k - \overline{Q}_{k-1}$ would make $\cond (\overline{B}
\, \overline{A})$ independent of
$\bar{J}$ but blow-up as $s\to0$. This is an unnatural dependence on
$s$ because $(-\Delta)^s$ tends to the identity as $s\to0$. We
circumvent this issue by a suitable rescaling of coarse levels and
redefinition of the smoothers $\overline{R}_k$.

Let $\wt{\gamma}\in (0,1)$ be a fixed constant; it can be taken equal to
$\gamma$ but this is not needed. For every $v_k \in \overline{V}_k, k
= 0, \ldots, \bar{J}$, we define $\overline{R}_k: \overline{V}_k \to
\overline{V}_k$ to be
\begin{equation} \label{eq:j-smoother2}
\overline{R}_k v_k:= 
\begin{cases}
(1 - \wt{\gamma}^{s}) \bar{h}_k^{2s} v_k  & \quad k
= 0, \ldots, \bar{J}-1, \\
\bar{h}_k^{2s} v_k &\quad k = \bar{J}.
\end{cases}
\end{equation}
We now introduce the BPX preconditioner and study its properties in
the sequel
\begin{equation} \label{eq:BPX-FPDE}
  \overline{B} := \sum_{k=0}^{\bar{J}} \overline{I}_k \overline{R}_k
  \overline{I}_k^t = \overline{I}_{\bar{J}}
  \bar{h}_{\bar{J}}^{2s} \overline{Q}_{\bar{J}} +
  (1-\wt{\gamma}^{s})\sum_{k=0}^{\bar{J} - 1} \overline{I}_k \bar{h}_k^{2s}
  \overline{Q}_k.
\end{equation}

Our next goal is to prove the following theorem, namely that
$\overline{B}$ satisfies the two necessary conditions
\eqref{eq:stable-decomp} and \eqref{eq:boundedness} of Lemma
\ref{lm:auxiliary} (estimate of $\cond (BA)$) uniformly in $\bar{J}$
and $s$ over quasi-uniform grids.  We observe that the scaling $(1 -
\wt{\gamma}^{s})^{-1}>1$ makes it easier to prove
\eqref{eq:boundedness} but complicates \eqref{eq:stable-decomp}. We
prove \eqref{eq:stable-decomp} in Section \ref{S:stable-decomposition}
and \eqref{eq:boundedness} in Section \ref{S:boundedness}. 

\begin{theorem}[uniform preconditioning on quasi-uniform grids]
\label{tm:uniform-cond}
Let $\Omega$ be a bounded Lipschitz domain and $s \in (0,1)$. Consider
discretizations to \eqref{eq:Dirichlet} using piecewise linear
Lagrangian finite elements on quasi-uniform grids. Then, the
preconditioner \eqref{eq:BPX-FPDE} satisfies $\cond (\overline{B}
\, \overline{A}) \lesssim 1$, where the hidden constant is uniform with
respect to both $\bar{J}$ and $s$.
\end{theorem}

%------------------------------------------------------------------------
\subsection{Stable decomposition: Proof of \eqref{eq:stable-decomp}
for quasi-uniform grids}\label{S:stable-decomposition}
%------------------------------------------------------------------------
We start with a norm equivalence for discrete functions. We rely
on operator interpolation and the decomposition for $s=1$
\cite{xu1992iterative,oswald1992norm,bornemann1993basic} , which was
proposed earlier in \cite{xu1989thesis, bramble1990parallel} with a
removable logarithmic factor. A similar result, for the interpolation
norm of $(L^2(\Omega), H_0^1(\Omega))_{s,2}$, was given in
\cite[Theorem 10.5]{Xu97}.

\begin{theorem}[norm equivalence] \label{tm:norm-equivalence}
Let $\Omega$ be a bounded Lipschitz domain and $s \in (0,1)$.
  If $\overline{Q}_k: \overline{V} \to \overline{V}_k$ denotes the
  $L^2$-projection operators onto discrete spaces $\overline{V}_k$,
  and $\overline{Q}_{-1} := 0$, 
  then for any $v\in \overline{V}$ the decomposition $v =
  \sum_{k=0}^{\bar{J}} (\overline{Q}_k-\overline{Q}_{k-1}) v$
  satisfies
\begin{equation}\label{eq:norm-equivalence}
|v|_{s}^2\simeq \sum_{k=0}^{\bar{J}} \bar{h}_k^{-2s}
\|(\overline{Q}_k-\overline{Q}_{k-1}) v\|_{0}^2.
\end{equation}
The equivalence hidden constant is independent of $s$ and $\bar{J}$.
\end{theorem}

\begin{proof}
Consider the operator $\Lambda = \sum_{k=0}^{\bar{J}} \bar{h}_k^{-1}
  (\overline{Q}_k - \overline{Q}_{k-1}): \wt{H}_0^1(\Omega) \to
  \wt{L}^2(\Omega)$, which happens to be self-adjoint and coercive in
  $\wt{L}^2(\Omega)$. Combining the Poincar\'e inequality in
  $\widetilde{H}_0^1(\Omega)$ with the $H^1$-norm equivalence
  \cite{xu1992iterative,oswald1992norm,bornemann1993basic} yields
$$
|v|_1^2 \simeq \|v\|_1^2 \simeq \sum_{k=0}^{\bar{J}} \bar{h}_k^{-2}\|
  (\overline{Q}_k-\overline{Q}_{k-1}) v \|_0^2  = \|\Lambda v\|_0^2
  \quad \forall v \in \widetilde{H}_0^1(\Omega),
$$
whence $|v|_s \simeq \|\Lambda^s v\|_0$ according to
  \eqref{eq:Hs-equiv-Ls}.  It remains to characterize $\|\Lambda^s
  v\|_0$. To this end, notice that $\wt{V}_k := (\overline{Q}_k -
  \overline{Q}_{k-1}) \overline{V}$ is an eigenspace of $\Lambda$ with
  eigenvalue $\lambda_k=\bar{h}_k^{-1}$, namely $\Lambda|_{\wt{V}_k} =
  \bar{h}_k^{-1}I$. Moreover, $\overline{V} = \oplus_{k=0}^{\bar{J}}
  \wt{V}_k$ is an $L^2$-orthogonal decomposition of $\overline{V}$ and
  $\Lambda^s|_{\wt{V}_k} = \bar{h}_k^{-s}I$ according to
  \eqref{eq:Lambda-s}. This implies $\|\Lambda^s v\|_0^2 =
  \sum_{k=0}^{\bar{J}} \bar{h}_k^{-2s}
  \|(\overline{Q}_k-\overline{Q}_{k-1}) v\|_{0}^2$ and thus
  \eqref{eq:norm-equivalence}, as asserted.
\end{proof}

\begin{corollary}[stable decomposition]\label{P:stable-decomposition}
For every $v \in \overline{V}=\overline{V}_{\bar{J}}$, there exists a
decomposition $(v_0,\ldots, v_{\bar{J}}) \in \overline{V}_0 \times
\ldots \times \overline{V}_{\bar{J}}$, such that $\sum_{k=0}^{\bar{J}}
v_k = v$ and
\begin{equation*} \label{eq:decomp-barJ-s}
\bar{h}_{\bar{J}}^{-2s}\|v_{\bar{J}}\|_0^2 +
\frac{1}{1-\wt{\gamma}^s}\sum_{k=0}^{\bar{J} - 1} \bar{h}_{k}^{-2s}
\|v_k\|_0^2 
\simeq  |v|_{s}^2.
\end{equation*}
\end{corollary}
\begin{proof}
This is a direct consequence of Lemma \ref{lm:s-decomp} ($s$-uniform
decomposition) and Theorem \ref{tm:norm-equivalence} (norm
equivalence) because $\bar{h}_k \simeq \gamma^k$.
\end{proof}

%------------------------------------------------------------------------
\subsection{Boundedness: Proof of \eqref{eq:boundedness} for quasi-uniform grids}
\label{S:boundedness}
%------------------------------------------------------------------------
We now prove the boundedness estimate in Lemma \ref{lm:auxiliary}
(estimate on $\cond (BA)$) with a constant independent of both
$\bar{J}$ and $s$. 

\begin{proposition}[boundedness]\label{P:boundedness}
The preconditioner $\overline{B}$ in \eqref{eq:BPX-FPDE} satisfies
\eqref{eq:boundedness}, namely  
\begin{equation}\label{eq:uniform-bd}
  \Big| \sum_{k=0}^{\bar{J}} v_k \Big|^2_s \le c_1
  \bigg(\bar{h}_{\bar{J}}^{-2s}\|v_{\bar{J}}\|_0^2 + 
  \frac{1}{1-\wt{\gamma}^s}\sum_{k = 0}^{\bar{J} - 1}
  \bar{h}_k^{-2s}\|v_k\|_0^2\bigg),
\end{equation}
where $\wt{\gamma} \in (0,1)$ can be taken arbitrarily and the
constant $c_1$ is independent of $\bar{J}$ and $s$.
\end{proposition}
\begin{proof}
Let $v:= \sum_{k=0}^{\bar{J}} v_k$. Then, we use Theorem
\ref{tm:norm-equivalence} (norm equivalence), the fact that $\bar{h}_k
\simeq \gamma^{k}$ and Lemma \ref{lm:s-decomp} ($s$-uniform
decomposition) to write
\[
\begin{aligned}
  \Big|\sum_{k=0}^{\bar{J}} v_k\Big|_s^2 &= |v|_s^2 \simeq
  \sum_{k=0}^{\bar{J}} \bar{h}_k^{-2s} \|(\overline{Q}_k -
  \overline{Q}_{k-1} )v\|_0^2 \\
  &\simeq \sum_{k=0}^{\bar{J}} \gamma^{-2sk} \|(\overline{Q}_k -
  \overline{Q}_{k-1} )v\|_0^2
  = \inf_{\substack{w_k \in \overline{V}_k \\ \sum_{k=0}^{\bar{J}} w_k = v}}
  \left[ \gamma^{-2s\bar{J}}\|w_{\bar{J}}\|_0^2 +  
  \sum_{k = 0}^{\bar{J}-1} \frac{\gamma^{-2sk}}{1 - \gamma^{2s}}\|w_k\|_0^2 
  \right].   
\end{aligned}
\]
Therefore, upon setting $w_k = v_k$ for $j = 0, \ldots \bar{J}$ above,
we deduce that
\[
\Big|\sum_{k=0}^{\bar{J}} v_k\Big|_s^2 \lesssim
\gamma^{-2s\bar{J}}\|v_{\bar{J}}\|_0^2 +  
  \sum_{k = 0}^{\bar{J}-1} \frac{\gamma^{-2sk}}{1 -
  \gamma^{2s}}\|v_k\|_0^2 \leq  c_1
  \bigg(\bar{h}_{\bar{J}}^{-2s}\|v_{\bar{J}}\|_0^2 + 
  \frac{1}{1-\wt{\gamma}^s}\sum_{k = 0}^{\bar{J} - 1}
  \bar{h}_k^{-2s}\|v_k\|_0^2\bigg).
\]
The proof is thus complete. 
\end{proof}

An alternative derivation of boundedness could be carried with the aid
of a {\em strengthened Cauchy-Schwarz inequality}, which plays an
important role in the analyis of multigrid methods (cf.
\cite{xu1989thesis, xu1992iterative, Xu97}). We provide a proof of
such an inequality together with a second proof of Proposition
\ref{P:boundedness} (boundedness) in  Appendix \ref{sec:SCS}.  This
tool also allows us to illustrate the need of the correction factor
$1-\wt{\gamma}^{s}$ in the coarser scales in \eqref{eq:BPX-FPDE}.
%The usual proof for
%second-order problems consists of an elementwise integration-by-parts
%argument, a local argument that quantifies the interaction between
%functions with different frequencies. This is not possible in the
%present context due to the nonlocal nature of the fractional norms.
%We resort instead to the well-known characterization of the fractional
%Sobolev space $\wt{H}^s(\Omega)$ as a Bessel potential space; see also
%\cite{xunotes,jiang2015multigrid}.

\begin{remark}[preconditioner \eqref{eq:standard-BPX}]\label{R:standard-BPX}
We wonder how the boundedness \eqref{eq:boundedness} changes if we
consider the standard preconditioner \eqref{eq:standard-BPX} instead
of \eqref{eq:BPX-FPDE}. Since $\|v_k\|_{\overline{R}_k^{-1}}^2 =
\bar{h}_k^{-2s}\|v_k\|_0^2$ in this case, Lemma
\ref{lm:generalized-SCS} (generalized strengthened Cauchy-Schwarz
inequality) and \eqref{eq:gamma-matrix} yield
\begin{equation*}
\Big|\sum_{k=0}^{\bar{J}} v_k\Big|_s^2 =
\sum_{k,\ell=0}^{\bar{J}} (v_k, v_\ell)_s 
\lesssim 
\sum_{k,\ell=0}^{\bar{J}} \gamma^{s|k-\ell|} \bar{h}_k^{-s}
\bar{h}_\ell^{-s} \|v_k\|_{0}\|v_\ell\|_0
\lesssim \frac{1}{1-\gamma^s} \sum_{k=0}^{\bar{J}}
\bar{h}_k^{-2s}\|v_k\|_0^2.
\end{equation*}
This together with \eqref{eq:boundedness-J-dependent} implies that the
constant $c_1 \lesssim \min\{ (1-\gamma^s)^{-1}, \bar{J} \}$ of
\eqref{eq:boundedness} blows up as $s\to0$ and $\bar{J} \to \infty$;
this is observed in the experimental results reported in Table
\ref{tab:comparison}.
\end{remark}

%--------------------------------------------------------------------
\subsection{Spectral equivalence: Spectral and censored
Laplacians}\label{sec:spectral-equivalence}
%--------------------------------------------------------------------
We now exploit the fact that Theorem \ref{tm:norm-equivalence} (norm
equivalence) is insensitive to the number of levels $\bar{J}$, whence
letting $\bar{J}\to\infty$ we obtain the multilevel decomposition of
any $v \in \wt{H}^s(\Omega)$
\begin{equation}\label{eq:J=infty}
|v|_{s}^2\simeq \sum_{k=0}^{\infty} \bar{h}_k^{-2s}
\|(\overline{Q}_k-\overline{Q}_{k-1}) v\|_{0}^2.
\end{equation}
An alternative definition of $|\cdot|_{s}$, but equivalent to
\eqref{eq:bilinear}, relies on the spectral decomposition of the
Laplacian $-\Delta$ in a bounded Lipschitz domain $\Omega$.
Recall from Section \ref{sec:interpolation} that if
$(\wh{\lambda}_k,\varphi_k)_{k=1}^\infty$ is the
sequence of eigenpairs of $-\Delta$ with zero Dirichlet boundary
condition and normalized in $L^2(\Omega)$, then \eqref{eq:Hs-equiv-Ls}
implies that the space
\begin{equation*}
  \widehat{H}^s(\Omega) := \Big\{
  v = \sum_{k=1}^\infty v_k \varphi_k \in L^2(\Omega): \quad
  |v|_{\widehat{H}^s(\Omega)}^2 = \sum_{k=1}^\infty \wh{\lambda}_k^s v_k^2 < \infty
  \Big\}
\end{equation*}
coincides with $\wt{H}^s(\Omega)$ and has equivalent norms. However,
these norms induce different fractional operators. Minima of the
functional $v \mapsto \frac12 |v|_{\widehat{H}^s(\Omega)}^2 -
\int_\Omega f v$ are weak solutions of the {\it spectral} fractional
Laplacian in $\Omega$ with homogeneous Dirichlet condition for
$0<s<1$, whose eigenpairs are
$(\wh{\lambda}_k^s,\varphi_k)_{k=1}^\infty$. If
$\{\mu_k^{(s)}\}_{k=1}^\infty$ are the eigenvalues of the integral
Laplacian \eqref{eq:defofLaps}, the following equivalence is derived
in \cite{ChenSong05}
\[
C(\Omega) \wh{\lambda}_k^s \le \mu_k^{(s)} \le \wh{\lambda}_k^s, \quad k \in
\mathbb{N}.
\]

There is yet a third family of fractional Sobolev spaces, namely
$H^s_0(\Omega)$, which are the completion of $C_0^\infty(\Omega)$ with the
$L^2$-norm plus the usual $H^s$-seminorm
\begin{equation}\label{eq:Hs-seminorm}
|v|_{H^s(\Omega)}^2 =
C(d,s) \int_\Omega\int_\Omega \frac{|v(x)-v(y)|^2}{|x-y|^{d+2s}} dx
  dy.
\end{equation}
If $\Omega$ is Lipschitz, it turns out that $H^s_0(\Omega) =
\wt{H}^s(\Omega)$ for all $0<s<1$ such that $s\ne\frac12$; in the
latter case $\wt{H}^{\frac12}(\Omega)=H^{\frac12}_{00}(\Omega)$ is the
so-called Lions-Magenes space.  The seminorm \eqref{eq:Hs-seminorm} is
a norm equivalent to $|\cdot|_s$ for $s\in (\frac12,1)$ but not for $s\in
(0,\frac12]$; note that $1 \in H^s_0(\Omega)$ and $|1|_{H^s(\Omega)}=0$
for $s\in (0,\frac12]$. 
%
%This norm equivalence is consistent with the earlier multilevel
%decomposition \eqref{eq:J=infty} of $H^s_0(\Omega)$ derived in
%\cite[Theorem (4.32)]{Xu97} over quasi-uniform grids.
%
%\jp{While the natural norms in the spaces $H^s_0(\Omega)$ and
%$\wt{H}^s(\Omega)$ are equivalent for all $s$, the uniformity of such
%an equivalence constant with respect to $s$ is not straightforward and
%is critical for our analysis.} 
Functions in $H^s_0(\Omega)$ for $s\in
(\frac12,1)$ admit a trace on $\partial\Omega$ and minima of the
functional $v \mapsto \frac12 |v|_{H^s(\Omega)}^2 - \int_\Omega f v$
are weak solutions of the {\it censored} fractional Laplacian. In view
of the norm equivalence $|v|_{H^s_0(\Omega)} \simeq |v|_s$ for $s\in
(\frac12,1)$, the multilevel decomposition \eqref{eq:J=infty} applies to
$H^s_0(\Omega)$ uniformly in $s$ as $s\to 1$ but not as $s\to\frac12$.
This is in agreement with the fact that in the inequality
\[
|v|^2_{\widetilde H^s(\Omega)} \le C |v|^2_{H^s(\Omega)}, \quad v \in
\widetilde H^s(\Omega) = H^s_0(\Omega), \ s > 1/2,
\]
the constant $C$ scales as $(s-1/2)^{-1}$. Indeed, splitting the
integration to compute $\widetilde H^s(\Omega)$ above, one readily
finds that
\[
|v|^2_{\widetilde H^s(\Omega)} = |v|^2_{H^s(\Omega)} + 2 C(d,s)
\int_\Omega \int_{\Omega^c} \frac{|u(x)|^2}{|x-y|^{d+2s}} \, dy dx
\simeq |v|^2_{H^s(\Omega)} + \frac{C(d,s)}{s} \int_\Omega
\frac{|u(x)|^2}{d(x,\partial\Omega)^{2s}} \, dx ,
\]
and it is therefore necessary to bound the last integral in the right
hand side in terms of the $H^s(\Omega)$-seminorm. Such is the purpose
of the {\em Hardy inequality} (cf. \cite[Theorem
1.4.4.4]{grisvard2011elliptic}), for which the optimal constant is of
order $(s-1/2)^{-1}$ \cite{BoDy11}. 

In spite of their spectral equivalence, the inner products that give rise
to the integral, spectral and censored fractional Laplacians are different
and yield a strikingly different boundary behavior \cite{Bonforte18}.
In contrast to
\eqref{eq:boundary-behavior} for the integral Laplacian, for a generic
right-hand side function $f \in L^\infty(\Omega)$ the boundary
behavior of solutions $u$ of the spectral Laplacian is roughly like
\[
u \simeq d(\cdot,\partial\Omega)^{\min\{2s,1\}},
\]
except for $s=\frac12$ that requires an additional factor $|\log
d(\cdot,\partial\Omega)|$, whereas solutions of the censored Laplacian
are quite singular at the boundary \cite{Bonforte18}
\[
u \simeq d(\cdot,\partial\Omega)^{s-\frac12}.
\]
Nevertheless, the above norm equivalences and Theorem
\ref{tm:uniform-cond} imply that the preconditioner $B$ in
\eqref{eq:BPX-FPDE} leads to $\cond \, (BA)$ being bounded
independently of either $s$ and $\bar{J}$ if $A$ is associated to
the spectral Laplacian operator. For the censored Laplacian, the
$\cond \, (BA)$ is uniform with respect to $\bar{J}$ for
$s\in(\frac12,1)$ but blows up as $s\to\frac12$.

%%%%%%%%%%%%%%%%%%%%%%%%%%%%%%%%%%%%%%%%%%%%%%%%%%%%%%%%%%%%%%%%%%%%%%%%%%%%
\section{Graded Bisection grids}\label{S:bisection}
%%%%%%%%%%%%%%%%%%%%%%%%%%%%%%%%%%%%%%%%%%%%%%%%%%%%%%%%%%%%%%%%%%%%%%%%%%%%

This section briefly reviews the bisection method with emphasis on
graded grids, following \cite{chen2012optimal}, and presents new notions.
We also refer to \cite{NoSiVe:2009,NochettoVeeser:12,XuChenNochetto:2009} for
additional details.

%---------------------------------------------------------------------------
\subsection{Bisection rules}
%---------------------------------------------------------------------------
For each simplex $\tau \in \cT$ and a
refinement edge $e$, the pair $(\tau, e)$ is called {\it labeled
simplex}, and $(\cT, \mathcal{L}) := \{(\tau, e): \tau \in \cT\}$ is
called a {\it labeled triangulation}.  For a labeled triangulation
$(\cT, \mathcal{L})$, and $\tau \in \cT$, a {\it bisection} $b_\tau:
\{(\tau, e)\} \mapsto \{(\tau_1, e_1), (\tau_2, e_2)\}$ is a map that
encodes the refinement procedure.  The formal addition is defined as follows: 
$$ 
\cT + b_\tau := (\cT, \mathcal{L}) \setminus \{(\tau, e)\} \cup
\{(\tau_1, e_1), (\tau_2, e_2)\}. 
$$ 
For an ordered sequence of bisections $\mathcal{B} = (b_{\tau_1},
b_{\tau_2}, \ldots, b_{\tau_N})$, we set 
$$ 
\cT + \mathcal{B} := ((\cT + b_{\tau_1}) + b_{\tau_2}) + \cdots +
b_{\tau_N}.
$$ 
Given an initial grid $\mathcal{T}_0$, the set of conforming grids
obtained from $\mathcal{T}_0$ using the bisection method is defined as 
\begin{equation*} \label{eq:bisec-conforming}
\begin{aligned}
\mathbb{T}(\cT_0) := \{\cT = \cT_0 + \mathcal{B} \colon \mathcal{B} \textrm{ is a bisection sequence and $\mathcal{T}$ is conforming} \}. 
\end{aligned}
\end{equation*}
The bisection method considered in this paper satisfies the following two
assumptions: 
\begin{enumerate}
\item[(A1)] Shape regularity: $\mathbb{T}(\cT_0)$ is shape regular. 
\item[(A2)] Conformity of uniform refinement: $\overline{\cT}_k
  := \overline{\cT}_{k-1} + \{b_\tau: \tau\in\overline{\cT}_{k-1}\}
\in\mathbb{T}(\cT_0) ~\forall k \geq 1$.
\end{enumerate}

%----------------------------------------------------------------------------
\subsection{Compatible bisections}
\begin{figure}[!htbp]
\centering 
\subfloat[Interior edge]{\centering 
   \includegraphics[width=0.45\textwidth]{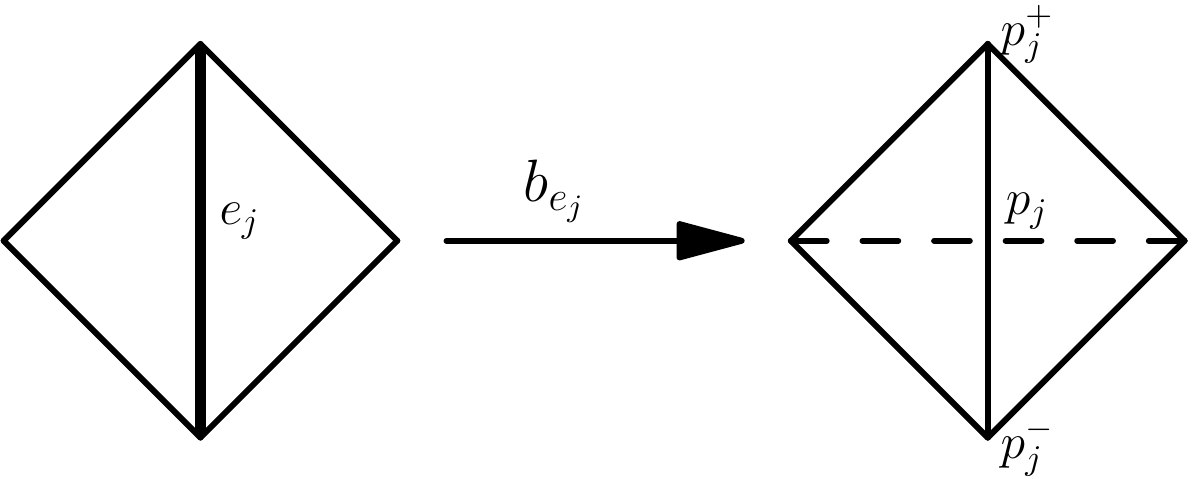} 
}\quad 
\subfloat[Boundary edge]{\centering 
   \includegraphics[width=0.36\textwidth]{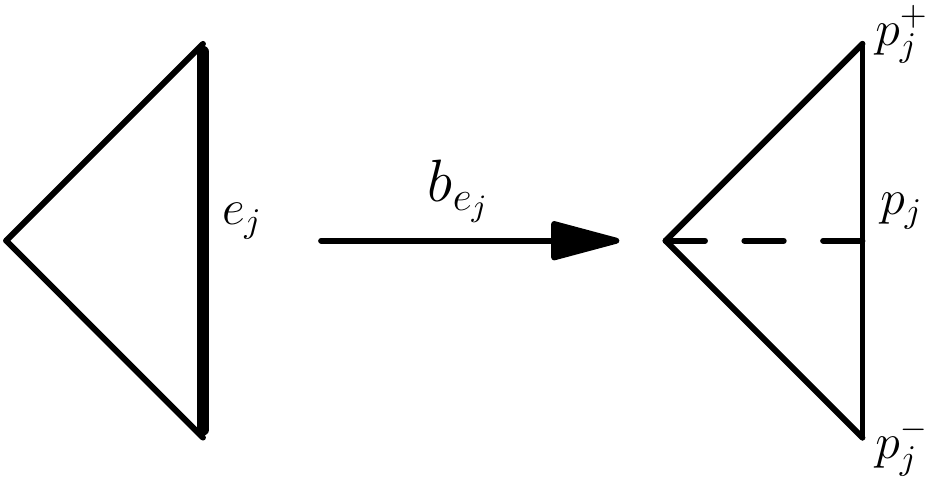} 
}% 
\caption{Two possible configurations of a compatible bisection
$b_{e_j}$ in 2D. The edge with boldface is the compatible refinement
edge, and the dash-line represents the bisection. }
\label{fg:bisec-2D}
\end{figure}

%----------------------------------------------------------------------------
We denote by $\mathcal{N}(\cT)$ the set of vertices of the mesh $\cT$, 
and define the {\it first ring} of either a vertex $p \in \mathcal{N}(\cT)$
or an edge $e \in \mathcal{E}(\cT)$ as 
$$ 
\mathcal{R}_{p} = \{\tau \in \cT~|~ p \in \tau\}, \quad 
\mathcal{R}_e = \{\tau \in \cT~|~ e \subset \tau\},
$$ 
and the {\it local patch} of either $p$ or $e$ as $\omega_{p} = \cup_{\tau\in
\mathcal{R}_{p}} \tau$, and $\omega_e = \cup_{\tau\in \mathcal{R}_e}
\tau$.  An edge $e$ is called {\it compatible} if $e$ is the
refinement edge of $\tau$ for all $\tau \in \mathcal{R}_e$.  Let $p$
be the midpoint of a compatible edge $e$ and $\mathcal{R}_{p}$ be the
ring of $p$ in $\mathcal{T} + \{b_\tau: \tau \in \mathcal{R}_e\}$.
Given a compatible edge $e$, a {\it compatible bisection} is a mapping
$b_e: \mathcal{R}_e \to \mathcal{R}_{p}$. The addition is thus defined
by 
$$ 
\mathcal{T} + b_e := \mathcal{T} + \{b_\tau: \tau \in \mathcal{R}_e\}
= \mathcal{T} \setminus\mathcal{R}_e \cup \mathcal{R}_{p},
$$ 
which preserves the conformity of triangulations. Figure
\ref{fg:bisec-2D} depicts the two possible configurations of a
compatible bisection $b_{e_j}$ in 2D.

We now introduce the concepts of {\it generation} and {\it level}. The
generation $g(\tau)$ of any element $\tau\in\mathcal{T}_0$ is set to be $0$,
and the generation 
of any subsequent element $\tau$ is $1$ plus the generation of 
its father. For any vertex $p$, the generation $g(p)$ of $p$ is defined
as the minimal integer $k$ such that $p \in
\mathcal{N}(\overline{\cT}_k)$. Therefore, $g(\tau)$ and $g(p)$ are the minimal number of
compatible bisections required to create $\tau$ and $p$ from $\mathcal{T}_0$. Once $p$
belongs to a bisection mesh, it will belong to all successive refinements;
hence $g(p)$ is a static quantity insensitive to the level of resolution around
$p$. To account for this issue, we define the {\it level} $\ell(p)$ of a vertex
$p$ to be the maximal generation of elements in the first ring $\mathcal{R}_p$;
this is then a dynamic quantity that characterizes the level of resolution
around $p$.

We then have the decomposition of
bisection grids in terms of compatible bisections; see
\cite[Theorem 3.1]{chen2012optimal}.

\begin{theorem}[decomposition of bisection grids]
\label{tm:bisec-decomp}
Let $\mathcal{T}_0$ be a conforming mesh with initial labeling
that enforces the
bisection method to satisfy assumption (A2), i.e. for all $k \geq 0$
all uniform refinements $\overline{\mathcal{T}}_k$ of $\mathcal{T}_0$
are conforming. Then for every $\mathcal{T} \in \mathbb{T}(\cT_0)$,
there exists a compatible bisection sequence $\mathcal{B} = (b_1, b_2,
\ldots, b_J)$ with $J = \#\mathcal{N}(\cT) - \#\mathcal{N}(\cT_0)$
such that 
\begin{equation} \label{eq:bisec-decomp}
\mathcal{T} = \mathcal{T}_0 + \mathcal{B}.
\end{equation} 
\end{theorem}

\begin{figure}[h!]
\centering
\begin{multicols}{2}
\qquad \includegraphics[width=1.4in]{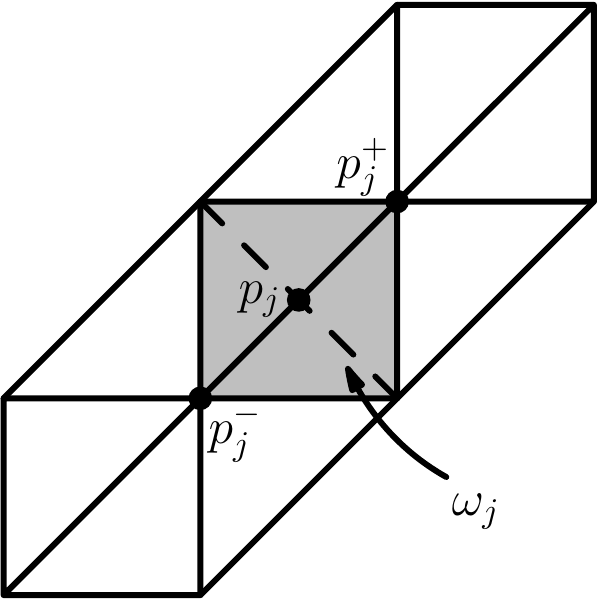} 
\begin{itemize}
\item $e_j$: the refinement edge;
\item $p_j$: the midpoint of $e_j$;
\item $p_j^-, p_j^+$: two end points of $e_j$;
\item $\omega_j$: the patch of $p_j$ (or $\omega_{p_j}$);
\item $\wt{\omega}_j = \omega_{p_j} \cup \omega_{p_j^-} \cup
\omega_{p_j^+}$; 
\item $h_j$: the local mesh size of $\omega_j$;
\item $\mathcal{T}_j = \mathcal{T}_0 + (b_1, \ldots, b_j)$;
\item $\mathcal{R}_j$: the first ring of $p_j$ in $\mathcal{T}_j$.
\end{itemize}
\end{multicols}
\caption{Plot of local patch $\omega_j$ associated to a bisection node $p_j$,
  enlarged local patch $\wt{\omega}_j$ and definition of related quantities.}
\label{fg:bisec-patch}
\end{figure}

For a compatible bisection $b_j$ with refinement edge $e_j$, we
introduce the {\it bisection triplet}
\begin{equation}\label{eq:triplet}
  T_j := \{p_j, p_j^+, p_j^- \},
\end{equation}
where $p_j^-$ and $p_j^+$ are the end points of $e_j$ and $p_j$ is its middle point;
see Figure \ref{fg:bisec-patch}. A vertex can be a middle point of a bisection solely once, when it is created, but instead it can be an end point of a refinement edge repeatedly; in fact this is the mechanism for the level to increment by $1$. In addition, since
$p_j^\pm$ already exist when $p_j$ is created, it follows that
\begin{equation*}\label{eq:gen-middle-end}
g_j :=g(p_j) \ge g(p_j^\pm).
\end{equation*}
The notion of generation of the bisection is well-defined due to the following lemma, see \cite[Lemma 3.3]{chen2012optimal}. 

\begin{lemma}[compatibility and generation] 
\label{lm:compatibility}
If $b_j \in \mathcal{B}$ is a compatible bisection, then all elements
in $\mathcal{R}_j := \mathcal{R}_{p_j}$ have the same generation $g_j$.
\end{lemma}

In light of the previous lemma, we say that $g_j$ is
the generation of the compatible bisection $b_j:\mathcal{R}_{e_j} \to
\mathcal{R}_{p_j}$. Because by assumption
$h(\tau) \simeq 1$ for $\tau \in \mathcal{T}_0$, we have the
following important relation between generation and mesh size:
\begin{equation*} \label{eq:bisec-scale}
h_j \simeq \gamma^{g_j}, \quad \text{with } \gamma = \left(\frac{1}{2}
\right)^{1/d} \in (0,1).  
\end{equation*} 

Moreover, there exists a constant $k_*$ depending on the shape
regularity of $\mathbb{T}(\cT_0)$ such that for every vertex
$p\in\cN(\cT_j)$
\begin{equation}\label{eq:qi-property}
\max_{\tau \in \mathcal{R}_p}g(\tau) - \min_{\tau \in
\mathcal{R}_p}g(\tau) \leq k_*, \qquad \#\mathcal{R}_p \leq k_*.
\end{equation}
Combining this geometric property with Lemma \ref{lm:compatibility}
(compatibility and generation), we deduce that
\begin{equation}\label{eq:gen-patch}
  g_j - k_* \le g(\tau) \le g_j + k_* \quad\forall
  \tau\in\wt{\cR}_j := \cR_{p_j}\cup\cR_{p_j^-}\cup\cR_{p_j^+}.
\end{equation}

Another ingredient for our analysis is the relation between the
generation of compatible bisections and their local or enlarged
patches \cite[Lemmas 3.4 and 3.5]{chen2012optimal}.

\begin{lemma}[generation and patches] \label{lm:nonoverlapping} 
Let $\mathcal{T}_J = \cT_0 + \mathcal{B} \in \mathbb{T}(\cT_0)$ with
compatible bisection sequence $\mathcal{B} = (b_1, \ldots, b_J)$. Then
the following properties are valid:
\begin{enumerate}[$\bullet$]
\item
{\it Nonoverlapping patches:} For any $j \neq k$ and $g_j = g_k$, we have
$$ 
\mathring{\omega}_j \cap \mathring{\omega}_k = \varnothing.
$$

\item
{\it Quasi-monotonicity:} For any $j > i$ and
$\mathring{\widetilde{\omega}}_j \cap \mathring{\widetilde{\omega}}_i
\neq \varnothing$, we have 
$$ 
g_j \geq g_i - 2k_*,
$$ 
where $k_*$ is the integer defined in \eqref{eq:qi-property}.
\end{enumerate}
\end{lemma}

We now investigate the evolution of the level $\ell(p)$ of a generic vertex
$p$ of $\cT$.
\begin{figure}[!htbp]
\centering 
%%\captionsetup{justification=centering}
\subfloat[$\mathcal{R}_q$]{\centering 
   \includegraphics[width=0.25\textwidth]{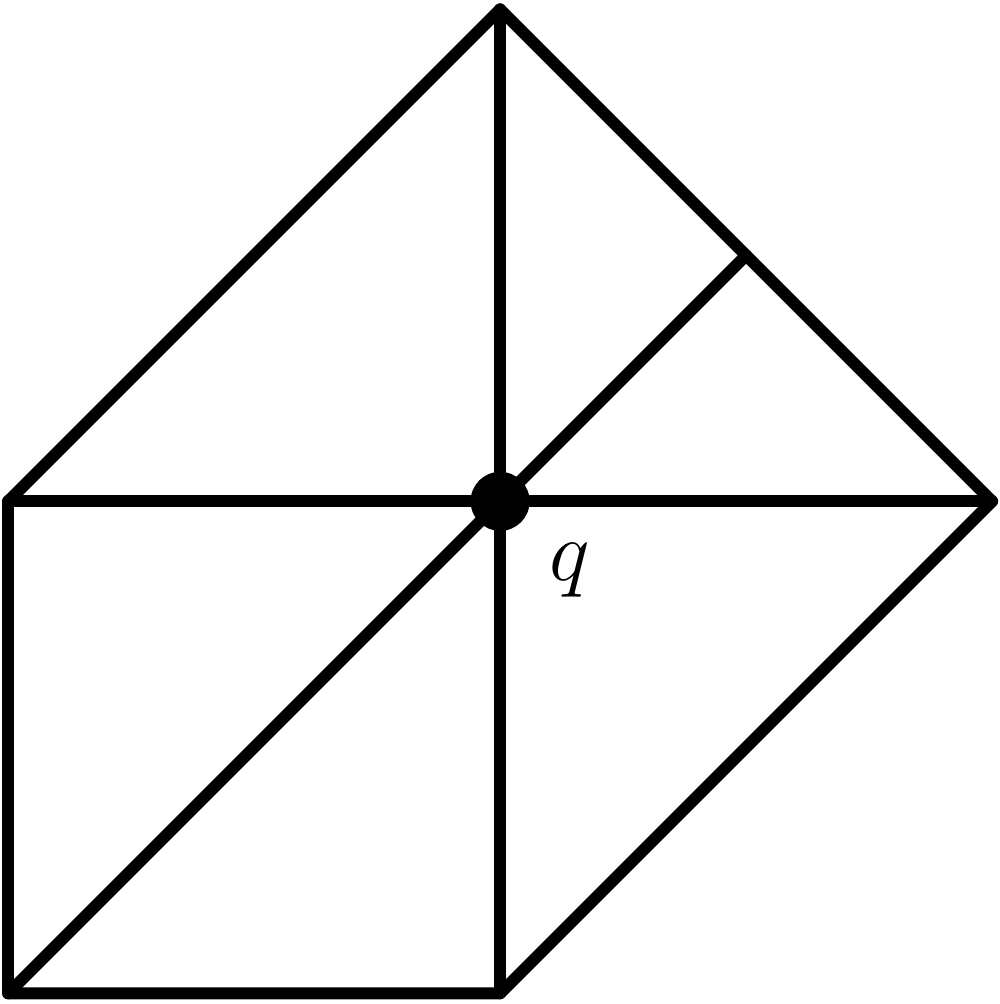} 
   \label{fg:vertex-inclusion1}
} \quad  
\subfloat[Case 1: $q\not\in e$]{\centering 
   \includegraphics[width=0.25\textwidth]{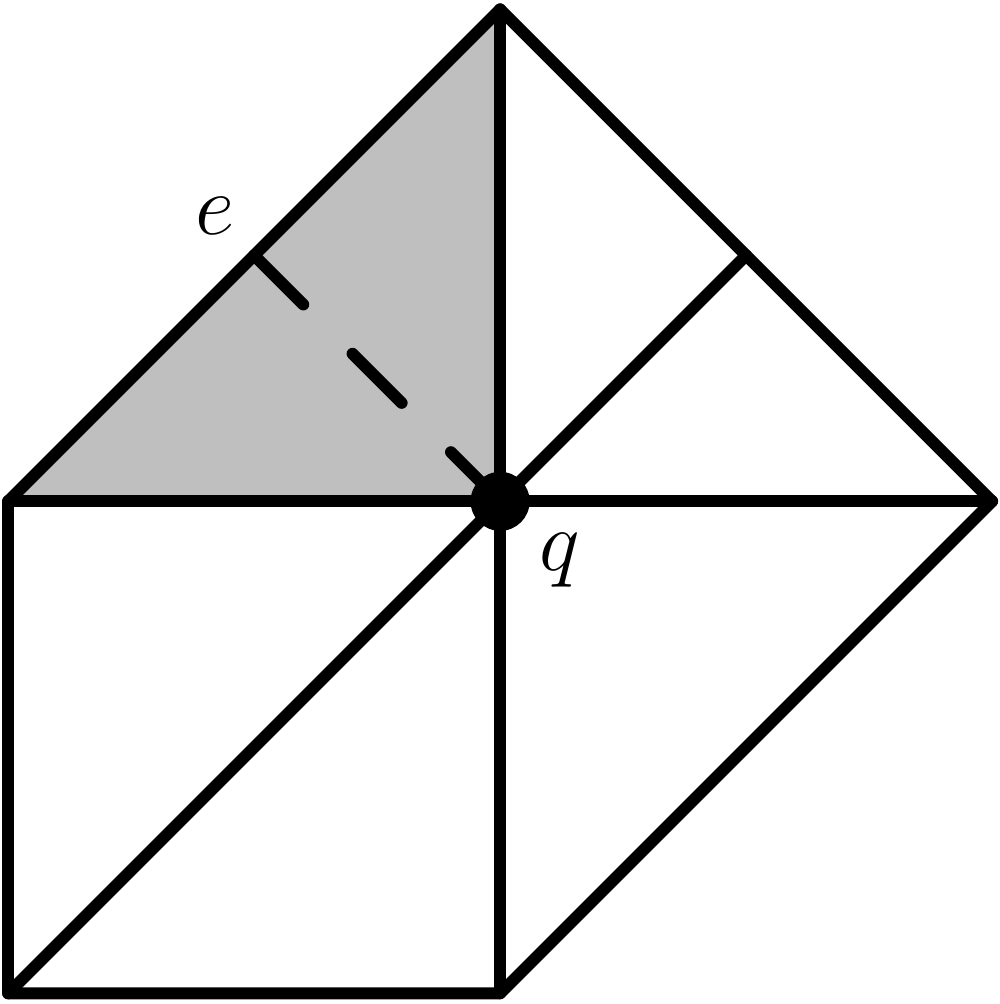} 
   \label{fg:vertex-inclusion2}
} \quad  
\subfloat[Case 2: $q\in e$]{\centering 
   \includegraphics[width=0.25\textwidth]{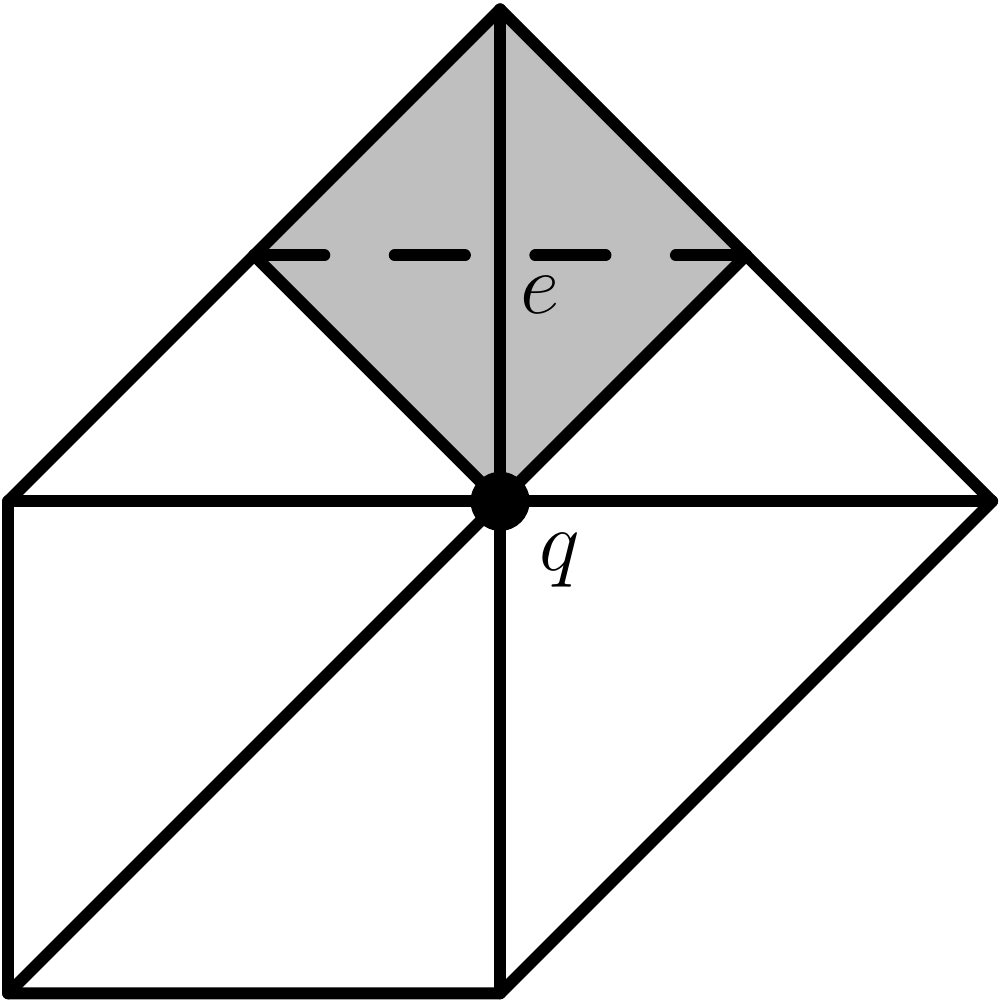} 
   \label{fg:vertex-inclusion3}
} 
\caption{Two cases of bisection in $\mathcal{R}_q$: the bisection edge
  $e$ is on the boundary of the patch and $q$ does not belong to the
  bisection triplet (middle); the node $q$ is an endpoint of $e$ and
  belongs to the bisection triplet (right). The former can happen a
  fixed number $k_*$ of times before the second takes place, where
  $k_*$ depends on the shape regularity of $\mathbb{T}(\cT_0)$.}
\label{fg:vertex-inclusion}
\end{figure}

\begin{lemma}[levels of a vertex]\label{L:level}
If $p\in T_j \cap T_k$, where $T_j$ is a bisection triplet and $T_k$ is the next one
to contain $q$ after $T_j$, and $\ell_j(p)$ and $\ell_k(p)$ are the corresponding
levels, then
\begin{equation*}\label{eq:level}
\ell_k(p) - \ell_j(p) \le k_*
\end{equation*}
where $k_*$ is the integer given in \eqref{eq:qi-property}.
\end{lemma}
\begin{proof}
  Every time a bisection changes the ring $\cR_q$, the level of $q$ may increase
  at most by 1. If the refinement edge $e$ of the bisection is on the boundary of
  the patch $\omega_q$, then $q$ does not belong to the bisection triplet;
  see Figure \ref{fg:vertex-inclusion} (middle). The number of such edges
  is smaller than a fixed integer $k_*$ that only depends on the shape regularity of
  $\mathbb{T}(\cT_0)$. Therefore, after at most $k_*$ bisections the vertex $q$
  is an endpoint of a bisection triplet $T_k$; see Figure \ref{fg:vertex-inclusion}
  (right). This implies $\ell_k(p) \le \ell_j(p)+k_*$ as asserted.
\end{proof}

We conclude this section with the following sequence of auxiliary meshes 
\begin{equation}\label{eq:aux-mesh}
  \wh{\cT}_j := \wh{\cT}_{j-1} + \{b_i\in\cB: g_i = j \} \quad j \ge 1, \qquad
  \wh{\cT}_0 := \cT_0,
\end{equation}
where $\cB$ is the set of compatible bisections \eqref{eq:bisec-decomp}. Note that
each bisection $b_i$ in \eqref{eq:bisec-decomp} does not require
additional refinement beyond the refinement patch $\omega_i$ when incorporated in
the order of the subscript $i$ according to \eqref{eq:bisec-decomp}. This is not
obvious in \eqref{eq:aux-mesh} because the bisections are now ordered by generation.
The mesh $\wh{\cT}_j$ contains all elements $\tau$ of generation $g(\tau)\le j$
leading to the finest graded mesh $\cT=\cT_J$.
The sequence $\{\wh{\cT}_j\}_{j=1}^{\bar{J}}$
is never constructed but is useful for theoretical purposes in Section
\ref{sec:BPX-bisect}.

\begin{lemma}[conformity of $\wh{\cT}_j$] \label{L:conformity-Tj}
The meshes $\wh{\cT}_j$ are conforming for all $j\ge0$.
\end{lemma}
\begin{proof}
We argue by induction. The starting mesh $\wh{\cT}_0$ is conforming by construction.
Suppose that $\wh{\cT}_{j-1}$ is conforming.
We observe that the bisections $b_i$ with $g_i=j$
are disjoint according to Lemma \ref{lm:nonoverlapping} (generation and patches).
Suppose that adding $b_i$ does lead to further refinement beyond the
refinement patch $\omega_i$. If this were the case, then recursive bisection
refinement would end up adding compatible bisections of generation strictly less
than $j$ that belong to the refinement chains emanating from $\omega_i$
\cite{NoSiVe:2009,NochettoVeeser:12}. But such
bisections are all included in $\wh{\cT}_{j-1}$ by virtue of \eqref{eq:aux-mesh}.
This shows that all bisections $b_i$ with $g_i=j$ are compatible with $\wh{\cT}_{j-1}$
and yield local refinements that keep mesh conformity.
\end{proof}

%%%%%%%%%%%%%%%%%%%%%%%%%%%%%%%%%%%%%%%%%%%%%%%%%%%%%%%%%%%%%%%%%%%%%%%%%%%%
\section{Robust BPX preconditioner for graded bisection grids} 
\label{sec:BPX-bisect}
%%%%%%%%%%%%%%%%%%%%%%%%%%%%%%%%%%%%%%%%%%%%%%%%%%%%%%%%%%%%%%%%%%%%%%%%%%%%

In this section, we design and analyze a BPX preconditioner for the
integral fractional Laplacian \eqref{eq:defofLaps} on graded
bisection grids that it is uniform with respect to both number of
levels $J$ and fractional order $s$.  We combine the BPX preconditioner on
quasi-uniform grids of Section \ref{sec:BPX-uniform} with the theory
for graded bisection grids from
\cite{NoSiVe:2009,chen2012optimal,NochettoVeeser:12}, summarized in
Section \ref{S:bisection}, that bridges the gap between graded and
quasi-uniform grids.  Building on Section
\ref{sec:spectral-equivalence}, the results in this section
--especially Theorem \ref{thm:bisection-conditioning}-- apply to the
spectral and censored fractional Laplacians, the latter if
$s>\frac12$ (with a blow up as $s \to \frac{1}{2}$), because of
their spectral equivalence.

%----------------------------------------------------------------------------
\subsection{Space decomposition and BPX preconditioner}
%----------------------------------------------------------------------------
Let $\cT_j = \cT_0 + \{b_1, \cdots,b_j\} \in \mathbb{T}(\cT_0)$
be a conforming bisection grid obtained
from $\cT_0$ after $j\le J$ compatible bisections $\{b_i\}_{i=1}^j$ and 
let $\mathcal{N}_j=\mathring{\mathcal{N}}(\cT_j)$ denote the set of interior
vertices of $\cT_j$.
Let $\mathbb{V}(\cT_j)$ be the finite element space of $C^0$ piecewise
linear functions over $\cT_j$ that vanish on $\partial\Omega$ and its nodal basis
functions be
$\phi_{j,p}$, namely $\mathbb{V}(\cT_j)=\mathrm{span}\{\phi_{j,p}: p\in\cN_j\}$.
We define the {\it local} spaces
\begin{equation} \label{eq:bisec-subspace}
V_j = \mathrm{span}\{\phi_{j,q}: q \in T_j \cap \cN_j\}, \quad j = 1,
\cdots, J
\end{equation}
associated with each bisection triplet $T_j$. We observe that $\mathrm{dim}V_j \leq 3$
and $\mathrm{supp} \, \phi\subset\wt\omega_j$ for $\phi\in V_j$ and $1 \le j \le J$;
see Figure \ref{fg:bisec-patch}. We indicate by $V:= \mathbb{V}(\mathcal{T}_J)$
the finite element space over the finest graded grid $\cT_J$,
with interior nodes $\cP = \cN_J$ and nodal basis functions $\phi_p$
\begin{equation}\label{eq:Vp}
  V = \mathrm{span} \, \{\phi_p: p\in\cP \},
  \quad
  V_p = \mathrm{span} \, \{\phi_p\};
\end{equation}
hence $\mathrm{dim} \ V_p = 1$. Adding the spaces $V_p$
and $V_j$ yields the {\it space decomposition} of $V$
\begin{equation} \label{eq:bisec-decomp-1}
V = \sum_{p\in \cP} V_p +  \sum_{j=0}^J V_j.
\end{equation}
We stress that the spaces $V_j$
appear in the order of creation and not of generation, as is typical of
adaptive procedures. Remarkably, the
functions $\phi_{j,q}$ with $q = p_j^\pm$ depend on the order of creation of
$V_j$ (see Figure \ref{fg:bisec-patch}). Consequently, 
reordering of $V_j$ by generation, which is convenient for analysis,
must be performed with caution; see Sections
\ref{S:boundedness-graded} and \ref{S:stability-graded}.

Let $Q_p$ (resp. $Q_j$) and $I_{p}$ (resp. $I_j$) be the
$L^2$-projection and inclusion operators to and from the discrete spaces $V_p$
(resp. $V_j$), defined in Section \ref{S:space-decomp}.  Inspired by the definition
\eqref{eq:j-smoother2}, we now define the subspace smoothers to be
\begin{equation*} \label{eq:bisec-norm}
\begin{aligned}
R_j v_j & := (1-\wt{\gamma}^s)h_j^{2s} v_j \quad \forall v_j \in V_j,
\\
R_p v_p & := h_p^{2s} v_p ~\qquad\qquad \forall v_p \in V_p, 
\end{aligned}
\end{equation*}
where $R_p$ plays the role of the finest scale whereas $R_j$
represents the intermediate scales. This in turn
induces the following BPX preconditioner on graded bisection grids
\begin{equation}\label{eq:bisec-BPX} 
B = \sum_{p\in \cP} I_p R_p I_p^t + \sum_{j=0}^J I_j
R_j I_j^t = \sum_{p\in \cP} I_p h_p^{2s} Q_p +
(1-\wt{\gamma}^s)\sum_{j=0}^J I_j h_j^{2s} Q_j.
\end{equation}

%---------------------------------------------------------------------  
\subsection{Boundedness: Proof of \eqref{eq:boundedness} for graded bisection grids}
\label{S:boundedness-graded}
%--------------------------------------------------------------------- 

Let $\bar{J} = \max_{\tau \in \mathcal{T}_J} g_\tau$ denote the maximal generation
of elements in $\mathcal{T}_J$. This quantity is useful next to reorder the spaces
$V_j$ by generation because $g_j \le \bar{J}$.

\begin{proposition}[boundedness]\label{P:boundedness-graded}
  Let $v=\sum_{p\in\cP} v_p + \sum_{j=0}^J v_j$ be a decomposition of $v\in V$
  according to \eqref{eq:bisec-decomp-1}. Then, there exists a constant $c_1>0$
  independent of $J$ and $s$ such that
  \begin{equation}\label{eq:boundedness-graded}
    |v|_s^2 \le c_1 \left( \sum_{p \in \cP} h_p^{-2s}\|v_p\|_0^2
    + \frac{1}{1-\wt{\gamma}^s} \sum_{j=0}^J h_j^{-2s}\|v_j\|_0^2 \right),
  \end{equation}
  whence the preconditioner $B$ in \eqref{eq:bisec-BPX} satisfies
  $\lambda_{\max} (BA) \le c_1$.
\end{proposition}
\begin{proof}
We resort to Lemma \ref{lm:inv_ineq} (local inverse inequality) with $\sigma=s$
  and $\mu=0$, which is valid on the graded grid $\cT_J$, to write
\begin{equation} \label{eq:bisec-bd1} 
|v|_s^2 = \bigg|\sum_{p \in \cP} v_p + \sum_{j = 0}^J v_j\bigg|_s^2 \lesssim 
\bigg|\sum_{p \in \cP} v_p\bigg|_s^2 + \bigg|\sum_{j = 0}^J v_j\bigg|_s^2 \lesssim 
\sum_{p \in \cP} h_p^{-2s}\|v_p\|_0^2 + \bigg|\sum_{j = 0}^J v_j\bigg|_s^2.
\end{equation}
In order to deal with the last term, we reorder the functions $v_j$ by
generation and observe that $\mathrm{supp} \, v_j \subset\wt\omega_j$.
We thus define $w_k = \sum_{g_j = k} v_j$ and use
\eqref{eq:gen-patch} to infer that $w_k \in \overline{V}_{k+k_*} =
\mathbb{V}(\overline{\cT}_{k+k_*})$. Similar to the proof of
Proposition \ref{P:boundedness}, using Theorem
\ref{tm:norm-equivalence} (norm equivalence), the fact that $\bar{h}_k
\simeq \gamma^{k}$ and Lemma \ref{lm:s-decomp} ($s$-uniform
decomposition), we have 
\begin{equation*}
\begin{aligned}
\bigg|\sum_{j=0}^J v_j\bigg|_s^2 &=
\bigg|\sum_{k=0}^{\bar{J}}\sum_{g_j = k} v_j\bigg|_s^2
= \bigg|\sum_{k=0}^{\bar{J}} w_k\bigg|_s^2 \\ 
&\simeq \sum_{\ell=0}^{\bar{J}+k_*} \gamma^{-2s\ell}
  \bigg\|(\overline{Q}_{\ell} - \overline{Q}_{\ell-1}) \sum_{k=0}^{\bar{J}}
  w_k\bigg\|_0^2 \\
  & = \inf_{\substack{z_\ell \in \overline{V}_\ell \\ 
  \sum_{\ell=0}^{\bar{J}+k_*} z_\ell = \sum_{k=0}^{\bar{J}} w_k 
  }} \left[ \gamma^{-2s(\bar{J}+k_*)}\|z_{\bar{J}+k_*}\|_0^2 
  + \sum_{\ell = 0}^{\bar{J}+k_*-1}
  \frac{\gamma^{-2s\ell}}{1-\gamma^{2s}}\|z_\ell\|_0^2
  \right].
\end{aligned}
\end{equation*}
Choosing $z_\ell = 0$ for $\ell \le k_*-1$ and
$z_\ell = w_{\ell - k_*} \in \overline{V}_\ell$ for $\ell \ge k_*$ we get
\[
\bigg|\sum_{j=0}^J v_j\bigg|_s^2
  \lesssim \frac{\gamma^{-2sk_*}}{1-\gamma^s}
  \sum_{k=0}^{\bar{J}}\gamma^{-2sk}\|w_k\|_0^2.
\]
In view of Lemma \ref{lm:nonoverlapping}, we see that the enlarged
patches $\wt\omega_j$ and $\wt\omega_i$ have finite overlap depending
only on shape regularity of $\mathbb{T}(\cT_0)$ provided $g_j=g_i$,
whence
\[
\|w_k\|_0^2 \lesssim \sum_{g_j=k} \|v_j\|_0^2.
\]
This in conjunction with $\frac{1-\wt{\gamma}^s}{1-\gamma^s} \simeq 1$
and the fact that $k_*$ is uniformly bounded yields
\begin{equation}\label{eq:bisec-bd2}
\bigg|\sum_{j=0}^J v_j\bigg|_s^2 \lesssim
\frac{1}{1-\wt{\gamma}^s} \sum_{k=0}^{\bar{J}}
\gamma^{-2sk}\sum_{g_j = k}\|v_j\|_0^2 \simeq
\frac{1}{1-\wt{\gamma}^s} \sum_{j=0}^J h_j^{-2s}\|v_j\|_0^2.
\end{equation}
Combining \eqref{eq:bisec-bd1} and \eqref{eq:bisec-bd2} leads to
\eqref{eq:boundedness-graded} as asserted. Finally, the estimate
$\lambda_{\max} (BA) \le c_1$ follows directly from Lemma
\ref{lm:auxiliary} (estimate on $\cond (BA)$).
\end{proof}

%----------------------------------------------------------------------------
\subsection{Stable decomposition: Proof of \eqref{eq:stable-decomp}
for graded bisection grids}\label{S:stability-graded}
%----------------------------------------------------------------------------

We start with a review of the case of quasi-uniform grids in
Corollary \ref{P:stable-decomposition} (stable decomposition) and a
roadmap of our approach. We point out that robustness with respect to
both $J$ and $s$, most notably the handling of factor
$(1-\wt{\gamma}^s)^{-1}$ on coarse levels, is due to the combination
of Lemma \ref{lm:s-decomp} ($s$-uniform decomposition) and Theorem
\ref{tm:norm-equivalence} (norm equivalence), which in turn relies on
Lemma \ref{lm:norm-equiv1} ($s$-uniform interpolation). Since Lemma
\ref{lm:norm-equiv1} fails on graded bisection grids, applying Lemma
\ref{lm:s-decomp} to such grids faces two main difficulties:
(a) Theorem \ref{tm:norm-equivalence} does not hold even for $s=1$;
(b) the spaces $V_j$ in \eqref{eq:bisec-subspace} and $V_p$ in \eqref{eq:Vp}
are locally supported, while the $s$-uniform interpolation requires nested spaces
(see in Lemma \ref{lm:norm-equiv1}). To overcome these
difficulties, we create a family of nested spaces
$\{W_k\}_{k=0}^{\bar{J}}$ with $W_{\bar{J}} = V$ upon grouping indices
according to generation and level around $k$: if
\begin{equation}\label{eq:indices}
\cJ_k:=\{0\le j \le J: \quad g_j\le k\},
\qquad
\cP_k := \{p\in\cP: \quad \ell(p) \le k\},
\end{equation}
then we define $W_k$ to be
\begin{equation} \label{eq:bisect-Wk1}
W_k := \sum_{j \in \cJ_k} V_j + \sum_{p \in \cP_k} V_p.
\end{equation}
Our approach consists of three steps.  The first step, developed in
Section \ref{S:decomp-Wk}, is to derive a global decomposition based
on $W_k$. Since the levels within $W_k$ are only bounded above, to
account for coarse levels we invoke a localization argument
based on a slicing Scott-Zhang operator as in  \cite{chen2012optimal},
which gives the stability result \eqref{eq:interpolation-s1} on
$\{W_k\}_{k=0}^{\bar{J}}$ via Lemma \ref{lm:s-decomp} ($s$-uniform decomposition)
for $s=1$; we bridge the gap to $0<s<1$ via
Lemma \ref{lm:norm-equiv1} ($s$-uniform interpolation). The
space $W_k$ is created for theoretical convenience, but never
constructed in practice, because there is no obvious underlying graded
bisection grid on which the functions of $W_k$ are piecewise linear.
This complicates the stable decomposition of $W_k$ into local spaces
and requires a characterization of $W_k$ in terms of the space
$\wh{V}_k=\mathbb{V}(\wh{\cT}_k)$ of piecewise linear
functions over $\wh{\cT}_k$. The second step in Section \ref{S:charact-Wk}
consists of proving
\[  
\wh{V}_k \subset W_k \subset \wh{V}_{k+k_*},
\]
where $k_*$ is constant. Therefore, the space $W_k$ of unordered
bisections of generation and level $\le k$ is equivalent, up to level $k_*$, to the space 
$\wh{V}_k$ of ordered bisections of generation $\le k$; note that the
individual spaces $V_j$ might not coincide though. In the last step, performed in Section
\ref{S:stable-decomp-graded}, we construct a stable decomposition for
graded bisection grids and associated BPX preconditioner $\wh{B}$. We
also show that $\wh{B}$ is equivalent to $B$ in \eqref{eq:bisec-BPX}.

%----------------------------------------------------------------------------
\subsubsection{Global $L^2$-orthogonal decomposition of
$W_k$}\label{S:decomp-Wk}
%----------------------------------------------------------------------------
%
We recall that the Scott-Zhang quasi-interpolation operator $S_j: V
\to \mathbb{V}(\cT_j)$ can be defined at a node $p\in\cP$ through
the dual basis function on arbitrary elements $\tau \subset \mathcal{R}_p$
\cite{scott1990finite, chen2012optimal}. We exploit this flexibility
to define a suitable quasi-interpolation operator $S_j$ as follows provided
$S_{j-1}: V \to \mathbb{V}(\cT_{j-1})$ is already known. Since
$\cT_j=\cT_{j-1}+b_j$ and the compatible bisection $b_j$ changes $\cT_{j-1}$ locally
in the bisection patch $\omega_{p_j}$ associated with the new vertex $p_j$, we
set $S_j v(p) := S_{j-1}v(p)$ for all $p\in\cN_j\setminus T_j$, where $T_j$ is
the bisection triplet \eqref{eq:triplet}. We next define $S_j v(p_j)$ using a simplex
$\tau\in\cR_j$ newly created by the bisection $b_j$. If $p=p_j^{\pm}\in T_j$
and $\tau \in \cT_{j-1}$ is the simplex used to define
$S_{j-1}v(p)$, then we define $S_j v(p)$ according to the following
rules:
\begin{enumerate}
\item if $\tau \subset \omega_p(\cT_j)$ we keep the nodal value of $S_{j-1}v$,
  i.e. $S_jv(p) = S_{j-1}v(p)$;
\item otherwise we choose a new $\tau \subset \omega_p(\cT_j) \cap \omega_p(\cT_{j-1})$
  to define $S_j v(p)$;
\end{enumerate}
note that $\tau\in\cR_j$ in case (2). Once $\tau\in\cT_j$ has been chosen, then
definition of $S_j v(p)$ for $p\in T_j$ is the same as in \cite{scott1990finite,BrSc07}.
This construction guarantees the local stability bound \cite{scott1990finite}
\begin{equation}\label{equ:scott-zhang3}
h_p^{d/2} |S_j v(p)| \lesssim \|v\|_{\omega_p} \quad\forall p\in\cN_j,
\end{equation}
and that the slicing operator $S_j - S_{j-1}$ is supported in the enlarged
patch $\wt\omega_j$, namely
\begin{equation}\label{equ:scott-zhang2}
(S_j - S_{j-1}) \, v \in V_j \quad \forall 1 \leq j \leq J.
\end{equation}

\begin{lemma}[stable $L^2$-orthogonal decomposition] \label{lm:bisect-bound1}
Let $\wh Q_k:V\to W_k$ be the $L^2$-orthogonal projection operator onto $W_k$ and
$\wh Q_{-1}=0$. For any $v\in V$, the global $L^2$-orthogonal decomposition
$v = \sum_{k=0}^{\bar{J}} (\wh Q_k - \wh Q_{k-1})v$ satisfies
\begin{equation} \label{eq:bisect-bound1}
\sum_{k=0}^{\bar{J}} \gamma^{-2sk} \|(\wh Q_k- \wh Q_{k-1})v\|_0^2 \lesssim 
|v|_s^2, 
\end{equation}
where the hidden constant is independent of $0\leq s \leq 1$
and $\bar{J}$.
\end{lemma}
\begin{proof}
We rely on the auxiliary spaces $\overline{V}_k = \mathbb{V}(\overline{\mathcal T}_k)$
defined over uniformly refined meshes
$\overline{\mathcal T}_k$ of $\mathcal{T}_0$ for $0 \leq k \leq \bar{J}$.
Let $\overline{Q}_k:\overline{V}_{\bar{J}} \to \overline{V}_k$ denote the
$L^2$-orthogonal projection operator onto $\overline{V}_k$ and consider the
global $L^2$-orthogonal decomposition $v = \sum_{k=0} \bar{v}_k$ of any
$v\in V \subset \overline{V}_{\bar{J}}$, where
$\bar{v}_k := (\overline{Q}_k - \overline{Q}_{k-1})v$.
This decomposition is stable in $H^1$
\cite{xu1992iterative, oswald1992norm, bornemann1993basic}
$$ 
\sum_{k=0}^{\bar{J}} \gamma^{-2k} \|\bar{v}_k\|_0^2 \lesssim |v|_1^2.
$$ 

If $g_j$ is the generation of bisection $b_j$ and $g_j>k$,
then $\bar{v}_k$ is piecewise linear in $\omega_{e_j}$ (the patch
of the refinement edge $e_j$), whence $(S_j - S_{j-1})\bar{v}_k = 0$
and the slicing operator detects frequencies $k \ge g_j$.
Consider now the decomposition $v=\sum_{k=0}^{\bar{J}} v_k$ of $v\in V$ where
\begin{equation} \label{eq:decomp-vk-Wk}
  v_k := \sum_{g_j = k} (S_j - S_{j-1}) v = \sum_{g_j =
k} (S_j - S_{j-1})\sum_{\ell=k}^{\bar{J}} \bar{v}_\ell \in W_k.
\end{equation}
In view of Lemma \ref{lm:nonoverlapping} (generation and patches)
and shape regularity of $\mathbb{T}(\cT_0)$,
enlarged patches $\wt\omega_j$ with the same generation $g_j=k$ have
a finite overlapping property. This, in conjunction with \eqref{equ:scott-zhang3} and
\eqref{equ:scott-zhang2} as well as the $L^2$-orthogonality of
$\{\bar{v}_\ell\}_{\ell=k}^{\bar{J}}$, yields
$$ 
\|v_k\|_{0}^2 \lesssim 
\sum_{g_j = k} \Big\|(S_j - S_{j-1})\sum_{\ell=k}^{\bar{J}}\bar{v}_\ell \Big\|_{0,
  \wt{w}_j}^2 \lesssim 
\sum_{g_j = k} \Big\|\sum_{\ell=k}^{\bar{J}}\bar{v}_\ell \Big\|_{0,
  \wt{w}_j}^2 \lesssim 
\Big\|\sum_{\ell=k}^{\bar{J}} \bar{v}_\ell \Big\|_0^2 = 
\sum_{\ell=k}^{\bar{J}} \|\bar{v}_\ell\|_0^2.
$$ 
We use Lemma \ref{lm:s-decomp} ($s$-uniform decomposition) with $s = 1$, together
with \eqref{eq:decomp-vk-Wk}, to obtain
\begin{align*}
\sum_{k=0}^{\bar{J}} \gamma^{-2k}\|(\wh{Q}_k-\wh{Q}_{k-1}) v\|_0^2  & =
  \inf_{\substack{w_k \in W_k \\ \sum_{k=0}^{\bar J} w_k = v}} \Big[
      \gamma^{-2\bar{J}}\|w_{\bar J}\|_0^2 +  \sum_{k = 0}^{\bar{J}-1}
      \frac{\gamma^{-2k}}{1 - \gamma^{2}}\|w_k\|_0^2 \Big]
  \\
    & \leq \gamma^{-2\bar{J}}\|v_{\bar J}\|_0^2 +  \sum_{k
    = 0}^{\bar{J}-1} \frac{\gamma^{-2k}}{1 - \gamma^{2}}\|v_k\|_0^2.
\end{align*}
Employing the preceding estimate of $\|v_k\|_{0}^2$ and reordering the sum
implies
\begin{align*}
\sum_{k=0}^{\bar{J}} \gamma^{-2k}\|(\wh{Q}_k -\wh{Q}_{k-1}) v\|_0^2  &
\lesssim \gamma^{-2\bar{J}}\|\bar{v}_{\bar J}\|_0^2 + \sum_{k = 0}^{\bar{J}-1}
    \frac{\gamma^{-2k}}{1-\gamma^2} \sum_{\ell=k}^{\bar{J}}
    \|\bar{v}_\ell\|_0^2 \\
    & = \gamma^{-2\bar{J}}\|\bar{v}_{\bar J}\|_0^2 + 
    \sum_{\ell=0}^{\bar{J}} \sum_{k=0}^{\ell}
    \frac{\gamma^{-2k}}{1-\gamma^2} 
    \|\bar{v}_\ell\|_0^2 \\
    & = \gamma^{-2\bar{J}}\|\bar{v}_{\bar J}\|_0^2 + 
    \sum_{\ell=0}^{\bar{J}}
    \frac{\gamma^{-2\ell} - \gamma^2}{(1-\gamma^2)^2} 
    \|\bar{v}_\ell\|_0^2 
    \lesssim \sum_{\ell=0}^{\bar J}
    \gamma^{-2\ell}\|\bar{v}_\ell\|_0^2 \lesssim |v|_1^2.
\end{align*}
Hence, we have shown that \eqref{eq:bisect-bound1} holds for $s = 1$.
The desired estimate for arbitrary $0\leq s \leq 1$ follows by
Lemma \ref{lm:norm-equiv1} ($s$-uniform interpolation).
\end{proof}

As a consequence of Lemma \ref{lm:s-decomp} ($s$-uniform decomposition)
and Lemma \ref{lm:bisect-bound1} (stable $L^2$-orthogonal decomposition), 
we deduce the following property.
\begin{corollary}[$s$-uniform decomposition on $W_k$]
\label{co:bisec-decomp-W} 
For every $v \in V$, there exists a decomposition $v = \sum_{k=0}^{\bar{J}} w_k$
with $w_k \in W_k$ for all $k=0,1,\ldots,\bar{J}$ and 
$$
  \gamma^{-2s\bar{J}} \|w_{\bar{J}}\|_0^2 + \sum_{k=0}^{\bar{J}-1}
  \frac{\gamma^{-2sk}}{1-\gamma^{2s}} \|w_{k}\|_0^2 \lesssim |v|_s^2.
$$
\end{corollary}

%----------------------------------------------------------------------------
\subsubsection{Characterization of $W_k$}\label{S:charact-Wk}
%----------------------------------------------------------------------------
We now study the geometric structure of the spaces $W_k$, defined in
\eqref{eq:bisect-Wk1}, which is useful in the construction of a
stable decomposition of $V$. Recalling definition \eqref{eq:aux-mesh},
our first goal is to compare $W_k$ with the space
\[
\wh{V}_k:=\mathbb{V}(\wh{\cT}_k)
\]
of $C^0$ piecewise linear
functions over $\wh{\cT}_{k}$ that have vanishing trace. We will show below
\begin{equation}\label{eq:Vh-Wh}
  \wh{V}_k\subset W_k;
\end{equation}
see Lemmas \ref{lm:nodal-whP} and \ref{lm:nodal-whPc}.
We start with the set of interior
vertices of $W_k$
\begin{equation*} \label{eq:Wk-vertex}
  \cV_k := \cB_k \cup \cP_k, \qquad
  \cB_k := \bigcup  \Big\{T_j: j\in\cJ_k \Big\}, \quad
  \cP_k = \Big\{p\in\cP: \ell(p) \le k \Big\}.
\end{equation*}

\begin{lemma}[geometric structure of $W_k$]\label{L:geom-Wk}
Functions in $W_k$ are $C^0$ piecewise linear on the auxiliary mesh
  $\wh{\cT}_{k+k_*}$, where $k_*$ is given in \eqref{eq:qi-property}. Equivalently,
  $W_k \subset \wh{V}_{k+k_*}$.
\end{lemma}
\begin{proof}
  We examine separately each vertex $q\in\cV_k$. If $q\in\cP_k$, then
  $\ell(q) \le k$ and all elements $\tau\in\cR(q)$ have generation $g(\tau)\le k$
  by definition of level; hence $\tau\in\wh{\cT}_k$ for all $\tau\in\cR(q)$.
  If $q\in\cB_k\setminus\cP_k$
  instead, then the patch of $q$ shares elements with that of the
  bisection node $p_j$
  \[
   \min_{\tau\in\cR_j(q)} g(\tau) \le g(p_j)=g_j \le k,
  \]
  where $\cR_j(q)$ is the ring of elements containing $q$ in the mesh $\cT_j$.
  Property \eqref{eq:qi-property} yields
  \[
   \max_{\tau\in\cR_j(q)} g(\tau) \le \min_{\tau\in\cR_j(q)} g(\tau) + k_* \le k + k_*.
  \]
  It turns out that all elements $\tau\in\wt{\cR}_j$, the enlarged ring around $p_j$,
  have generation $g(\tau)\le k + k_*$, whence $\tau\in \wh{\cT}_{k+k_*}$.
  It remains to realize that any function $w\in V_j$ is thus piecewise linear over
  $\wh{\cT}_{k+k_*}$ and vanishes outside $\wt{\omega}_j$.
\end{proof}

We next exploit the $L^2$-stability of the {\it nodal basis}
$\{\wh{\phi}_q\}_{q\in\wh\cV}$ of $\wh{V}_{k+k_*}$, where $\wh\cV = \wh{\cV}_{k+k_*}$
is the set of interior vertices of $\wh{\cT}=\wh{\cT}_{k+k_*}$.
In fact, if $w = \sum_{q\in\wh{\cV}} w(q) \, \wh{\phi}_q$, then
\begin{equation} \label{eq:nodal-stab}
\|w\|_0^2 = \sum_{\tau \in \wh{\cT}} \|w\|_{0,\tau}^2 \simeq \sum_{\tau \in \wh{\cT}} |\tau|\sum_{q\in \tau} w(q)^2
= \sum_{q \in \wh{\cV}} w(q)^2 \sum_{\tau \ni q}|\tau| \simeq \sum_{q\in \wh{\cV}}
w(q)^2 \|\wh{\phi}_q\|_0^2 .
\end{equation}
Our goal now is to represent each function $\wh{\phi}_q\in\wh{V}_{k+k_*}$ in terms of
functions of $W_{k+k_*}$, which in turn shows $\wh{V}_{k+k_*} \subset W_{k+k_*}$
and thus \eqref{eq:Vh-Wh}.  We start with a partition of $\wh{\cV}_{k+k_*}$,
\begin{equation*} \label{eq:whT-vertex}
  \wh{\cP}_{k+k_*} := \{q \in \wh{\cV}_{k+k_*}: \wh{\ell}(q) \leq k +
  k_*-1\}, \quad  
  \wh{\cP}_{k+k_*}^c := \wh{\cV}_{k+k_*} \setminus \wh{\cP}_{k+k_*},
\end{equation*}
where $\wh{\ell}(q) \le k+k_*$ is the level of $q$ on $\wh{\cT}_{k+k_*}$.
Consequently, $\wh{\ell}(q)= k+k_*$ for all $q \in \wh{\cP}_{k+k_*}^c$ and the
corresponding functions $\wh{\phi}_q$ have all the same scaling due to
shape regularity of $\mathbb{T}(\cT_0)$. In the next two lemmas we represent
the functions $\wh{\phi}_q$ in terms of $W_{k+k_*}$.

\begin{lemma}[nodal basis $\wh{\phi}_q$ with $q\in\wh{\cP}_{k+k_*}$] \label{lm:nodal-whP} 
For any $q \in \wh{\cP}_{k+k_*}$, there holds
\begin{equation*} \label{eq:nodal-whP}
\wh{\phi}_q = \phi_q \quad q\in \cP_{k+k_*-1},
\end{equation*}
where $\cP_k$ is defined in \eqref{eq:indices}; hence, $\wh{\phi}_q\in W_{k+k_*-1}$.
\end{lemma}
\begin{proof}
  Since $\wh{\ell}(q) \leq k+k_*-1$, all elements $\tau \in
  \cR(q)$ have generation $g(\tau) \leq k + k_* -1$. This implies that no further
  bisection is allowed in $\tau$ because all the bisections with
  generation lesser or equal than $k+k_*$ have been incorporated in
  $\wh{\cT}_{k+k_*}$ by definition. Therefore, $\cR(q)$ belongs to the finest
  grid $\cT$ and $\ell(q) = \wh{\ell}(q) \leq k+k_*-1$,
  whence $\wh{\phi}_q\in W_{k+k_*-1}$.
\end{proof}

Next, we consider a nodal basis function $\wh{\phi}_q$ corresponding to
$q \in \wh{\cP}_{k+k_*}^c$. There exists a bisection triplet $T_{j_q}$ that
contains $q$ and $ k \leq \ell_{j_q}(q) \leq k+k_*$, for otherwise
$\ell_{j_q}(q) < k$ would violate Lemma \ref{L:level} (levels of a
vertex).  We thus deduce
\begin{equation} \label{eq:gen-jq} 
k - k_* \leq \ell_{j_q}(q) - k_* \leq g_{j_q} \leq \ell_{j_q}(q) \leq
k + k_*.
\end{equation}
In accordance with \eqref{eq:bisec-subspace}, we denote by $\phi_{j_q, q}$
the nodal basis function of $V_{j_q}$ centered at $q$.  We next show that
$\wh{\phi}_q$ can be obtained by a suitable modification of $\phi_{j_q, q}$
within $W_{k+k_*}$.

\begin{lemma}[nodal basis $\wh{\phi}_q$ with $q\in\wh{\cP}_{k+k_*}^c$]
  \label{lm:nodal-whPc} 
  For any $q \in \wh{\cP}_{k+k_*}^c$, let
  \[
  \mathcal{S}_q := \{j \in \cJ_{k+k_*}: ~ j>j_q, ~ \omega_j \cap
  \supp \, \phi_{j_q,q} \neq \varnothing \}
  \]
  be the set of bisection indices $j>j_q$ such that $g_j \le k+k_*$,
  $\phi_{j,p_j}$ be the function of $V_j$ centered at the bisection vertex $p_j$
  and $\omega_j = \supp \, p_j$. Then there exist numbers $c_{j,q} \in (-1,0]$
  for $j \in \mathcal{S}_q$
such that the nodal basis function $\wh{\phi}_q\in V_{k+k_*}$ associated with $q$
can be written as
  \begin{equation} \label{eq:nodal-whPc}
    \wh{\phi}_q = \phi_{j_q, q} + \sum_{j\in
    \mathcal{S}_q}c_{j,q}\phi_{j,p_j},
  \end{equation}
  and the representation is $L^2$-stable, i.e.,
  \begin{equation} \label{eq:nodal-whPc-stab}
    \|\wh{\phi}_q\|_{0}^2 \simeq \|\phi_{j_q,q}\|_0^2 + \sum_{j\in
    \mathcal{S}_q} c_{j,q}^2 \|\phi_{j,p_j}\|_{0}^2.
  \end{equation}
\end{lemma}
\begin{proof}
  The discussion leading to \eqref{eq:gen-jq} yields
  $k \leq \ell_{j_q}(q) \leq k+k_*$ which, combined with \eqref{eq:qi-property},
  implies that all elements $\tau \in
  \cR_{j_q}(q)$ have generation between $k-k_*$ and $k+k_*$.
  The idea now is to start from the patch $\cR_{j_q}(q)$, the local conforming
  mesh associated with $\phi_{j_q,q}$, and successively refine
  it with compatible bisections in the spirit of the construction of
  $\wh{\cT}_j$ in \eqref{eq:aux-mesh} until we reach the level $k+k_*$;
  see Figure \ref{fg:bisection-basis}.
  To this end, let $\wh{\cT}_{k-k_*}(q) := \cR_{j_q}(q)$ and
  consider the sequence of local auxiliary meshes
  \begin{equation*} \label{eq:local-aux-mesh}
    \wh{\cT}_{j}(q) := \wh{\cT}_{j-1}(q) + \{b_i\in \cB: ~ i \in
    \mathcal{S}_q, ~ g_i = j\} 
  \quad k-k_*+1 \leq j \leq k+k_*,
  \end{equation*} 
  which are conforming according to Lemma \ref{L:conformity-Tj} (conformity
  of $\wh{\cT}_j$).
\begin{figure}[!htbp]
\centering 
  \subfloat[$\wh{\cT}_{k-2,q}$]{\centering 
   \includegraphics[width=0.19\textwidth]{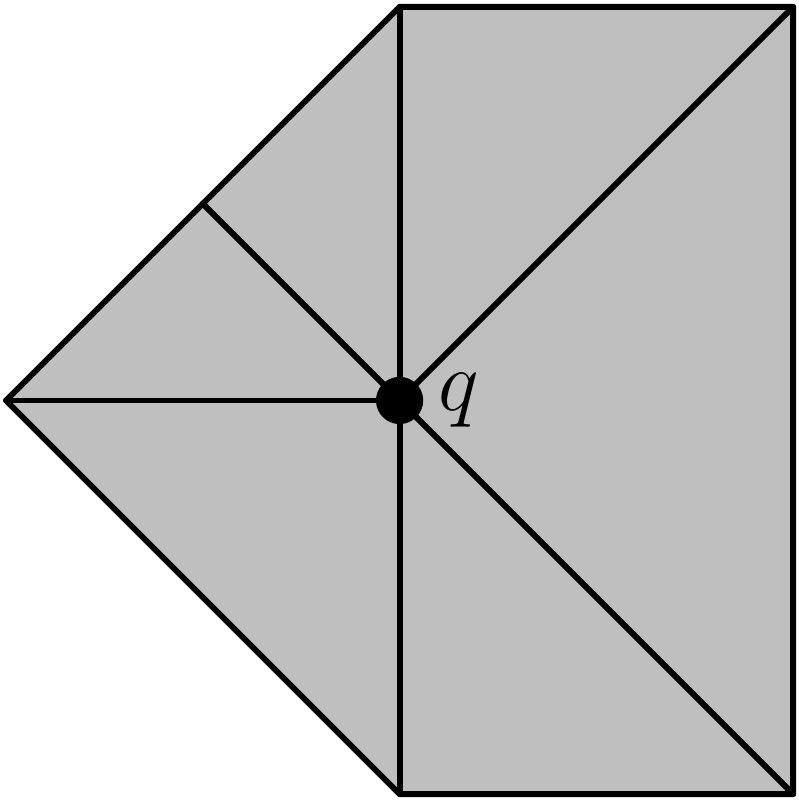} 
   \label{fg:bisection-basis1}
}   
  \subfloat[$\wh{\cT}_{k-1,q}$]{\centering 
   \includegraphics[width=0.19\textwidth]{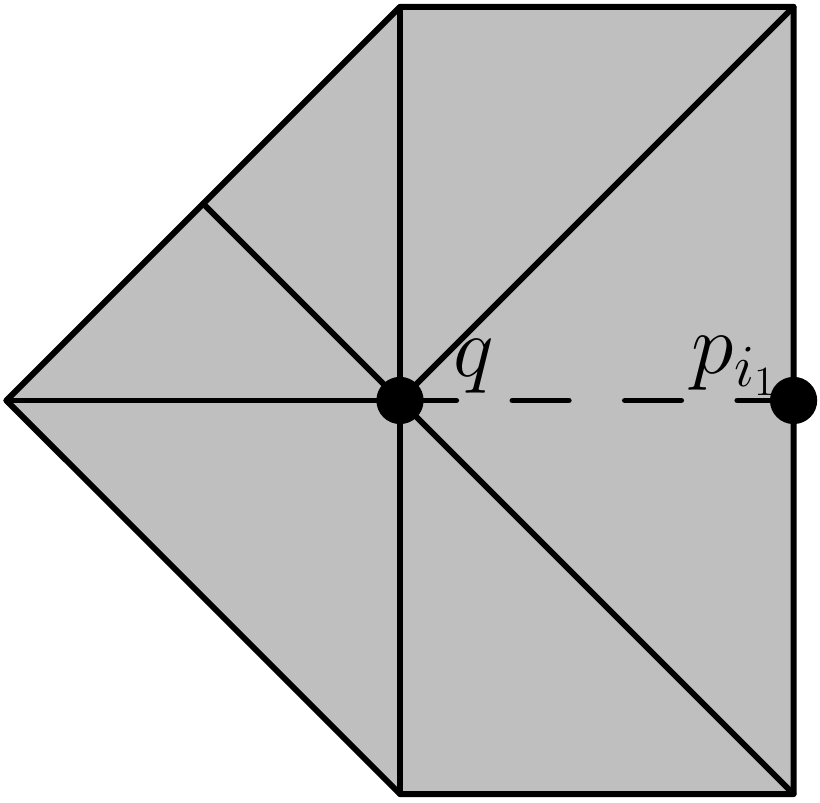} 
   \label{fg:bisection-basis2}
} 
  \subfloat[$\wh{\cT}_{k,q}$]{\centering 
   \includegraphics[width=0.19\textwidth]{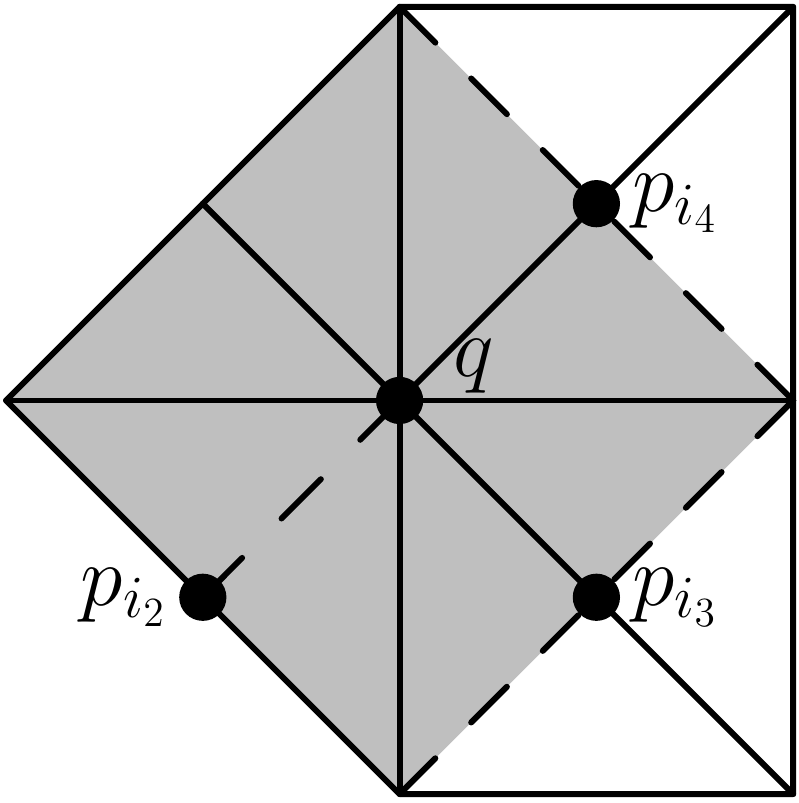} 
   \label{fg:bisection-basis3}
} 
  \subfloat[$\wh{\cT}_{k+1,q}$]{\centering 
   \includegraphics[width=0.19\textwidth]{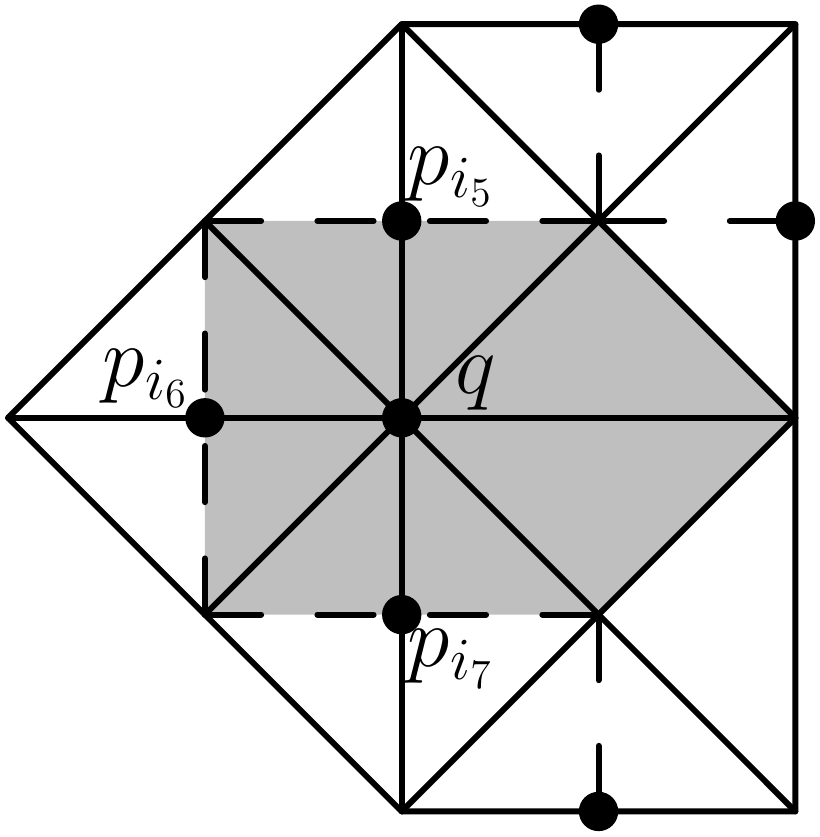} 
   \label{fg:bisection-basis4}
}   
  \subfloat[$\wh{\cT}_{k+2,q}$]{\centering 
   \includegraphics[width=0.19\textwidth]{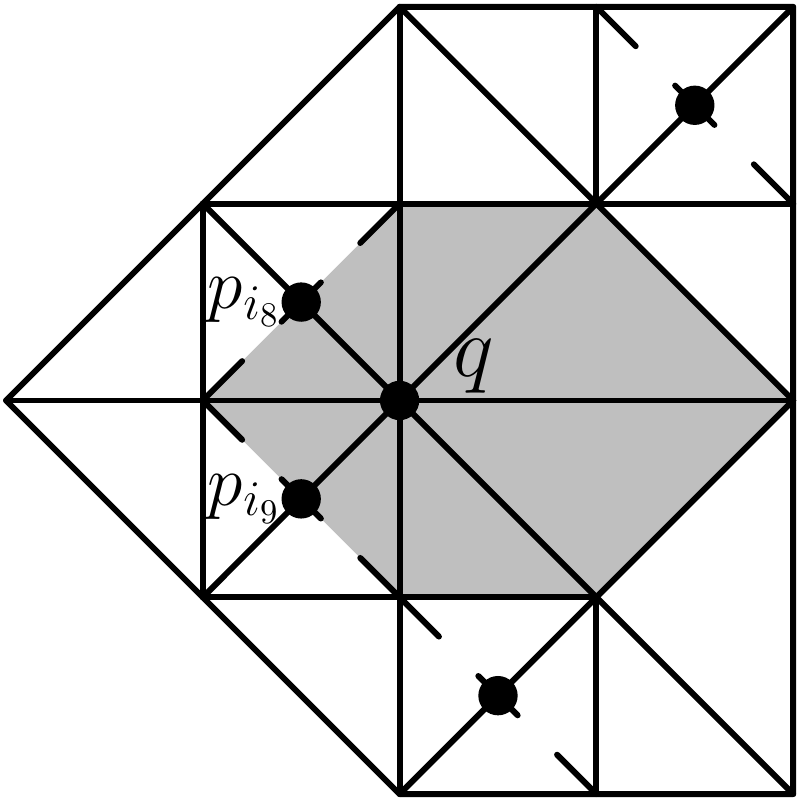} 
   \label{fg:bisection-basis5}
}
  \caption{Local auxiliary meshes $\wh{\cT}_{j,q}$ with $|j-k|\le k_*= 2$. Index sets
  $\mathcal{S}_{k-1,q} = \{i_1\}$, $\mathcal{S}_{k,q} = \{i_2, i_3,
  i_4\}$, $\mathcal{S}_{k+1,q} = \{i_5, i_6, i_7\}$, $\mathcal{S}_{k+2,q}
  = \{i_8, i_9\}$ of compatible bisections to transition from $\wh{\phi}_{j-1,q}$
  to $\wh{\phi}_{j,q}$. The support of $\wh{\phi}_{j,q}$ is monotone
  decreasing as $j$ increases and is plotted in grey.}
\label{fg:bisection-basis}
\end{figure}

We now consider the following recursive procedure: let
  $\wh{\phi}_{k-k_*,q} := \phi_{j_q,q}$ and
  \begin{equation}\label{eq:recursion-phi-j}
    \wh{\phi}_{j,q} := \wh{\phi}_{j-1,q} - \sum_{i \in \mathcal{S}_{j,q}}
    \wh{\phi}_{j-1,q}(p_i) \, \phi_{i,p_i} \quad k-k_*+1 \leq j \leq k+k_*,
  \end{equation}
  where $p_i$ is the bisection node of $b_i \in \cB$ and
  \[
  \mathcal{S}_{j,q} := \Big\{i \in\cJ_{k+k_*}: ~ 
  g_i = j, ~ \omega_i \cap \supp  \, \wh{\phi}_{j-1,q} \neq \varnothing \Big\}.
  \]
  Unless $p_i$ belongs to the boundary of $\supp \, \wh{\phi}_{j-1,q}$,
  the construction \eqref{eq:recursion-phi-j} always modifies $\wh{\phi}_{j-1,q}$;
  compare Figure \ref{fg:bisection-basis2} with Figures
  \ref{fg:bisection-basis3}--\ref{fg:bisection-basis5}. In view of Lemma
  \ref{lm:nonoverlapping} (generation and patches), the sets $\mathring{\omega}_i$
  for $i\in \mathcal{S}_{j,q}$ are disjoint, whence $\wh{\phi}_{j,q}(p)=\delta_{pq}$
  for all nodes $p$ of $\wh{\cT}_{j}(q)$ and $\wh{\phi}_{j,q}$ is the nodal basis
  function centered at $q$ on $\wh{\cT}_j(q)$. Morever,
  \[
  \wh{\phi}_{j,q} = \wh{\phi}_{j-1,q} + \sum_{i\in \mathcal{S}_{j,q}} c_{i,q} \, \phi_{i,p_i}
  \]
  with coefficients $c_{i,q}\in(-1,0]$. The scales of these functions being
  comparable yields
  \[
  \|\wh{\phi}_{j,q}\|_0^2 \simeq \|\wh{\phi}_{j-1,q}\|_0^2 +
  \sum_{i\in \mathcal{S}_{j,q}} c_{i,q}^2 \, \|\phi_{i,p_i}\|_0^2 \, .
  \]
  Since $k_*$ is uniformly bounded depending on shape regularity of $\mathbb{T}(\cT_0)$,
  iterating these two expressions at most $2k_*$ times leads to \eqref{eq:nodal-whPc}
  and \eqref{eq:nodal-whPc-stab}, and concludes the proof.
\end{proof}  

We are now in a position to exploit the representation of nodal basis
of $\wh{V}_{k+k_*}$, given in Lemmas \ref{lm:nodal-whP} and \ref{lm:nodal-whPc}, 
  to decompose functions in $W_k$. We do this next.

\begin{corollary}[$L^2$-stable decomposition of $W_k$] \label{lm:Wk-nodal2}
Given any $0\le k \le \bar{J}$ consider the sets
\begin{equation}\label{eq:indices-2}
  \cP_{k+k_*} = \{q\in\cP: ~ \ell(q) \le k + k_*\},
  \quad
  \cI_{k+k_*} = \{0\le i\le J: ~ k - k_* \le g_i \le k + k_* \}.
\end{equation}
Then, every function $w\in W_k$ admits a $L^2$-stable decomposition
\begin{equation} \label{eq:Wk-decomp}
    w = \sum_{q \in \cP_{k+k_*}} w_q + \sum_{j \in \cI_{k+k_*}} w_j,
    \quad 
    \|w\|_0^2 \simeq \sum_{q \in \cP_{k+k_*}} \|w_q\|_0^2 + \sum_{j
    \in \cI_{k+k_*}} \|w_j\|_0^2 ,
\end{equation}
where $w_q\in V_q$ for all $q \in \cP_{k+k_*}$ and $w_j \in V_j$ for
  all $j \in \cI_{k+k_*}$.
\end{corollary}
\begin{proof}
  Invoking Lemma \ref{L:geom-Wk} (geometric structure of $W_k$), we
  infer that $w\in \wh{V}_{k+k_*}$, which yields the $L^2$-stable
  decomposition of $w$ in terms of nodal basis of $\wh{V}_{k+k_*}$ 
$$ 
  w = \sum_{q \in \wh{\cV}_{k+k_*}} w(q) \, \wh{\phi}_q = \sum_{q \in
  \wh{\cP}_{k+k_*}} w(q) \, \wh{\phi}_q + \sum_{q \in \wh{\cP}_{k+k_*}^c}
  w(q) \, \wh{\phi}_q.
$$ 
  On the one hand, Lemma \ref{lm:nodal-whP} (nodal basis $\wh{\phi}_q$
  with $q\in\wh{\cP}_{k+k_*}$)
  implies that $\wh{\phi}_q = \phi_q$ and $\wh{\cP}_{k+k_*} \subset \cP_{k+k_*}$;
  hence we simply take $w_q := w(q)\phi_q$. On the other hand, using the
  representation \eqref{eq:nodal-whPc} of $\wh{\phi}_q$ from Lemma
  \ref{lm:nodal-whPc} (nodal basis $\wh{\phi}_q$ with $q\in\wh{\cP}_{k+k_*}^c$)
  and reordering, we arrive at
  \begin{align*}
  \sum_{q \in \wh{\cP}_{k+k_*}^c} w(q) \wh{\phi}_q = 
  \sum_{q \in \wh{\cP}_{k+k_*}^c} w(q) \bigg( 
  \phi_{j_q,q} + \sum_{j \in \mathcal{S}_q} c_{j,q} \, \phi_{j, p_j} 
  \bigg)
  = \sum_{j \in \cI_{k+k_*}} w_j,
  \end{align*}
  where
  \[
   w_j :=  \sum_{j_q = j} w(q) \, \phi_{j,q} +
  \sum_{\mathcal{S}_q \ni j} w(q) \, c_{j,q} \, \phi_{j,p_j} \in V_j.
  \]
  This gives the decomposition \eqref{eq:Wk-decomp}. The
  $L^2$-stability \eqref{eq:nodal-stab} of $\{\wh{\phi}_q\}_{q\in\wh{\cV}_{k+k_*}}$
  \[
  \|w\|_0^2 \simeq \sum_{q\in \wh{\cP}_{k+k_*}}
  w(q)^2 \, \|\wh{\phi}_q\|_0^2 + \sum_{q \in \wh{\cP}_{k+k_*}^c} w(q)^2 \,
  \|\wh{\phi}_q\|_0^2,
  \]
  in conjunction with \eqref{eq:nodal-whPc-stab}, gives
  \begin{align*}
  \|w\|_0^2 &
  \simeq \sum_{q\in \wh{\cP}_{k+k_*}} \|w(q)\wh{\phi}_q\|_0^2 + 
  \sum_{q \in \wh{\cP}_{k+k_*}^c} w^2(q) \bigg( %\left[ 
  \|\phi_{j_q,q}\|_0^2 + \sum_{j \in \mathcal{S}_q} c_{j,q}^2
  \|\phi_{j,p_j}\|_0^2 
  \bigg) \\ %\right] \\
  &= \sum_{q \in \cP_{k+k_*}} \|w_q\|_0^2 + \sum_{j \in
  \cI_{k+k_*}} \bigg( %\left[
    \sum_{j_q = j} w^2(q)\|\phi_{j,q}\|_0^2 + \sum_{\mathcal{S}_q \ni
    j} w^2(q) c_{j,q}^2 \|\phi_{j,p_j}\|_0^2
  \bigg ). %\right].
  \end{align*}
  To prove the $L^2$-stability in \eqref{eq:Wk-decomp},
  it remains to show that the term in parenthesis is
  equivalent to $\|w_j\|_0^2$ for any $j \in \cI_{k+k_*}$,
  which in turn is a consequence of the number of
  summands being bounded uniformly. We first observe
  that the cardinality of $\{q: j_q=j\}$ is at most three because this
  corresponds to $q\in T_j$, the $j$-th bisection triplet. Finally, the
  cardinality of the set
  $\{q\in\wh{\cP}_{k+k_*}^c: j \in \mathcal{S}_q\cap\cI_{k+k_*}\}$ is bounded uniformly by a
  constant that depends solely on shape regularity of $\mathbb{T}(\cT_0)$.
  To see this, note that $\wh{\ell}(q)=k+k_*$ yields $k\le g(\tau)\le k+k_*$
  for all elements
  $\tau$ within $\supp \, \phi_{j_q,q}$ and $k-k_* \le g_j \le k+k_*$, whence
  the number of vertices $q$ such that $\supp \,
  \phi_{j_q,q}\cap\omega_j\ne\varnothing$ is uniformly bounded as
  asserted. Hence
  \[
   \|w_j\|_0^2 \simeq \sum_{j_q = j} w^2(q)\|\phi_{j,q}\|_0^2 + \sum_{\mathcal{S}_q \ni
    j} w^2(q) c_{j,q}^2 \|\phi_{j,p_j}\|_0^2
   \]
   yields the norm equivalence in \eqref{eq:Wk-decomp} and finishes
   the proof.
\end{proof}

%----------------------------------------------------------------------------
\subsubsection{Construction of stable decomposition}\label{S:stable-decomp-graded}
%----------------------------------------------------------------------------
We first construct a BPX preconditioner that hinges on the space decomposition
of Section \ref{S:decomp-Wk} and the nodal basis functions just discussed in
Section \ref{S:charact-Wk}. We next show that this preconditioner is
equivalent to \eqref{eq:bisec-BPX}.

\begin{theorem}[stable decomposition on graded bisection grids]
  \label{tm:bisec-decomp}
  For every $v \in V$, there exist $v_p \in V_p$ with $p\in\cP$,
  $v_{p,k} \in V_p$ with
  $p\in\cP_{k+k_*}$, and $v_{j,k} \in V_j$ with $j\in\cI_{k+k_*}$, such that \looseness=-1
  \begin{equation}\label{eq:decomp-graded}
    v = \sum_{p \in \cP} v_p +
    \sum_{k=0}^{\bar{J}}\left(\sum_{q\in\cP_{k+k_*}} v_{q,k} +
    \sum_{j\in\cI_{k+k_*}} v_{j,k} \right),
  \end{equation}
  where $\cP_{k+k_*}$ and $\cI_{k+k_*}$ are given in \eqref{eq:indices-2}, and
  there exists a constant $c_0$ independent of $s$ and $J$ such that
\begin{equation} \label{eq:bisec-decomp-2}
  \gamma^{-2s\bar{J}}\sum_{p \in \cP} \|v_p\|_0^2 +
  \sum_{k=0}^{\bar{J}} \frac{\gamma^{-2sk}}{1-\gamma^{2s}} \left(
  \sum_{p\in\cP_{k+k_*}}\|v_{p,k}\|_0^2 + \sum_{j\in\cI_{k+k_*}}
  \|v_{j,k}\|_0^2 \right)
  \le c_0 |v|_s^2.
\end{equation}
\end{theorem}
\begin{proof}
  We construct the decomposition \eqref{eq:decomp-graded} in three steps.

\smallskip\noindent
    {\it Step 1: Decomposition on $W_k$.} Applying Corollary \ref{co:bisec-decomp-W}
($s$-uniform decomposition on $W_k$), we observe that there
exist $w_k \in W_k$, $k = 0, 1, \cdots, \bar{J}$ such that
$v = \sum_{k=0}^{\bar{J}} w_k$ and
\begin{equation} \label{eq:bisec-decomp-W} 
  \gamma^{-2s\bar{J}} \|w_{\bar{J}}\|_0^2 + \sum_{k=0}^{\bar{J}-1}
  \frac{\gamma^{-2sk}}{1-\gamma^{2s}} \|w_{k}\|_0^2 \lesssim |v|_s^2.
\end{equation}

\smallskip\noindent  
{\it Step 2: Finest scale.} We let $\{\phi_p\}_{p\in\cP}$ be the nodal basis of $V$
and set $v_p := w_{\bar{J}}(p) \phi_p$; hence $w_{\bar{J}}=\sum_{p\in\cP}v_p$. Applying
the $L^2$-stability \eqref{eq:nodal-stab} to $\{\phi_p\}_{p\in\cP}$ gives
\begin{equation} \label{eq:bisec-decomp-J}
\|w_{\bar{J}}\|_0^2 \simeq \sum_{p\in \cP} \|v_p\|_0^2.
\end{equation}
We also choose the finest scale of $v_{p,k}$ and $v_{j,k}$ to be
$v_{q,\bar{J}} = 0$ and $v_{j,\bar{J}} = 0$.

\smallskip\noindent
{\it Step 3: Intermediate scales.} By Corollary \ref{lm:Wk-nodal2}
($L^2$-stable decomposition of $W_k$), we have the $L^2$-stable decomposition
\eqref{eq:Wk-decomp} of $w_k\in W_k$ for every $k = 0 , \ldots , \bar{J}-1$.
Combining the stability bound \eqref{eq:bisec-decomp-W} with 
\eqref{eq:bisec-decomp-J} and \eqref{eq:Wk-decomp}, we
deduce the stable decomposition \eqref{eq:bisec-decomp-2}.
\end{proof}

In view of Theorem \ref{tm:bisec-decomp} above, we consider the BPX
preconditioner
\begin{equation} \label{eq:bisec-nBPX}
\wh{B} := \gamma^{2s\bar{J}} \sum_{p \in \cP} I_p Q_p
  + (1-\gamma^{2s})\sum_{k=0}^{\bar{J}} \gamma^{2sk} \left( 
  \sum_{p\in\cP_{k+k_*}} I_p Q_p +
  \sum_{j\in\cI_{k+k_*}} I_j Q_j
  \right).
\end{equation}

The following corollary is a direct consequence of \eqref{eq:bisec-decomp-2}
and \eqref{eq:stable-decomp}.

\begin{corollary}[uniform bound for $\lambda_{\min} (\wh{B}A)$]
\label{C:lambda-min-BhA}
The preconditioner $\wh{B}$ in \eqref{eq:bisec-nBPX} satisfies
\[
\lambda_{\min} (\wh{B} A ) \ge c_0^{-1}
\]
\end{corollary}

We are now ready to prove the main result of this section, namely that $B$
in \eqref{eq:bisec-BPX} is a robust preconditioner for $A$ on graded bisection
grids. To this end, we need to show that $\wh{B}$ in \eqref{eq:bisec-nBPX} is
spectrally equivalent to $B$.

\begin{theorem}[uniform preconditioning on graded bisection grids]
\label{thm:bisection-conditioning}
Let $\Omega$ be a bounded Lipschitz domain and $s \in (0,1)$.
Let $V$ be the space of continuous piecewise linear finite elements
over a graded bisection grid $\cT$, and consider the space decomposition
\eqref{eq:bisec-decomp-1}. The corresponding BPX preconditioner $B$ in
\eqref{eq:bisec-BPX}, namely
\begin{equation*}\label{eq:graded-BPX}
B = \sum_{p\in \cP} I_p h_p^{2s} Q_p +
(1-\wt{\gamma}^s)\sum_{j=0}^J I_j h_j^{2s} Q_j,
\end{equation*}
is spectrally equivalent to $\wh{B}$ in \eqref{eq:bisec-nBPX},
whence $\lambda_{\min} (BA) \gtrsim c_0^{-1}$. Therefore, the
condition number of $BA$ satisfies
\[
\cond \, (BA) \lesssim c_0 c_1,
\]
where the constants $c_0$ and $c_1$ are independent of $s$ and $J$
and given in \eqref{eq:bisec-decomp-2} and \eqref{eq:boundedness-graded}.
\end{theorem}
\begin{proof}
We show that the ratio $\frac{(Bv,v)}{(\wh{B}v,v)}$ is bounded below and above
by constants independent of $s$ and $J$ for all $v\in V$. We first observe
that for $p\in\cP$ with level $\ell(p)$, we have
\[
h_p^{2s} \simeq \gamma^{2s\ell(p)} = 
\gamma^{2s(\bar{J}+1)} + (1-\gamma^{2s}) \sum_{k=\ell(p)}^{\bar{J}}
\gamma^{2sk},
\]
whence $B_1 :=  \sum_{p\in \cP} I_p h_p^{2s} Q_p$ and $v_p
= Q_p v$ satisfy
\begin{equation*}
    (B_1 v,v) \simeq \gamma^{2s\bar{J}} \sum_{p \in \cP} \|v_p\|_0^2
    + (1-\gamma^{2s}) \sum_{p \in \cP} \sum_{k=\ell(p)}^{\bar{J}} \gamma^{2sk} \|v_p\|_0^2.
\end{equation*}
The rightmost sum can be further decomposed as follows:
\begin{align*}
  \sum_{p \in \cP} \sum_{k=\ell(p)}^{\bar{J}} \gamma^{2sk} \|v_p\|_0^2 &=
  \sum_{j=0}^{\bar{J}} \sum_{\ell(p) = j} \sum_{k =
    j}^{\bar{J}} \gamma^{2sk} \|v_p\|_0^2 \\
  & = \sum_{k=0}^{\bar{J}} \gamma^{2sk} \sum_{\ell(p) \leq k} \|v_p\|_0^2
  \le \sum_{k=0}^{\bar{J}} \gamma^{2sk} \sum_{\ell(p) \leq k+k_*} \|v_p\|_0^2 \\
  &= \gamma^{-2sk_*} \sum_{k=0}^{\bar{J}} \gamma^{2s(k+k_*)}
  \sum_{\ell(p) \leq k+k_*} \|v_p\|_0^2 \le \gamma^{-2sk_*}
  \sum_{k=0}^{\bar{J}} \gamma^{2sk} \sum_{\ell(p) \leq k} \|v_p\|_0^2.
\end{align*}
Since $\gamma^{-2sk_*}\simeq1$, there exist equivalence constants independent of $s$ and $J$ such that
\begin{equation} \label{eq:bisec-nBPX1} 
(B_1 v,v) \simeq \gamma^{2s\bar{J}} \sum_{p \in \cP} \|v_p\|_0^2
    + (1-\gamma^{2s}) \sum_{k=0}^{\bar{J}} \gamma^{2sk} \sum_{p\in\cP_{k+k_*}}\|v_p\|_0^2.
\end{equation}

We now consider the bisection triplets $T_j$ and $3$-dimensional spaces $V_j$,
  for which $h_j \simeq \gamma^{g_j}$. We let $\wh{B}_2 :=
  \sum_{k=0}^{\bar{J}} \gamma^{2sk} \sum_{j\in\cI_{k+k_*}} I_j Q_j$,
  $B_2 := \sum_{j=0}^J I_j h_j^{2s}Q_j$ and $v_j := Q_j v$, to write
\begin{equation}\label{eq:bisec-nBPX2} 
\begin{aligned}
  (\wh{B}_2 v,v) &= \sum_{k=0}^{\bar{J}} \gamma^{2sk}
  \sum_{k-k_* \le g_j \le k+k_* } \|v_j\|_0^2
  \\
  & = \sum_{k=0}^{\bar{J}} \gamma^{2sk}
  \sum_{i=-k_*}^{i=k_*} \gamma^{2si} \sum_{g_j=k}\|v_j\|_0^2
  \\
  & \simeq \sum_{k=0}^{\bar{J}} \gamma^{2sk} \sum_{g_j=k}\|v_j\|_0^2
  = \sum_{j=1}^J \gamma^{2s g_j} \|v_j\|_0^2 \simeq (B_2 v,v),
\end{aligned}
\end{equation}
because $\sum_{i=-k_*}^{i=k_*} \gamma^{2si} \simeq 1$ due to the fact
that $k_*$ is a fixed integer depending solely on shape regularity
of $\mathbb{T}(\cT_0)$.
Combining \eqref{eq:bisec-nBPX1} and \eqref{eq:bisec-nBPX2} we
obtain
\[
  (Bv,v) = (B_1v,v) + (1-\wt{\gamma}^{s}) (B_2v,v) \simeq (\wh{B} v,v)
\quad\forall \, v\in V,
\]
whence the operators $B$ and $\wh{B}$ are spectrally equivalent. Invoking
Corollary \ref{C:lambda-min-BhA} (uniform bound for $\lambda_{\min} (\wh{B}A)$),
we readily deduce $\lambda_{\min} (BA) \gtrsim c_0^{-1}$. We finally
recall that $\lambda_{\max} (BA) \lesssim c_1$, according to
Proposition \ref{P:boundedness-graded} (boundedness), to infer the
desired uniform bound $\cond (BA) = \lambda_{\max} (BA)
\lambda_{\min} (BA)^{-1}\lesssim c_0 c_1$.
\end{proof}

%%%%%%%%%%%%%%%%%%%%%%%%%%%%%%%%%%%%%%%%%%%%%%%%%%%%%%%%%%%%%%%%%%%%%%%%%
\section{Numerical Experiments} \label{sec:experiments}
%%%%%%%%%%%%%%%%%%%%%%%%%%%%%%%%%%%%%%%%%%%%%%%%%%%%%%%%%%%%%%%%%%%%%%%%%
This section presents some experiments in both uniform and graded
bisection grids. We provide some details about the implementation of
the BPX preconditioners and their matrix representations in Appendix
\ref{sec:matrix-representation}.

In the sequel, we solve
\eqref{eq:Dirichlet} with $\Omega = (-1,1)^2$ and $f = 1$, and $s =
0.9$, $s=0.5$ or $s = 0.1$.  In all numerical experiments, the
stopping criterion is
\[
\frac{\|b- {\bf A}x \|_2}{\|b\|_2} \leq
1\times 10^{-6}.
\]

%------------------------------------------------------------------------
\subsection{Uniform grids}
%------------------------------------------------------------------------

In first place, we perform computations on a family of nested,
uniformly refined meshes.
Table \ref{tab:FL-BPX} lists the number of iterations performed when
solving the linear systems using either the Gauss-Seidel (GS), Conjugate Gradient (CG) and Preconditioned Conjugate Gradient (PCG) methods. 
Limited by computational capacity, the largest $\bar{J}$ we take in
our computations is $6$, which corresponds to $16129$ degrees of freedom
(DOFs).  Even though this is a small-scale problem, the BPX
preconditioner \eqref{eq:BPX-FPDE-M} performs well. 
\begin{table}[!ht]
\centering
\begin{tabular}{|m{0.3cm}|m{0.5cm}|m{1.0cm}||m{1cm}|m{0.75cm}|
m{0.75cm}||m{1cm}|m{0.75cm}|m{0.75cm}||m{1cm}|m{0.75cm}|m{0.75cm}|}
  \hline
  \multirow{2}{*}{$\bar{J}$} & \multirow{2}{*}{$h_{\bar{J}}$} &
  \multirow{2}{*}{DOFs} 
  & \multicolumn{3}{c||}{$s=0.9$}
  & \multicolumn{3}{c||}{$s=0.5$} 
  & \multicolumn{3}{c|}{$s=0.1$} \\ \cline{4-12}
  && & GS & CG  & PCG  & GS & CG & PCG & GS & CG  & PCG 
  \\ \hline
  1 & $2^{-1}$ & 9 & 18 & 4 & 4 
                     & 8 & 4 & 4 
                     & 7 & 4 & 4 \\ \hline
  2 & $2^{-2}$ & 49 & 64 & 12 & 12 
                     & 16 & 8 & 8 
                     & 7 & 8 & 9 \\ \hline
  3 & $2^{-3}$ & 225 & 222 & 25 & 16 
                      & 33 & 11 & 10 
                     & 8 & 8 & 10 \\ \hline
  4 & $2^{-4}$ & 961 & 772 & 46 & 19 
                     & 68 & 17 & 11 
                     & 8 & 8 & 10 \\ \hline
  5 & $2^{-5}$ & 3969 & 2689 & 84 & 21 
                     & 139 & 24 & 12 
                     & 9 & 8 & 10 \\ \hline
  6 & $2^{-6}$ & 16129 & 9363 & 157 & 22 
                     & 279 & 32 & 13
                     & 9 & 8 & 10 \\ \hline
\end{tabular}
\medskip
\caption{Number of iterations: GS, CG, and PCG with BPX
  preconditioner \eqref{eq:BPX-FPDE-M}, $\tilde{\gamma} = 0.5$.}
\label{tab:FL-BPX}
\end{table}

%------------------------------------------------------------------------
\subsection{Graded bisection grids} \label{sec:experiments-graded}
%------------------------------------------------------------------------
We next consider graded bisection grids.  As described in Proposition
\ref{prop:bdry-regularity} (regularity in weighted spaces) and Remark
\ref{rm:optimal-parameters}, the solution $u$ to \eqref{eq:Dirichlet}
satisfies $u \in \cap_{\eps > 0} \, \widetilde
H^{1+s-2\eps}_{\frac12-\eps}(\Omega)$ and this regularity can be optimally
exploited by considering grids graded according to \eqref{eq:grading}
with $\mu = 2$. In the energy norm, one obtains linear convergence
rates with this strategy.

\begin{figure}[!htbp]
\centering 
%%\captionsetup{justification=centering}
\subfloat[$\bar{J}=6$]{\centering 
   \includegraphics[width=0.28\textwidth]{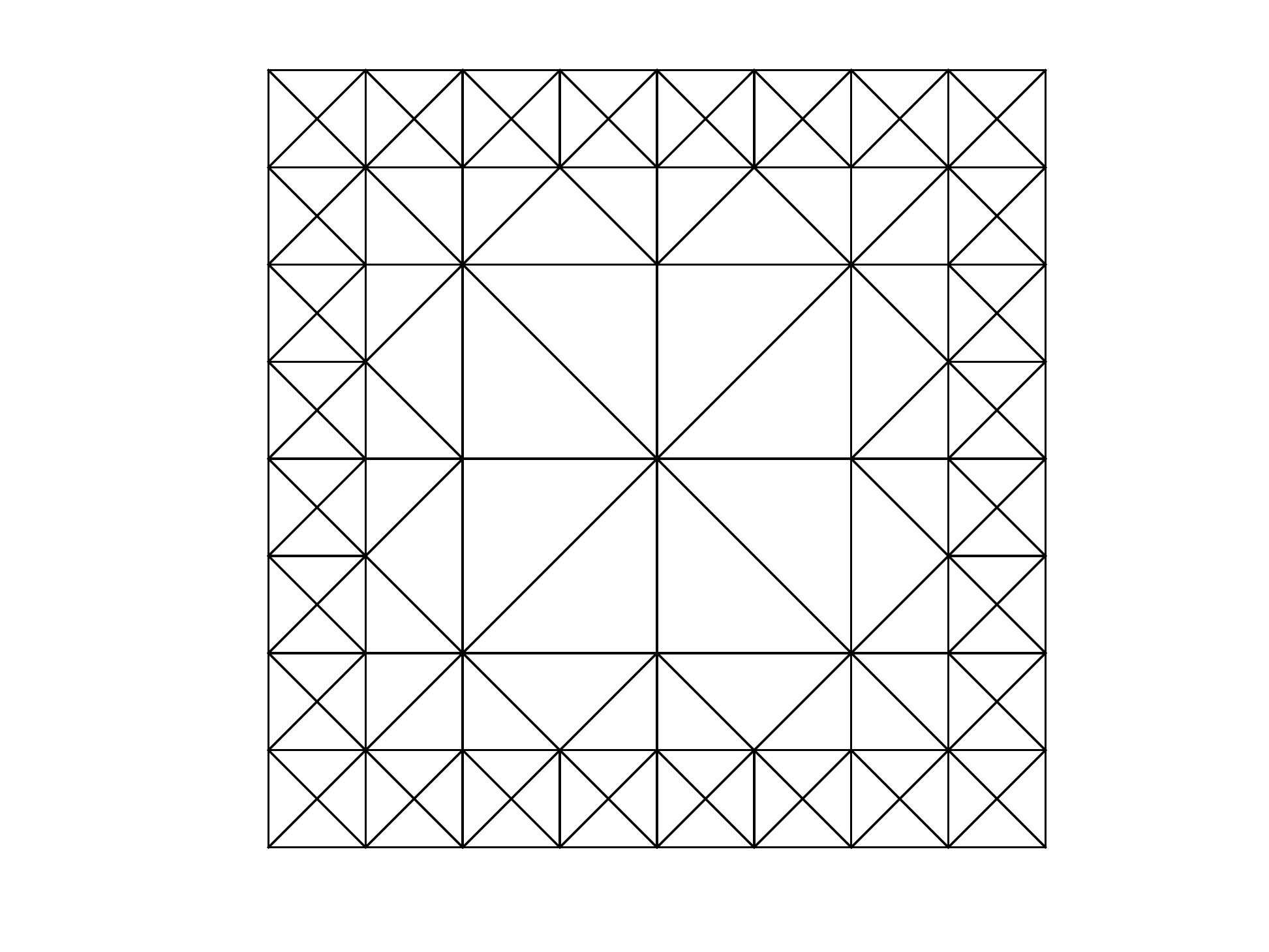} 
   \label{fig:J6}
} \hspace{-0.9cm} 
\subfloat[$\bar{J}=9$]{\centering 
   \includegraphics[width=0.28\textwidth]{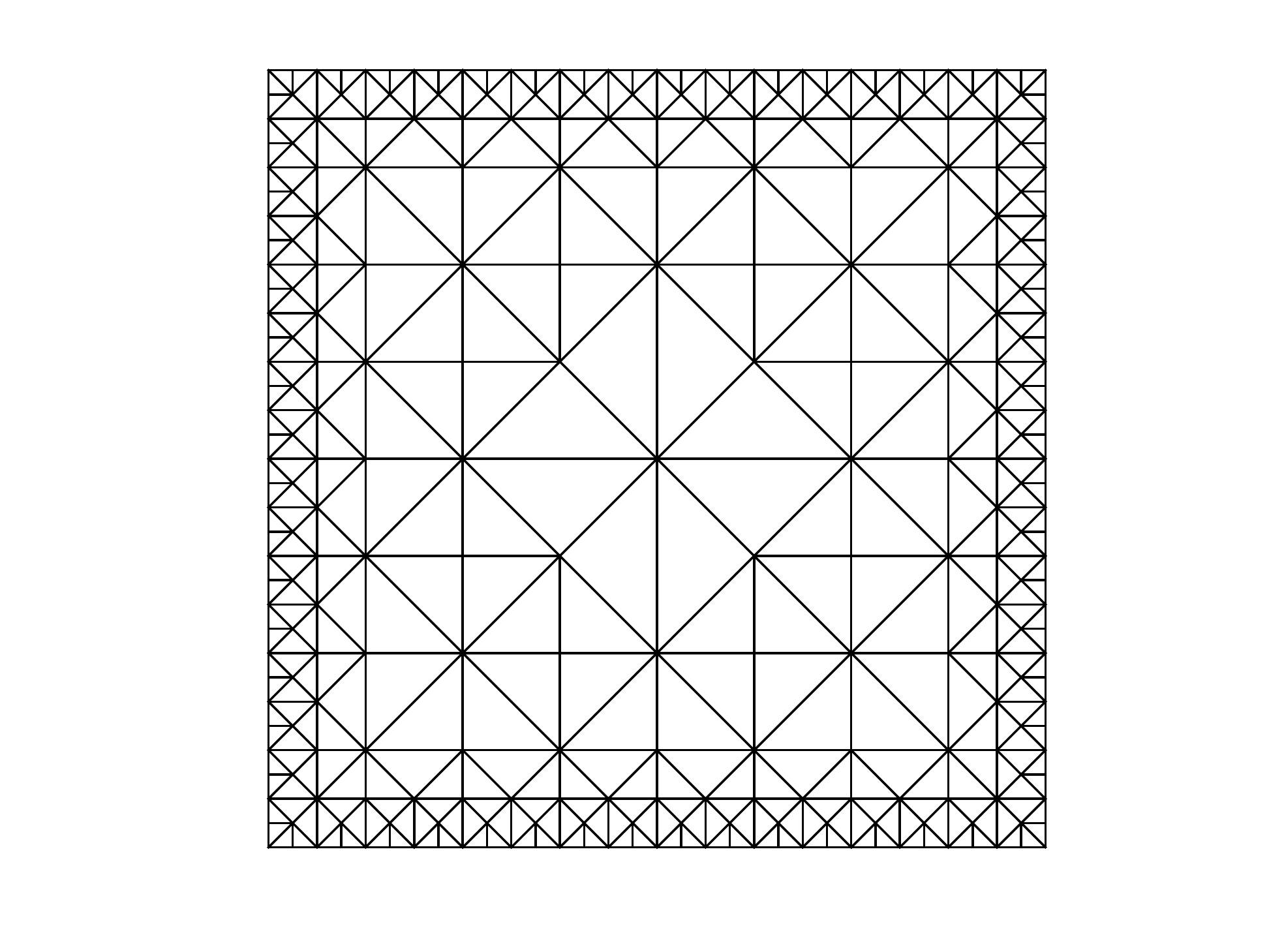} 
   \label{fig:J9}
}% 
\hspace{-0.9cm}
\subfloat[$\bar{J}=12$]{\centering 
   \includegraphics[width=0.28\textwidth]{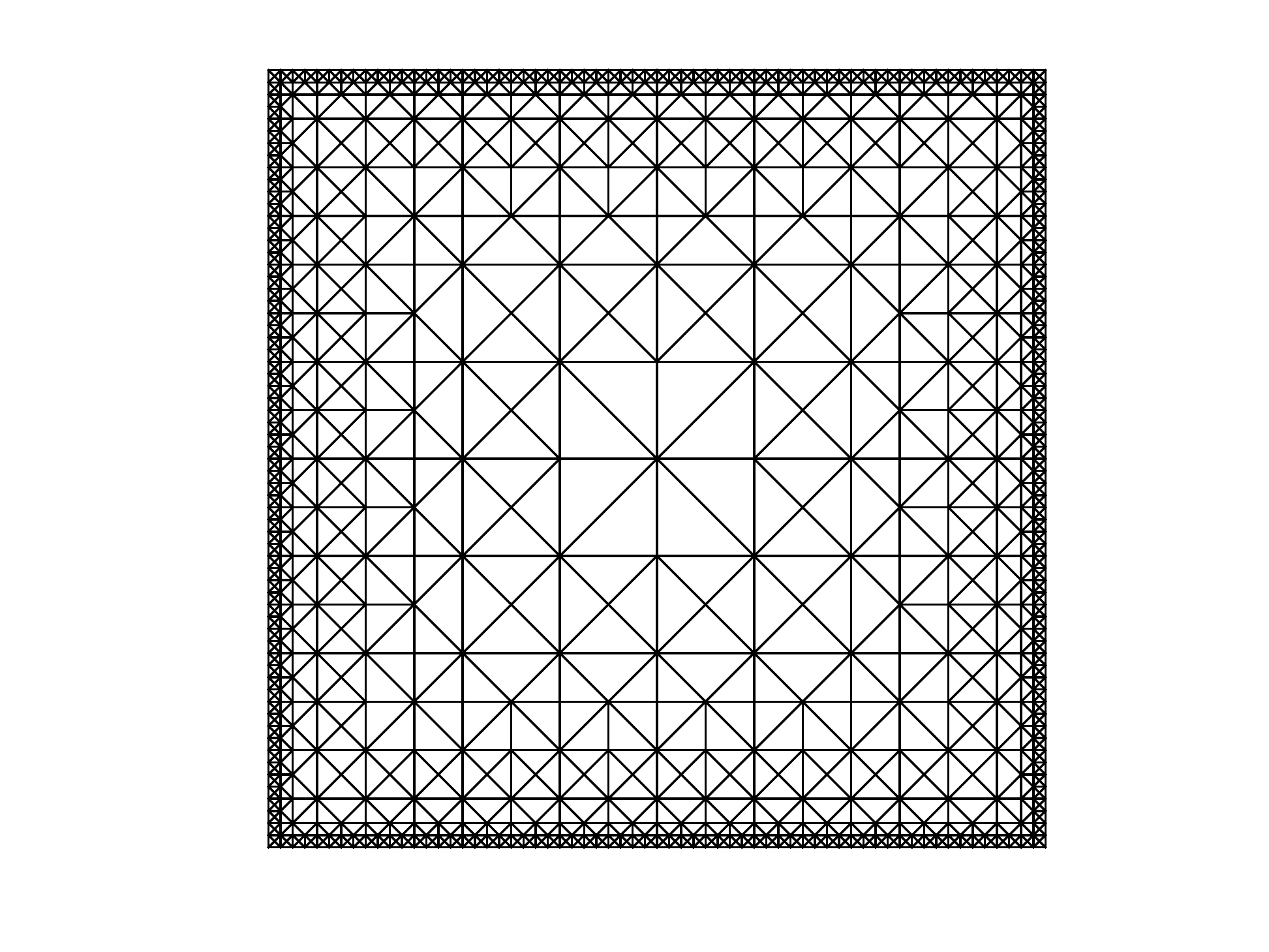} 
   \label{fig:J12}
}\hspace{-0.9cm} 
\subfloat[$\bar{J}=15$]{\centering 
   \includegraphics[width=0.28\textwidth]{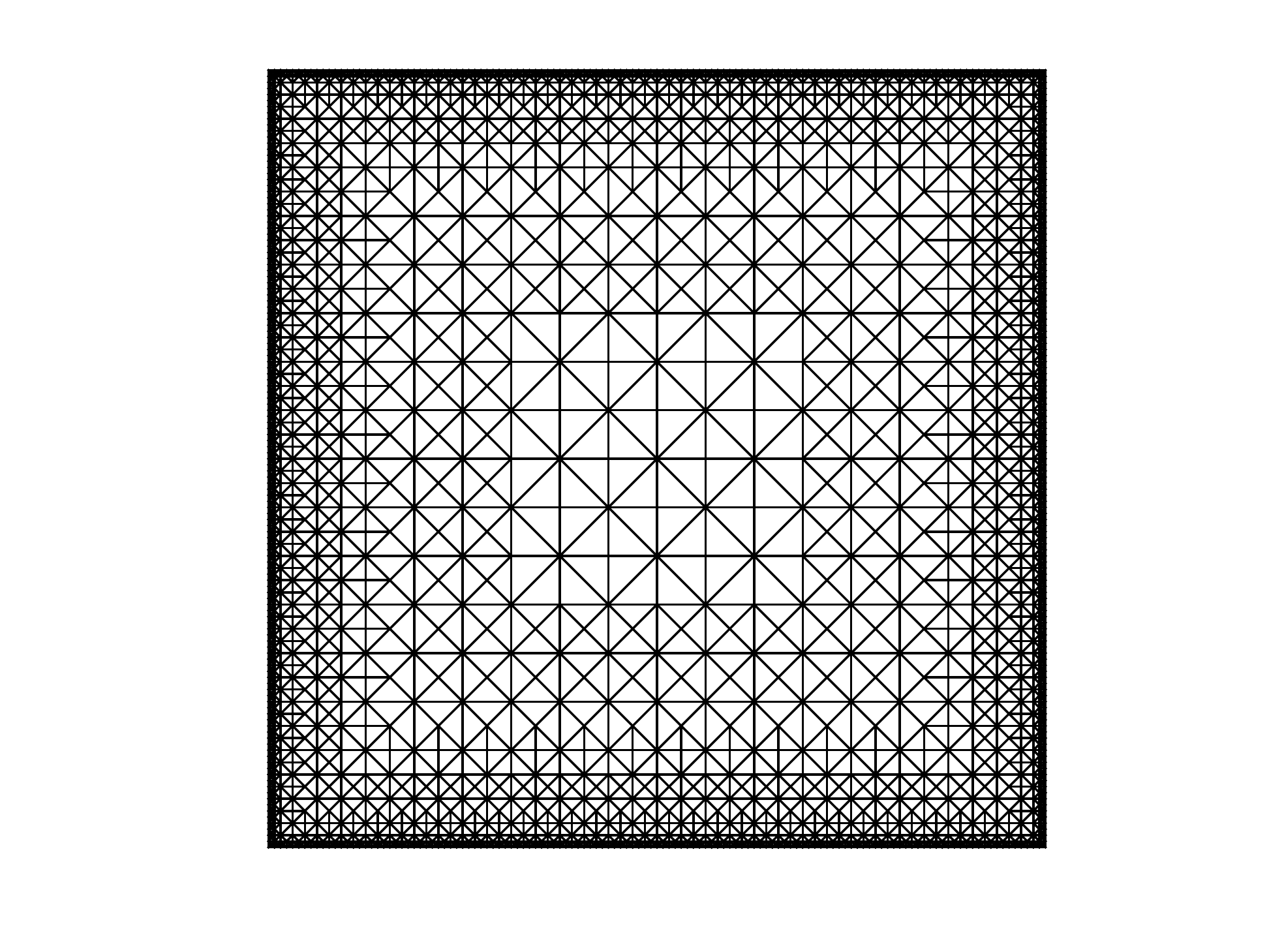} 
   \label{fig:J15}
}%
\caption{Graded bisection grids on $(-1,1)^2$, using strategy
\eqref{eq:marking} with $\theta = 4$ and $\mu = 2$.}
\label{fig:bisect}
\end{figure}

In order to obtain the graded refinement \eqref{eq:grading} when using
bisection grids, we consider the following strategy. Given an element
$\tau \in \cT$, let $x_\tau$ be its barycenter. Our strategy is based
on choosing a number $\theta > 1$ and marking those elements $\tau$
such that
\begin{equation} \label{eq:marking}
|\tau| > \theta N^{-1} \log N \cdot d(x_\tau,\partial\Omega)^{2(\mu-1)
  / \mu},
\end{equation}
where $N = \mbox{dim} \mathbb{V}(\cT)$ is the number of degrees of
freedom. We use the newest vertex marking strategy. Figure
\ref{fig:bisect} shows some graded bisection grids by using the
marking strategy \eqref{eq:marking} with $\theta = 4$, $\mu = 2$.

We document the number of iterations needed when solving the
linear systems using GS, CG and PCG over graded bisection grids, for
the same example as in the previous subsection, and with the same
stopping criterion. As shown in Table \ref{tab:FL-BPX-bisect}, the BPX
preconditioner \eqref{eq:bisec-BPX-M} performed satisfactorily in the
experiments we have carried out.

\begin{table}[!htbp]
\centering
\begin{tabular}{|m{0.3cm}|m{1.4cm}||m{1cm}|m{0.75cm}|m{0.75cm}
||m{1cm}|m{0.75cm}|m{0.75cm}||m{1cm}|m{0.75cm}|m{0.75cm}|}
  \hline
  \multirow{2}{*}{$\bar{J}$} & \multirow{2}{*}{DOFs($N$)} 
  & \multicolumn{3}{c||}{$s=0.9$}
  & \multicolumn{3}{c||}{$s=0.5$} 
  & \multicolumn{3}{c|}{$s=0.1$}\\ \cline{3-11}
  & & GS & CG  & PCG  & GS & CG & PCG & GS & CG & PCG \\ \hline
  7&  61 & 35 & 10 & {10} & 12 & 10 & {7} 
  & 9 & 13 & {8}\\ \hline
  8 &  153 & 51 & 15 &{13} & 15 & 15 & {9} 
  & 9 & 21 & {10}\\ \hline
  9 &  161 & 62 & 15 &{14} & 17  & 15 & {9} 
  & 10 & 21 & {10}\\ \hline
  10 &  369 & 93  & 20 &{16} & 21 & 20 & {11} 
  & 9 & 34 & {11}\\ \hline
  11 & 405 & 125  & 21 &{16} & 25 & 19 & {11} 
  & 9 & 31 & {12}\\ \hline
  12 & 853 & 160 & 26 &{18} & 29 & 26 & {12} 
  & 9 & 48 & {12}\\ \hline
  13 & 973 & 224 & 30 &{19} & 35 & 26 & {12} 
  & 9 & 47 & {12}\\ \hline
  14 & 1921 & 282 & 34 &{20} & 41 & 33 & {13} 
  & 9 & 72 & {12}\\ \hline
  15 & 2265 & 407 & 40 &{21} & 50 & 32 & {13} 
  & 9 & 65 & {12}\\ \hline
  16 & 4269 & 532 & 46 &{22} & 58 & 39 & {14} 
  & 9 & 97 & {13}\\ \hline
   17 & 5157 & 745 & 55 &{22} & 70 & 40 & {14} 
  & 10 & 92 & {12}\\ \hline
   18 & 9397 & 997 & 64 &{24} & 83 & 48 & {14} 
  & 9 & 135 & {13}\\ \hline
\end{tabular}
\medskip
\caption{Number of iterations: GS, CG and PCG with BPX preconditioner
  \eqref{eq:bisec-BPX-M}, $\tilde{\gamma} = \sqrt{2}/2$.}
\label{tab:FL-BPX-bisect}
\end{table}

%% end of file 

\bibliographystyle{plain}
\bibliography{bpx.bib}

%------------------------------------------------------------------------
%------------------------------------------------------------------------
\appendix

%------------------------------------------------------------------------
%------------------------------------------------------------------------

\section{Parallel Subspace Correction preconditioners: proof of Lemmas \ref{lm:PSC} and \ref{lm:auxiliary}} \label{sec:PSC-lemmas}

%------------------------------------------------------------------------
%------------------------------------------------------------------------
In this appendix we prove Lemmas \ref{lm:PSC} (identity
for PSC) and \ref{lm:auxiliary} (estimate on $\cond(BA)$) that
are instrumental for the theory developed in this paper. These results are
well known in the theory of fictitious or auxiliary spaces
\cite{nepomnyaschikh1992decomposition, griebel1995abstract,
xu1996auxiliary, xu2002method}, and implicitly given in \cite{xu1992iterative}.

\begin{proof}[Proof of Lemma \ref{lm:PSC}]
Since $R_j$ is SPD for all $0\le j\le J$, we have
$$ 
(Bv, v) = \sum_{j=0}^J (R_j Q_j v, v) = \sum_{j=0}^J (R_j Q_j v, Q_jv)
  \geq 0
$$
for all $v \in V$ and equality holds provided $Q_j v = 0$. This implies $v = 0$
because $\sum_{j=0}^J V_j = V$.

For any $v \in V$, we define $v_j^\star = R_jQ_jB^{-1}v$ and note that
$\sum_{j=0}^J v_j^\star = v$. Then
$$ 
  \sum_{j=0}^J (R_j^{-1} v_j^{\star}, v_j^\star ) =
  \sum_{j=0}^J (Q_jB^{-1}v, v_j^\star) = (B^{-1}v, v), 
$$ 
which gives a special decomposition that satisfies \eqref{eq:PSC}.  For
any other decomposition $v=\sum_{j=0}^J v_j$, we write $v_j = v_j^\star + w_j$
with $\sum_{j=0}^J w_j = 0$ and observe that
\begin{align*}
  \sum_{j=0}^J (R_j^{-1}v_j, v_j) &= \sum_{j=0}^J
  (R_j^{-1}(v_j^\star + w_j), v_j^\star + w_j) \\
  &= (B^{-1}v, v) + 2\sum_{j=0}^J (R_j^{-1}v_j^\star, w_j) +
  \sum_{j=0}^J (R_j^{-1}w_j, w_j).
\end{align*}
Since
\[
\sum_{j=0}^J (R_j^{-1} v_j^\star, w_j ) 
= \sum_{j=0}^J (B^{-1} v, w_j ) 
= (B^{-1} v, \sum_{j=0}^J w_j ) = 0,
\]
we deduce
\[
\sum_{j=0}^J (R_j^{-1}v_j, v_j) = 
 (B^{-1}v, v) + \sum_{j=0}^J (R_j^{-1}w_j, w_j) \geq  (B^{-1}v, v).
\]
This gives \eqref{eq:PSC} and concludes the proof.
\end{proof}

\begin{proof}[Proof of Lemma \ref{lm:auxiliary}]
  We note that $BA: V \to V$ is SPD with the inner product $(A \cdot, \cdot)$.
  If $(\lambda,v)$ is an eigenpair of $(BA)^{-1}$, then $B^{-1}v = \lambda Av$.
  The stable decomposition \eqref{eq:stable-decomp} thus yields
  $$ 
  \lambda_{\max}((BA)^{-1}) = \sup_{\|v\|_A = 1} (B^{-1}v, v) 
  = \sup_{\|v\|_A = 1}\inf_{\sum_{j=0}^J v_j = v}
  \|v_j\|_{R_j^{-1}}^{2} \leq c_0, 
  $$
  whence $\lambda_{\min} (BA) \ge c_0^{-1}$. On the other
  hand, the boundedness \eqref{eq:boundedness} implies
  $$ 
  \lambda_{\min}((BA)^{-1}) = \inf_{\|v\|_A = 1}(B^{-1}v, v) =
  \inf_{\|v\|_A = 1}\inf_{\sum_{j=0}^J v_j = v}
  \|v_j\|_{R_j^{-1}}^{2} \geq c_1^{-1}, 
  $$ 
  which gives $\lambda_{\max}(BA)\le c_1$. Applying the definition
  $\cond (BA) = \lambda_{\max}(BA)\lambda_{\min}(BA)^{-1}$ concludes the proof.
\end{proof}

%------------------------------------------------------------------------
%------------------------------------------------------------------------

\section{Two auxiliary results: proof of Lemmas \ref{lm:norm-equiv1}
and \ref{lm:inv_ineq}}
\label{sec:auxiliary-lemmas}

%------------------------------------------------------------------------
%------------------------------------------------------------------------

The proof of Lemma \ref{lm:norm-equiv1} involves interpolation of
weighted $L^2$ spaces \cite[Lemma 23.1]{Tartar07}. 

\begin{lemma}[interpolation of weighted $L^2$ spaces] \label{lem:tartar}
Given $w : \Omega \to (0, \infty)$ measurable, let 
\[
E(w) := \left\{ v : \Omega \to \polR \colon \int_\Omega |v|^2 w <
  \infty \right\}, \quad \| v \|_{E(w)} := \left( \int_\Omega |v|^2 w
  \right)^{\frac12}.
\]
If $w_0, w_1$ are two functions as above, then for $s \in (0,1)$ one
has
$
(E(w_0), E(w_1))_{s,2} = E(w_s) ,
$
where $w_s = w_0^{1-s} w_1^s$. Moreover, the interpolation norm
\eqref{eq:interpolation_norms} is equivalent to the $E(w_s)$ norm,
with a equivalence constant independent of $s$.
\end{lemma}
  
Upon invoking the modified $K$-functional \eqref{eq:modified-K}, Lemma
\ref{lm:norm-equiv1} ($s$-uniform interpolation) is a consequence of
interpolation theory.  See \cite{xu1989thesis} for a non-optimal
version of this result.

\begin{proof}[Proof of Lemma \ref{lm:norm-equiv1}]
We consider the spaces
\[
\begin{aligned}
X^0 = (V_{J}, \| \cdot \|_0 ) , \qquad &  Y^0 =
(V_0 \times V_1 \times \ldots \times V_{J}, \| \cdot \|_{Y^0}),  \\
X^1 = (V_{J}, | \cdot |_1 ), \qquad & Y^1 =
(V_0 \times V_1 \times \ldots \times V_{J}, \| \cdot \|_{Y^1}),
\end{aligned}
\]
where, for $\utilde v = (v_0, \ldots , v_{J}) \in V_0
\times V_1 \times \ldots \times V_{J},$
\[ \begin{aligned}
& \| \utilde{v} \|_{Y^0} := \left( \sum_{j=0}^{J} \| v_j \|_0^2 \right)^{\frac12}, 
& \| \utilde v \|_{Y^1} := \left( \sum_{j=0}^{J} \gamma^{-2j} \| v_j \|_0^2. \right)^{\frac12}.
\end{aligned} \]
Furthermore, we shall denote, for $i = 0,1$, $j = 0, \ldots, J$, $(Y_j^i, \| \cdot \|_{Y_j^i}) = (V_j, \gamma^{-ij} \| \cdot \|_0)$.

We now consider the map $Tv = (\wt{Q}_0 v , \ldots,
\wt{Q}_{J} v )$. By $L^2$-orthogonality, this map satisfies
\[
\| Tv \|_{Y^0} = \| v \|_{X^0} \quad \forall v \in X^0.
\]
The assumption guarantees that
\[
\| Tv \|_{Y^1} \lesssim \| v \|_{X^1} \quad \forall v \in X^1.
\]
Therefore, by interpolation theory, the map $T$ satisfies
\[
T : (X^0, X^1)_{s,2} \to (Y^0, Y^1)_{s,2} \quad \forall s \in (0,1),
\]
with a continuity constant independent of $s$. 
As discussed in Section \ref{sec:interpolation}, we have that $(X^0,
X^1)_{s,2} = (V_{J}, \| \cdot \|_{X^s})$, and that the
interpolation norm is equivalent to the $| \cdot |_s$ norm, with
an equivalence constant independent of $s$. 

We need to verify that the interpolation norm in $(Y^0, Y^1)_{s,2}$
coincides with the left hand side in \eqref{eq:norm-equiv1}. For that
purpose, given $\utilde w \in Y^0 + Y^1$ with $\utilde w = (w_0,
\ldots, w_{J})$, we have
\[ 
  \begin{aligned}
  K_2(t,\utilde w)^2 & = \inf_{\substack{\utilde w^0 \in X^0, \utilde
    w^1 \in X^1 \\ \utilde w = \utilde w^0 + \utilde w^1}} \| \utilde
    w^0 \|_{Y^0}^2 + t^2 \| \utilde w^1 \|_{Y^1}^2 \\
%& = \inf_{\substack{\utilde w^0 \in X^0, \utilde w^1 \in X^1 \\
%    \utilde w = \utilde w^0 + \utilde w^1}} \sum_{j=0}^{J} \| w^0_j
%    \|_0^2 + t^2 \gamma^{-2j} \| w^1_j \|_0^2 \\
& = \sum_{j=0}^{J} \inf_{\substack{w^0_j \in V_j, w^1_j \in V_j \\ w_j
    = w^0_j + w^1_j}}  \| w^0_j \|_0^2 + t^2 \gamma^{-2j} \| w^1_j
    \|_0^2 = \sum_{j=0}^{J} K_2 (t, w_j)^2.
\end{aligned} 
\]
Therefore, we can write the interpolation norm as $\| \utilde w \|_{(Y^0,
  Y^1)_{s,2}}^2 = \sum_{j=0}^{J} \| w_j \|_{(Y^0_j,Y^1_j)_{s,2}}^2$.
By Lemma \ref{lem:tartar}, we have $\| w_j \|_{(Y^0_j, Y^1_j)_{s,2}}
\simeq \gamma^{-sj} \| w_j \|_0$, with an equivalence constant independent
of $s$. Thus, we have proved the desired result \eqref{eq:norm-equiv1}.
\end{proof}
%\[ \begin{aligned}
%\| \utilde w \|_{(Y^0, Y^1)_{s,2}}^2 & = \frac{2 \sin(\pi s)}{\pi} \int_0^\infty t^{-1-2s} K_2(t,\utilde w)^2 dt \\
%& = \sum_{j=0}^{J} \frac{2 \sin(\pi s)}{\pi} \int_0^\infty t^{-1-2s}
%K_2 (t, w_j)^2 dt = \sum_{j=0}^{J} \| w_j \|_{(Y^0_j,
%Y^1_j)_{s,2}}^2.
%\end{aligned}\]

The proof of Lemma \ref{lm:inv_ineq} (local inverse inequality)
exploits a localization property of fractional Sobolev spaces (cf.  
\cite[Lemma 3.2]{xu1989thesis} and \cite{Faermann2,Faermann}) and
standard local estimates.

\begin{proof}[Proof of Lemma \ref{lm:inv_ineq}]
We distinguish between $\sigma \in [0,1]$ and $\sigma \in (1, 3/2)$.

\smallskip\noindent
{\it Step 1: $\sigma \in [0,1]$.}
We decompose the seminorm $|v|_\sigma$ locally according to
\cite{Faermann2,Faermann} for $\sigma<1$
\begin{equation*} \label{eq:patches}
|v|_\sigma^2 \le \sigma
\sum_{\tau \in \mathcal T} \left( \iint_{\tau \times
    S_\tau} \frac{|v(x) - v(y)|^2}{|x-y|^{d+2\sigma}}~{\rm d}y~{\rm d}x
  + \frac{C}{\sigma \,
  h_\tau^{2\sigma}} \, \| v \|^2_{L^2(\tau)} \right) \quad \forall v\in \mathbb{V}(\cT),
\end{equation*}
where the prefactor $\sigma$ comes from \eqref{eq:defofLaps} and
the constant $C$ depends only on the spatial dimension and
shape-regularity constant of $\mathcal{T}$. We next exploit
the local quasi-uniformity of $\mathcal{T}$ and
operator interpolation theory applied to the estimates
\[
|v|_{H^1(S_\tau)} \lesssim h_\tau^{-1} \|v\|_{L^2(S_\tau)},
\qquad
\|v\|_{L^2(S_\tau)} \lesssim \|v\|_{L^2(S_\tau)},
\]
to deduce the local inverse estimate
\[
|v|_{H^\sigma(S_\tau)} \lesssim h_\tau^{-\sigma} \|v\|_{L^2(S_\tau)}.
\]
Combining this estimate with
$|v|_{H^\sigma(S_\tau)} \lesssim |v|_{H^\sigma(S_\tau)}$, operator interpolation theory gives
\[
|v|_{H^\sigma(S_\tau)} \lesssim h_\tau^{\mu-\sigma} |v|_{H^\mu(S_\tau)},
\quad \mu \in [0, \sigma]
\]
and leads to the desired estimate \eqref{eq:inv_ineq} for $\sigma<1$. The case
$\sigma = 1$ is similar and hinges on the local inverse estimate
$|v|_{H^1(S_\tau)} \lesssim h_\tau^{\mu-1} |v|_{H^\mu(S_\tau)}$, which in turn
results from operator interpolation between the estimates
$|v|_{H^1(S_\tau)} \lesssim h_\tau^{-1} \|v\|_{L^2(S_\tau)}$ and
$|v|_{H^1(S_\tau)} \lesssim |v|_{H^1(S_\tau)}$.

\smallskip\noindent
{\it Step 2: $\sigma \in (1,3/2)$.}
Let $1\le \mu \le \sigma$ and apply Step 1 to $\nabla v$
\[
|v|_\sigma^2 = |\nabla v|_{\sigma-1}^2 \lesssim
\sum_{\tau\in\mathcal{T}} h_\tau^{2(\mu-\sigma)} |\nabla v|_{H^{\mu-1}(S_\tau)}^2
= \sum_{\tau\in\mathcal{T}} h_\tau^{2(\mu-\sigma)} |v|_{H^\mu(S_\tau)}^2
\]
If $0<\mu<1$, instead, we concatenate the preceding estimate for $\mu=1$ with
the inverse estimate $|v|_{H^1(S_\tau)} \lesssim h_\tau^{\mu-1} |v|_{H^\mu(S_\tau)}$
of Step 1. We finally observe that $v \notin H^{3/2}(\Omega)$ because it is piecewise
linear, which implies that the constant hidden in \eqref{eq:inv_ineq} blows up
as $\sigma \to 3/2$. This concludes the proof.
\end{proof}

%------------------------------------------------------------------------
\section{Generalized strengthened Cauchy-Schwarz inequality} \label{sec:SCS}
%------------------------------------------------------------------------
This appendix offers a proof of an inequality in the spirit of 
the well-known strengthened Cauchy-Schwarz inequality, that is amenable
for applications in the analysis of fractional-order problems.
The usual proof for second-order problems consists of an elementwise 
integration-by-parts argument, a local argument that quantifies the 
interaction between functions with different frequencies. This is not 
possible in the present context due to the nonlocal nature of the 
fractional norms. We resort instead to the well-known characterization 
of the fractional Sobolev space $\wt{H}^s(\Omega)$ as a Bessel potential 
space; see also \cite{xunotes,jiang2015multigrid}.

\begin{lemma}[generalized strengthened Cauchy-Schwarz inequality]
\label{lm:generalized-SCS}
Let $\sigma \in [0,3/2)$ and $k \leq \ell$. Then, given $\beta>0$ such that
$\beta \leq \sigma$ and $\beta < \frac32-\sigma$, there holds
\begin{equation*} \label{eq:generalized-SCS}
(v_k, v_\ell)_\sigma \lesssim \gamma^{\beta|\ell-k|}
\bar{h}_\ell^{-\sigma}|v_k|_\sigma \|v_\ell\|_0 
\quad \forall v_k\in \overline{V}_k, v_\ell \in \overline{V}_\ell,
\end{equation*}
where the hidden constant only blows up as $\sigma\to 3/2$.
\end{lemma}
\begin{proof}
Fix $\beta$ as in the statement of the lemma, and recall that we
denote the Fourier transform by $\mathcal{F}$.  Applying Parseval's
identity, we deduce
\begin{align*}
(v_k, v_\ell)_\sigma &= \int_{\mathbb{R}^d} |\xi|^\sigma
\mathcal{F}({v_k}) |\xi|^\sigma
\overline{\mathcal{F}({v_\ell})}\,\rm{d}\xi \\
  & = \int_{\mathbb{R}^d} |\xi|^{\sigma+\beta} \mathcal{F}({v_k})
|\xi|^{\sigma-\beta}\overline{\mathcal{F}({v_\ell})} \, {\rm{d}} \xi
  \leq  |v_k|_{\sigma+\beta} |v_\ell|_{\sigma-\beta}.
\end{align*}
The assertion follows upon applying the inverse inequality \eqref{eq:inv_ineq}
on quasi-uniform meshes and recalling that
$\bar{h}_k \simeq \gamma^{k}$:
\begin{align*}
(v_k, v_\ell)_\sigma  \lesssim \bar{h}_k^{-\beta}|v_k|_{\sigma}
\bar{h}_\ell^{-\sigma+\beta} \|v_\ell\|_0
\simeq \gamma^{\beta|\ell-k|} \bar{h}_\ell^{-\sigma} |v_k|_{\sigma}
\|v_\ell\|_0.
\end{align*}
This completes the proof.
\end{proof}

Relying on the generalized strengthened Cauchy-Schwarz inequality, 
we next offer a second proof of the boundedness property in 
Lemma \ref{lm:auxiliary} (estimate on $\cond (BA)$).

\begin{proof}[Alternative proof of Proposition \ref{P:boundedness} (boundedness)]
Combining the inverse inequality \eqref{eq:inv_ineq} with 
$\sigma = \beta = s$ and $\mu = 0$, with Lemma \ref{lm:generalized-SCS} 
(generalized strengthened Cauchy-Schwarz inequality), we obtain
\begin{equation*} %\label{eq:uniform-bd1}
\begin{aligned}
\Big|\sum_{k=0}^{\bar{J}} v_k\Big|_s^2 & \lesssim 
|v_{\bar{J}}|_s^2 + \Big|\sum_{k=0}^{\bar{J}-1} v_k\Big|_s^2 
= |v_{\bar{J}}|_s^2 + \sum_{k,\ell=0}^{\bar{J}-1} (v_k, v_\ell)_s \\ 
&\lesssim \bar{h}_{\bar{J}}^{-2s}\|v_{\bar{J}}\|_0^2 +
\sum_{k,\ell=0}^{\bar{J}-1} \gamma^{s|k-\ell|} \bar{h}_k^{-s}
\bar{h}_\ell^{-s} \|v_k\|_{0}\|v_\ell\|_0,
\end{aligned}
\end{equation*}
provided $\beta+\sigma=2s<\frac32$. We recall the elementary
inequality for $\theta < 1$,
\begin{equation}\label{eq:gamma-matrix}
\sum_{i,j=1}^n\theta^{|i-j|}x_i \, y_j\le\frac{2}{1-\theta}
\bigg(\sum_{i=1}^nx_i^2\bigg)^{\frac12}\bigg(\sum_{j=1}^n y_j^2\bigg)^{\frac12}
\quad \forall (x_i)_{i=1}^n, (y_i)_{i=1}^n \in \mathbb{R}^n.
\end{equation}
We next apply \eqref{eq:gamma-matrix} with $\theta = \gamma^s$ to obtain
\begin{equation*} %\label{eq:uniform-bd2}
\sum_{k,\ell=0}^{\bar{J} - 1} \gamma^{s|k-\ell|}
\bar{h}_k^{-s} \bar{h}_\ell^{-s}\|v_k\|_{0}\|v_\ell\|_0 \lesssim 
\frac{1}{1-\gamma^s} \sum_{k=0}^{\bar{J} - 1}
\bar{h}_k^{-2s}\|v_k\|_0^2 
\simeq
\frac{1}{1-\wt{\gamma}^s} \sum_{k=0}^{\bar{J} - 1}
\bar{h}_k^{-2s}\|v_k\|_0^2.
\end{equation*}
Combining the two preceding estimates, we arrive at
the desired bound \eqref{eq:uniform-bd}, which is \eqref{eq:boundedness}
for $\overline{B}$ in accordance with the
definitions \eqref{eq:j-smoother2} and \eqref{eq:BPX-FPDE}, provided
$s<\frac34$. If $\frac34 \le s < 1$, then our choices of $\sigma, \beta$
in Lemma \ref{lm:generalized-SCS} are restricted: we take $\sigma=s$
and $\beta=\frac12 < \frac32-\sigma$. The previous argument still works but
the prefactor in \eqref{eq:uniform-bd} becomes
$(1-\wt\gamma^\beta)^{-1}$ instead. Since there is a constant $C>0$, independent of
$s$, such that
\[
\frac{1-\wt\gamma^s}{1-\wt\gamma^\beta} \le C
\]
the expression \eqref{eq:uniform-bd} is still valid in this case. This proof is
thus complete.
\end{proof}

%------------------------------------------------------------------------
\section{Matrix representation and implementation} \label{sec:matrix-representation}
%------------------------------------------------------------------------
In this appendix we briefly discuss the
implementation of BPX preconditioners. Denoting the nodal basis
functions of $V$ by $\{\phi_i\}_{i=1}^N$, we have the following
matrices:  
\begin{itemize}
  \item Stiffness matrix ${\bf K} = (k_{ij})_{ij=1}^N \in
    \mathbb{R}^{N\times N}$, where $k_{ij} = a(\phi_j, \phi_i)$;
  \item Mass matrix ${\bf M} = (m_{ij})_{ij=1}^N \in \mathbb{R}^{N\times
    N}$, where $m_{ij} = (\phi_j, \phi_i)$;
  \item Matrix representation ${\bf A} = (a_{ij})_{ij=1}^N$ of $A$:
    $A\phi_i = \sum_{j=1}^N a_{ji} \phi_j$ or equivalently
    $$
    A[\phi_1, \ldots, \phi_N] = [\phi_1, \ldots, \phi_N]{\bf A}.
    $$ 
\end{itemize}
Recalling the definition $(A\phi_j,\phi_i) = a(\phi_j,\phi_i)$ of $A$,
we deduce ${\bf K} = {\bf M}{\bf A}$ or ${\bf A} = {\bf M}^{-1}{\bf K}$.
Next we derive the matrix presentation of the BPX
preconditioner \eqref{eq:BPX-FPDE} on quasi-uniform grids. We denote
the nodal basis functions of $\overline{V}_k$ by
$\{\phi_i^k\}_{i=1}^{N_k}$ and then have the following matrices:
\begin{itemize}
\item Matrix representation $\overline{\bf I}_k\in\mathbb{R}^{N\times N_k}$ of the
  inclusion $\overline{I}_k$, often called prolongation matrix:
    $$ 
    \overline{I}_k[\phi_1^k, \ldots, \phi_{N_k}^k] = [\phi_1, \ldots,
    \phi_N]\overline{\bf I}_k;
    $$ 
  \item Matrix representation $\overline{\bf Q}_k\in\mathbb{R}^{N_k\times N}$
    of the $L^2$-projector $\overline{Q}_k$:
    $$ 
    \overline{Q}_k[\phi_1, \ldots, \phi_{N}] = [\phi_{1}^k, \ldots,
    \phi_{N_k}^k]\overline{\bf Q}_k
    $$ 
\end{itemize}
If $\overline{\bf M}_k \in \mathbb{R}^{N_k\times N_k}$ denotes the
mass matrix on $\overline{V}_k$, the definition of $L^2$-projection yields 
$$
\begin{aligned}
  \overline{\bf I}_k^T {\bf M} & =
  \left(
  \overline{\bf I}_k^T
  \begin{bmatrix}
    \phi_1 \\ \vdots \\ \phi_{N}
  \end{bmatrix}, 
  [\phi_1, \ldots, \phi_{N}]
  \right) =
\left( 
\begin{bmatrix}
  \phi_1^k \\ \vdots \\ \phi_{N_k}^k
\end{bmatrix}, 
[\phi_1, \ldots, \phi_{N}]
\right) \\
&=
\left( 
\begin{bmatrix}
  \phi_1^k \\ \vdots \\ \phi_{N_k}^k
\end{bmatrix}, 
\overline{Q}_k [\phi_1, \ldots, \phi_{N}]
\right) = 
\left(
\begin{bmatrix}
  \phi_1^k \\ \vdots \\ \phi_{N_k}^k
\end{bmatrix}, 
[\phi_{1}^k, \ldots,
    \phi_{N_k}^k]\overline{\bf Q}_k
\right) = \overline{\bf M}_k \overline{\bf Q}_k,
\end{aligned}
$$ 
Consequently, the matrix representation $\overline{\bf B}$ of
$\overline{B}$ in \eqref{eq:BPX-FPDE} reads 
$$ 
\begin{aligned}
  \overline{\bf B} &= \overline{\bf I}_{\bar{J}} \bar{h}_{\bar{J}}^{2s} +
(1-\wt{\gamma}^s)\sum_{k=0}^{\bar{J}-1} \overline{\bf I}_k
\bar{h}_k^{2s}\overline{{\bf Q}}_k \\ 
  & = \left[ \overline{\bf I}_{\bar{J}}\bar{h}_{\bar{J}}^{2s}
  \overline{\bf M}_{\bar{J}}^{-1} +
  (1-\wt{\gamma}^s)\sum_{k=0}^{\bar{J}-1} \overline{\bf I}_k
  \bar{h}_k^{2s}{\bf M}_k^{-1}\overline{{\bf I}}_k^T \right] {\bf M}
  \simeq \overline{\bf P}{\bf M}, 
\end{aligned}
$$
where we have used the equivalence ${\bf M}_k^{-1} \simeq \bar{h}_k^{-d}{\bf I}_k$
to avoid inverting ${\bf M}_k$ and
\begin{equation} \label{eq:BPX-FPDE-M}
  \overline{{\bf P}} := \bar{h}_{\bar{J}}^{2s-d} \overline{\bf I}_{\bar{J}} +
  (1-\wt{\gamma}^s)\sum_{k=0}^{\bar{J}-1} \bar{h}_k^{2s-d}\overline{\bf
  I}_k\overline{{\bf I}}_k^T.
\end{equation}
This implies
$\cond (\overline{\bf P}{\bf K}) = \cond \big((\overline{\bf P}{\bf
M})({\bf M}^{-1}{\bf K}) \big) \simeq \cond (\overline{\bf B}{\bf A}) =
\cond (\overline{B} A) \lesssim 1$, whence $\overline{\bf P}$ is a
robust preconditioning matrix for the stiffness matrix ${\bf K}$. 

On graded bisection grids, our implementation of \eqref{eq:bisec-BPX} 
hinges on the local scaling of $V_j$
$$ 
h_{j,q} := \left(\frac{|\omega_q|}{\#\cR_q} \right)^{1/d}
\quad q \in T_j\cap \cN_j;
$$ 
note that $h_{j,q} \simeq h_{j}$ in view of shape regularity.
The robust preconditioning matrix reads
\begin{equation} \label{eq:bisec-BPX-M}
{\bf P} = \sum_{p\in \cP} h_p^{2s-d} \vI_p\vI_p^T +
(1-\wt{\gamma}^s) \sum_{j=0}^J \sum_{q \in T_j\cap\cN_j} h_{j,q}^{2s-d} \vI_{j,q}
\vI_{j,q}^T,
\end{equation} 
where $\vI_{j,q}$ is the prolongation matrix from
$\mathrm{span}\{\phi_{j,q}\}$ to $V$.

%%%%%%%%%%%%%%%%%%%%%%%%%%%%%%%%%%%%%%%%%%%%%%%%%%%%%%%%%%%%%%%%%%%%%%%%%%%%%
\end{document}